%% file: MIQMC.tex
\documentclass[11pt, oneside]{article}

\input{paper_template}

\TheTitle{A Multi-Index Quasi-Monte Carlo Algorithm for Lognormal Diffusion Problems}
\TheAuthors{P. Robbe, D. Nuyens, and S. Vandewalle}

\title{\TheTitle}
\author{
  Pieterjan Robbe\thanks{Department of Computer Science, KU Leuven - University of Leuven, Celestijnenlaan 200A, B3001 Leuven, Belgium
    (\texttt{\{\href{mailto:pieterjan.robbe@kuleuven.be}{pieterjan.robbe},
    					\href{mailto:dirk.nuyens@kuleuven.be}{dirk.nuyens},
						\href{mailto:stefan.vandewalle@kuleuven.be}{stefan.vandewalle}\}@kuleuven.be}).}
  \and
  Dirk Nuyens\footnotemark[1]
  \and
  Stefan Vandewalle\footnotemark[1]
}
\usepackage{multirow}
\begin{document}

\maketitle

% abstract
\begin{abstract}
We present a Multi-Index Quasi-Monte Carlo method for the solution of elliptic partial differential equations with random coefficients. By combining the multi-index sampling idea with randomly shifted rank-1 lattice rules, the algorithm constructs an estimator for the expected value of some functional of the solution. The efficiency of this new method is illustrated on a three-dimensional subsurface flow problem with lognormal diffusion coefficient with underlying Mat\'ern covariance function. This example is particularly challenging because of the small correlation length considered, and thus the large number of uncertainties that must be included. We show numerical evidence that it is possible to achieve a cost inversely proportional to the requested tolerance on the root-mean-square error, \review{for problems with a smoothly varying random field}.
\end{abstract}

%%%%%%%%%%%%%%%%%%%%%%%%%%%%%%%%%%%%%%%%%%%%%%%%%%%%%%%%%%%
% 1 - INTRODUCTION
%%%%%%%%%%%%%%%%%%%%%%%%%%%%%%%%%%%%%%%%%%%%%%%%%%%%%%%%%%%
\section{Introduction}
\label{sec:introduction}

In a mathematical model for a real-life process, the parameters are often unknown or subject to uncertainty. These models often show up in \emph{Uncertainty Quantification} (UQ) in engineering applications. Notable examples are partial differential equations (PDEs) with random coefficients, random initial or boundary values or an uncertain geometry. UQ aims at developing rigorous methods to characterize the impact of these uncertainties on the model outputs.

Randomized UQ methods, such as the Monte Carlo method, continue to draw a lot of attention, because they allow \sr{us} to compute statistics of the model output in a non-intrusive way. However, the classical Monte Carlo (MC) method is often viewed as impractical due to the large number of expensive realizations required. It is a notorious result that the error of the MC method converges as $\mathcal{O}(1/\sqrt{N})$, where $N$ is the number of independent realizations.

The cost of MC simulation can be reduced by lowering the required number of samples by using, e.g., \emph{variance reduction techniques}, or by switching to the Quasi-Monte Carlo (QMC) method. QMC methods first became popular in 1995, when a 360-dimensional integral was computed very efficiently by Paskov and Traub~\cite{paskov1995faster}. The key to reducing the cost of the estimator lays in choosing the realizations carefully, as opposed to the random realizations in the MC method. A recent overview paper of QMC methods for PDEs with random coefficients can be found in~\cite{kuo2016application}.

\review{In 2008, the Multilevel Monte Carlo (MLMC) method was \sr{reinvented} as a very effective variance reduction technique~\cite{giles2008multilevel,giles2015multilevel,heinrich2001multilevel,cliffe2011multilevel,barth2011multi}}. MLMC is based on a multigrid idea, by assuming that  realizations with a different accuracy are available. By estimating successive differences between these approximations, the method reduces the computational cost of the estimator compared to standard MC. A recent generalization of MLMC, called Multi-Index Monte Carlo (MIMC), was proposed in~\cite{haji2016multi}. This method extends the one-dimensional level to a multi-index, allowing \sr{us} to achieve better convergence rates compared to MLMC. \review{However, the method requires more regularity of the underlying solution compared to MLMC.} The goal of this work is to combine the MIMC method with QMC methods. In this sense, our work can be viewed as a multi-index extension of~\cite{kuo2015multilevel}, or a Quasi-Monte Carlo extension of~\cite{haji2016multi}.

The text is organized as follows. In~\secref{sec:problem_formulation}, we introduce a typical application of PDEs with random coefficients that originates from geophysics. After \sr{recalling} the Multi-Index Monte Carlo method in~\secref{sec:MIMC} and its Quasi-Monte Carlo counterpart  in~\secref{sec:MIQMC}, we investigate numerically the performance of these methods in~\secref{sec:numerical_results}. We end the discussion with some conclusions and ideas for further work.

% figure: realizations
\setlength{\figureheight}{3cm}
\setlength{\figurewidth}{0.25cm}
\begin{figure}[t]
	\centering
	\subfloat[$\{1,1,2.5\}$]{\includegraphics[width=0.28\textwidth]{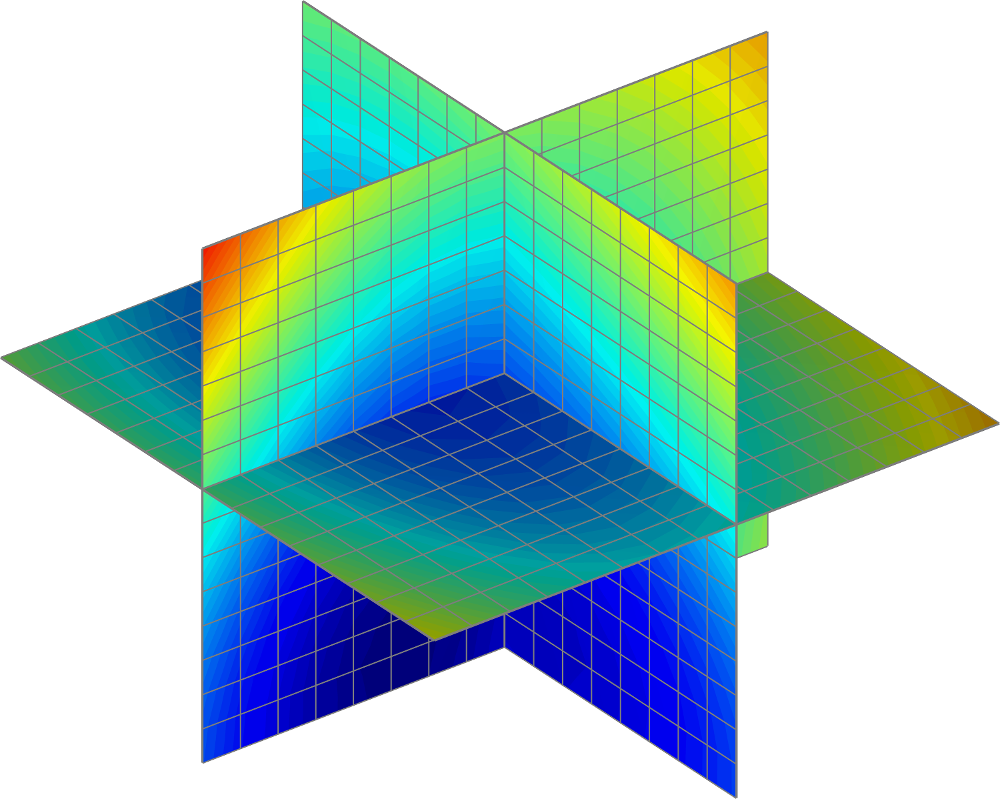}}\hfill
	\subfloat[$\{0.3,1,1\}$]{\includegraphics[width=0.28\textwidth]{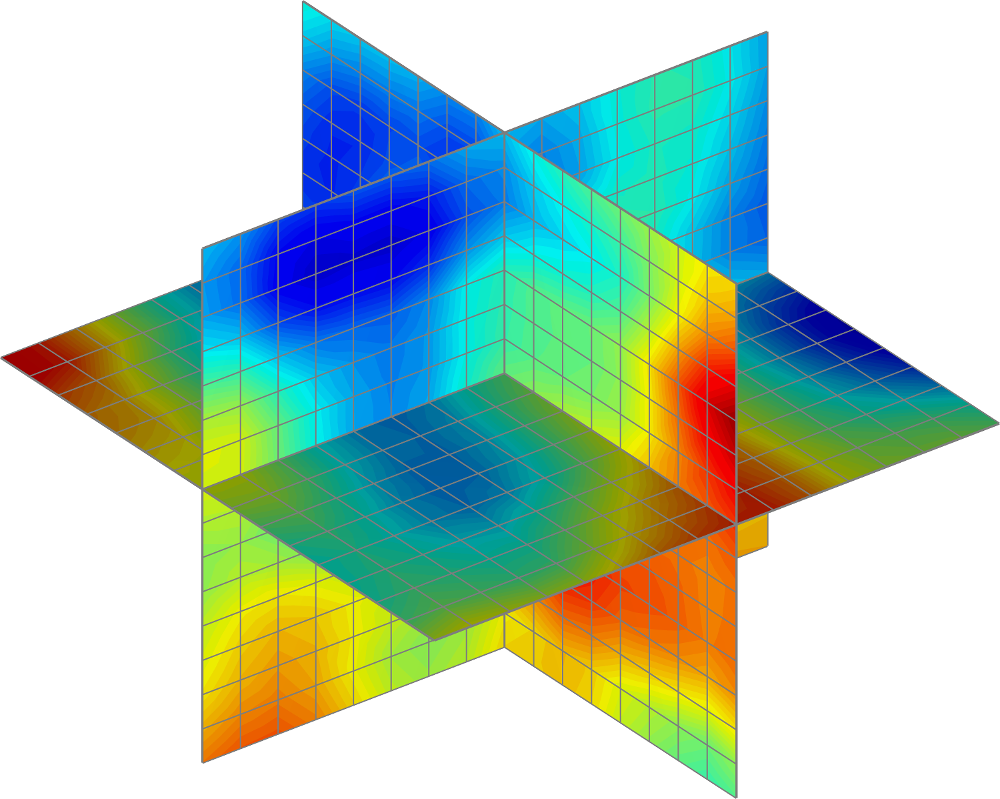}}\hfill
	\subfloat[$\{0.075,1,0.5\}$]{\includegraphics[width=0.28\textwidth]{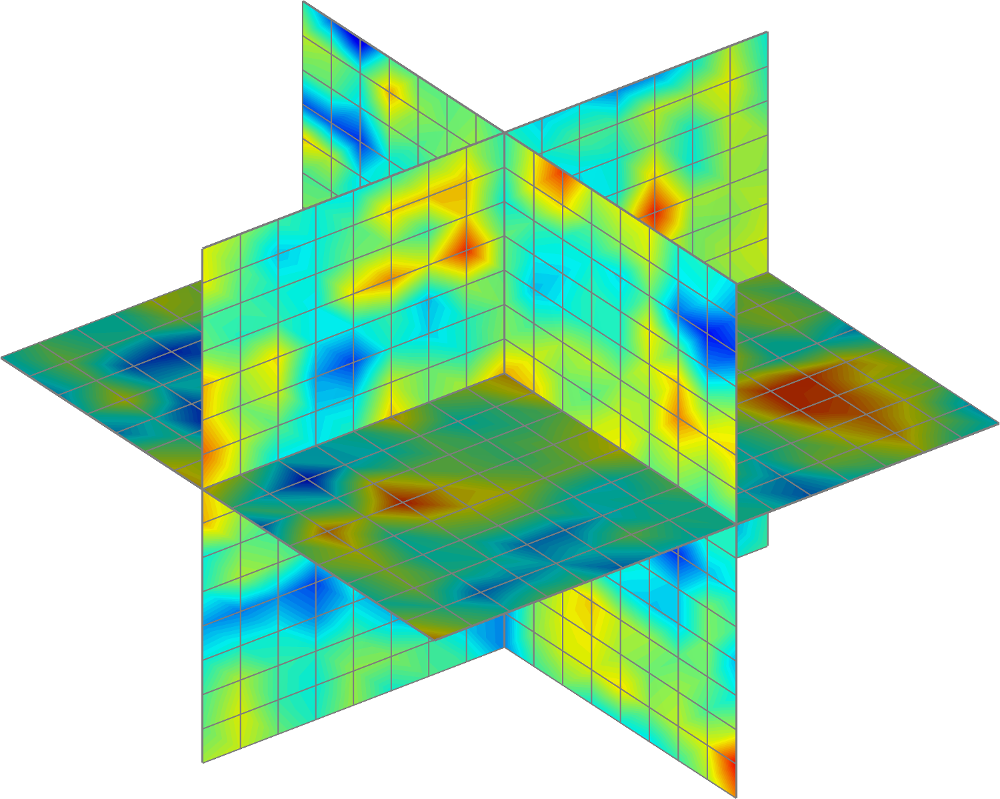}}
	\subfloat{\raisebox{-0.25cm}{%
	\tikzexternalenable
	\tikzsetnextfilename{colorbar_3}%
	\input{figures/colorbar_3.tikz}%
	\tikzexternaldisable
}}\hfill
	\caption{\review{Typical realizations of the Gaussian random field $Z(\bx,\omega)$ defined on $D=[0,1]^3$ for three sets of parameters $\{\lambda,\sigma^2,\nu\}$ in the Mat\'ern kernel. The realizations are computed using a truncated KL-expansion with 1000 terms.}}
	\label{fig:realisations}
\end{figure}

%%%%%%%%%%%%%%%%%%%%%%%%%%%%%%%%%%%%%%%%%%%%%%%%%%%%%%%%%%%
% 2 - PROBLEM FORMULATION
%%%%%%%%%%%%%%%%%%%%%%%%%%%%%%%%%%%%%%%%%%%%%%%%%%%%%%%%%%%
\section{Problem Formulation}
\label{sec:problem_formulation}

A central topic in groundwater studies is the steady-state flow through random porous media~\cite{cliffe2011multilevel}. This flow is described by Darcy's law, coupled with an incompressibility condition, leading to the \emph{parameterized PDE}
\begin{equation}\label{eq:SPDE}
	-\nabla \cdot ( k(\bx,\omega) \nabla p(\bx,\omega)) = f(\boldsymbol{x}) \quad\mathrm{for}\quad\bx \in D, \;\omega \in \Omega,
\end{equation}
where $D$ is a bounded domain in $\mathbb{R}^d$, with $d\in\{1,2,3\}$, and $\Omega$ is the sample space of a probability space $(\Omega,\mathcal{A},P)$. \review{We consider deterministic mixed Dirichlet--Neumann boundary conditions
\begin{align}
p(\bx,\omega) &= p_D(\bx) \quad &\mathrm{for\;} \bx \in \Gamma_D & \mathrm{\;and} \nonumber\\
n(\bx) \cdot (k(\bx,\omega)\nabla p(\bx,\omega)) &= p_N(\bx) \quad &\mathrm{for\;} \bx \in \Gamma_N&,\nonumber
\end{align}
with $\Gamma = \Gamma_D \cup \Gamma_N$ the boundary of the domain $D$, and $n(\bx)$ is the outward normal on the boundary $\Gamma_N$}. The diffusion coefficient $k(\bx,\omega)$ represents the permeability of the porous medium. In practice, this permeability is not known at every location $\bx$, and, in geophysics, it is commonly modeled as a \emph{random field} on $D\times\Omega$, i.e., $k:D\times\Omega\rightarrow\R:(\bx,\omega)\mapsto k(\bx,\omega)$. For a fixed sample $\omega\in\Omega$, the associated \emph{realization} of the random field is a deterministic function from $D$ to $\R$, denoted as $k(\cdot,\omega)$. Each such realization then corresponds to a deterministic version of the parameterized PDE. As a consequence, the solution of~\eqref{eq:SPDE}, the unknown hydrostatic pressure head $p(\bx,\omega)$, must itself be a random field on $D\times\Omega$. The source term $f(\bx)$ is assumed to be deterministic.

A commonly used model for the permeability $k(\bx,\omega)$ is a lognormal distribution,
\begin{equation*}
k(\bx,\omega)=\exp(Z(\bx,\omega)),
\end{equation*}
where $Z$ is an underlying Gaussian random field with given mean and covariance. The exponential ensures that the permeability remains positive throughout the domain $D$.

\review{A Gaussian random field $Z(\bx,\omega)$ is a random field where for every $M\in\N$ and $\bx_i\in D$, the vector $\boldsymbol{Z}=(Z(\boldsymbol{x}_i,\omega))_{i=1}^M$ follows a multivariate Gaussian distribution with mean $\mu_i=\mu(\bx_i)$ and covariance function
\begin{align*}
C(\bx_i,\bx_j)&\coloneqq\mathrm{cov}(Z(\bx_i,\omega),Z(\bx_j,\omega)) \\
&=\E[(Z(\bx_i,\omega)-\mu(\bx_i))(Z(\bx_j,\omega)-\mu(\bx_j))], \qquad \bx_i,\bx_j\in D.
\end{align*}
Specifically, we write $\boldsymbol{Z}\sim\mathcal{N}(\boldsymbol{\mu},\Sigma)$ with $\Sigma_{i,j}=C(\bx_i,\bx_j)$.} A Gaussian random field is fully characterized by its mean $\bmu$ and covariance function $C$. An important special case are the so-called \emph{stationary} random fields, where $\bmu$ is constant and the covariance function $C$ only depends on the difference $\bx_i-\bx_j$. Throughout this text, we will use the stationary Whittle--Mat\'ern covariance function, given by\review{
\begin{equation}\label{eq:Matern}
C(\bx_i,\bx_j) = \sigma^2\frac{1}{2^{\nu-1}\Gamma(\nu)}\left(\sqrt{2\nu}\frac{\|\bx_i-\bx_j\|_p}{\lambda}\right)^\nu K_\nu\left(\sqrt{2\nu}\frac{\|\bx_i-\bx_j\|_p}{\lambda}\right),
\end{equation}}
where $\Gamma$ is the Gamma function, $K_\nu$ is the modified Bessel function of the second kind and $\|\bx_i-\bx_j\|_p$ denotes the $\ell_p$ distance between the points $\bx_i$ and $\bx_j$. There are three parameters in this model: \sr{the correlation length $\lambda$, the (marginal) variance $\sigma^2$, and the smoothness parameter $\nu$}. By varying this set of parameters $\{\lambda,\sigma^2,\nu\}$ we can model a broad range of materials with different permeabilities, see~\figref{fig:realisations}. Note that for $\nu=1/2$, the Mat\'ern covariance reduces to the well-known exponential covariance function,
\begin{equation*}\label{eq:exponential}
	C(\bx_i,\bx_j) = \sigma^2\exp\left(-\frac{\|\bx_i-\bx_j\|_p}{\lambda}\right).
\end{equation*}

Several techniques exist to produce samples of a random field, such as the \emph{polynomial chaos expansion}~\cite{xiu2002wiener}, the \emph{circulant embedding} technique~\cite{graham2011quasi}, \review{a factorization based on \emph{H-matrices}~\cite{feischl2017fast}}, or the \emph{Karhunen--Lo\`eve} (KL) \emph{expansion}~\cite{ghanem1991stochastic}. We will focus on this last approach. The KL-expansion  
\begin{equation}\label{eq:KL}
	Z(\bx,\omega)=\mu(\bx) + \sum_{r=1}^\infty \sqrt{\theta_r}f_r(\bx) \xi_r(\omega)
\end{equation}
represents the Gaussian random field $Z(\bx,\omega)$ as a linear combination of a product of a number of eigenvalues $\theta_r$ and eigenfunctions $f_r$, with $\mathcal{N}(0,1)$-distributed random numbers $\xi_r(\omega)$ as coefficients. The eigenvalues $\theta_r$ and eigenfunctions $f_r$ are the eigenvalues and eigenfunctions of the integral operator $\mathscr{C}$ associated with the covariance function,
\begin{equation*}\label{eq:integraloperator}
(\mathscr{C}f)(\bx_i)=\int_D C(\bx_i,\bx_j)f(\bx_j) \,\mathrm{d}\bx_j,\quad\bx_i,\bx_j\in D.
\end{equation*}
For $\nu=1/2$ \review{and the $\ell_1$-distance}, analytic expressions are available for $\theta_r$ and $f_r$, see~\cite{cliffe2011multilevel}. For other $\nu$-values with $p=1$, one must solve the one-dimensional eigenvalue problem $\mathscr{C}f=\theta f$. In the numerical experiments later on, we will use the $\ell_1$-norm and discretize the operator $\mathscr{C}$ into a matrix and use its eigenvalues and \review{eigenvectors} as discrete approximations of $\theta_r$ and $f_r$. Solving this eigenvalue problem (EVP) is typically done only once, and the eigenvalues and \sr{eigenvectors} are stored for later reference. Therefore, we will ignore the cost of solving the EVP in our cost model later on.

\begin{figure}[t]
	\centering
	\tikzexternalenable
	\tikzsetnextfilename{decay}%
	\input{figures/decay.tikz}%
	\tikzexternaldisable

	\caption{\review{Decay of the three-dimensional eigenvalues with $p=1$ for three different sets of parameters $\{\lambda,\sigma^2,\nu\}$ in the Mat\'ern kernel.
	The dashed lines indicate the theoretical convergence rates.}}
	\label{fig:decayofeigenvalues}
\end{figure}
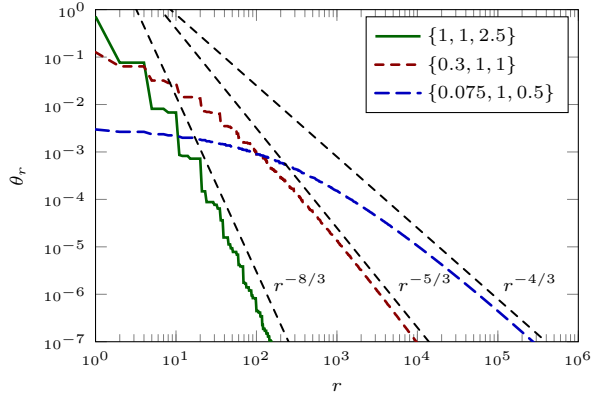

The KL-expansion is the continuous equivalent of the singular value decomposition (SVD) for matrices, and, in this sense, it is the unique expansion that minimizes the \emph{mean square error} (MSE) of the representation in $L^2(D)$ if the expansion is truncated after a finite number of terms:
\begin{equation}\label{eq:KL_trunc}
	Z_s(\bx,\omega)=\mu(\bx) + \sum_{r=1}^s \sqrt{\theta_r}f_r(\bx) \xi_r(\omega).
\end{equation}
An important question is how many terms should be retained in~\eqref{eq:KL_trunc} to accurately approximate the random field $Z(\bx,\omega)$.  If the eigenvalues $\theta_r$ decay fast, then, for large enough value of $r$, the relative contribution of $f_r$ to the sum in~\eqref{eq:KL} will be small. Hence, the faster the decay of $\theta_r$, the better an $s$-term approximation will be. Typically, the number of terms $s$ is chosen such that 95\% of the variance in the random field $Z(\bx,\omega)$ is captured by the first $s$ terms. \review{For the Mat\'ern covariance in $d$ dimensions, considered here, there is an analytic expression for the asymptotic convergence rate of the eigenvalues,
\begin{equation*}
\theta_r \sim \mathcal{O}\left(r^{-\frac{\sr{2}\nu+d}{d}}\right)\sr{,}
\end{equation*}}
\sr{see~\cite{bachmayr2017representations} or~\cite{graham2014quasi}}. It can be shown that, when $\nu=\infty$ in~\eqref{eq:Matern}, the eigenvalues decay at least exponentially, see~\cite{schwab2006karhunen}. It should be noted that, the smoother the underlying covariance function (determined by the smoothness parameter $\nu$), the faster the eigenvalues $\theta_r$ decay and thus the fewer terms are needed for an accurate representation of the random field, see~\figref{fig:decayofeigenvalues}. On the other hand, when the problem is non-smooth, a large number of initial eigenvalues have approximately the same magnitude, and a lot of terms are needed in the KL-expansion of the Gaussian random field.

In the remainder of this text, we will develop solution methods for PDEs with random coefficients such as~\eqref{eq:SPDE}, and show how to efficiently compute statistics of quantities derived from the solution of the PDE-model. For example, we will be interested in the expected value $\E[G(\omega)]= \E[\mathcal{G}(p(\bx,\omega))]$, where the \emph{quantity of interest} $G$ is a functional $\mathcal{G}$ applied to the solution $p(\bx,\omega)$. As such, we will quantify the uncertainty in the underlying PDE model.

% figure: decay of eigenvalues
\setlength{\figureheight}{6cm}
\setlength{\figurewidth}{8cm}

%%%%%%%%%%%%%%%%%%%%%%%%%%%%%%%%%%%%%%%%%%%%%%%%%%%%%%%%%%%
% 3 - MULTI-INDEX MONTE CARLO
%%%%%%%%%%%%%%%%%%%%%%%%%%%%%%%%%%%%%%%%%%%%%%%%%%%%%%%%%%%
\section{Multi-Index Monte Carlo Sampling}
\label{sec:MIMC}

In this section we review the main ideas of the Multi-Index Monte Carlo (MIMC) method, as introduced in~\cite{haji2016multi}. MIMC can be seen as an extension of the Multilevel Monte Carlo (MLMC) method \cite{giles2015multilevel} where the single scalar level is extended to a multi-index. \review{As a consequence, the hierarchy of scalar levels is extended to a larger, multi-dimensional hierarchy of indices. This allows more flexibility in choosing which grids are needed in the resulting estimator.} The method can also be seen as a combination of \emph{sparse grids} in its \emph{combination technique}-form~\cite{griebel1992combination,bungartz2004sparse} and Monte Carlo sampling.

%%%%%%%%%%%%%%%%%%%%%%%%%%%%%%%%%%%%%%%%%%%%%%%%%%%%%%%%%%%
\subsection{Derivation of the MIMC Estimator}
\label{sec:MIMC_derivation}

Consider the parameterized PDE from~\eqref{eq:SPDE}. For each realization of the random field $k(\bx,\omega)$, we must find a solution of a deterministic PDE using an appropriate numerical scheme. In our experiments later on, we will use a \review{second-order} \emph{finite volume} (FV) \emph{method}. This method is often used in the context of subsurface flow simulations, because of the mass conservation property. The FV method partitions the domain $D$ into cells with a finite volume, called \emph{control volumes}. \review{For ease of presentation, we limit ourselves to} \sr{the unit cube} $D=[0,1]^3$. Suppose we partition this domain into $m^d$ square cells. For every realization $k(\cdot,\omega)$, we compute the value of $k$ in each of the cell centers, and use the second-order FV method to find a solution $p(\cdot,\omega)$ in each of these points. From this solution, we then compute the value of a quantity of interest, such as a point evaluation on $D$ or a flux through a part of the boundary $\Gamma$. Let $G_m(\omega)$ denote the application of the quantity of interest to the discrete solution of a realization of the PDE associated with the sample $\omega$. The classical Monte Carlo method would then pick $m$ and $N$ large enough, to approximate $\E[G_m]$ by
\begin{equation}\label{eq:MC}
\mathcal{S}_N(G_m) \coloneqq \frac{1}{N}\sum_{n=0}^{N-1} G_m(\omega_n).
\end{equation}

\review{In contrast to this, the Multi-Index Monte Carlo method \cite{haji2016multi} uses discretizations with different mesh sizes for the different directions. For this, define
\begin{equation*}
\left(m_{0,i}M_i^{\ell_i}\right)_{i=1}^d\quad \mathrm{with\;integers\;} m_{0,i}>0 \mathrm{\;and\;} M_i>1,
\end{equation*}}
where $\bell\coloneqq(\ell_i)_{i=1}^d\in\N^d_0$, with $\N_0=\{0,1,2,\ldots\}$ and $d\geq1$, denotes a (multi-)index. Correspondingly, let $G_\bell(\omega)$ denote an approximation to the quantity of interest $G$ on such an $m_{0,1}M_1^{\ell_1}\times\cdots\times m_{0,d}M_d^{\ell_d}$-point mesh.

% figure: MIMC tensor grid
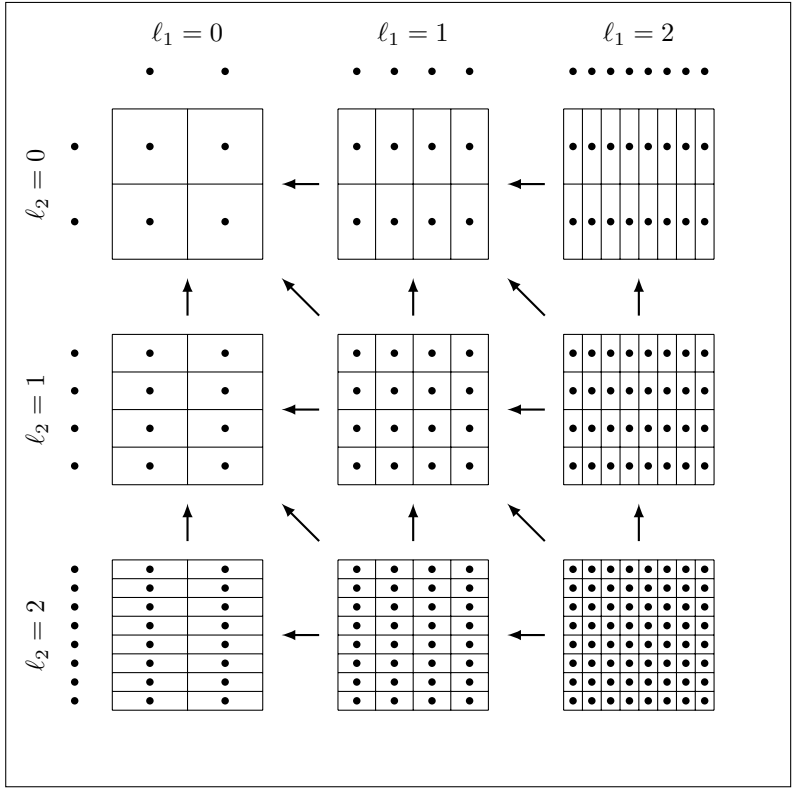
\begin{figure}[t]
	\centering
	\tikzexternalenable
	\tikzsetnextfilename{MIMCgrid}%
	\input{figures/MIMCgrid.tikz}%
	\tikzexternaldisable

	\caption{An example of multi-index grids in two dimensions. The arrows indicate which grids must be considered to compute a sample of the multi-index difference $\Delta G_\bell(\omega)$.}
	\label{fig:MIMCgrid}
\end{figure}

Instead of approximating the expected value of the quantity of interest directly on the fine mesh, the MIMC method finds approximations for the expected value of the differences $\Delta G_\bell(\omega)$ defined by
\begin{equation*}\label{eq:MIMCdiff}
\Delta G_\bell(\omega)
\coloneqq
\left(\bigotimes_{i=1}^d \Delta_i\right) G_\bell(\omega)
\end{equation*}
with
\begin{equation*}
\Delta_i G_\bell(\omega) =
\begin{cases}
G_\bell(\omega) - G_{\bell-\vec{\boldsymbol{e}}_i}(\omega) &\textrm{if}\quad\ell_i>0, \\
G_\bell(\omega) &\textrm{if}\quad\ell_i=0 .
\end{cases}
\end{equation*}
Here, $\vec{\boldsymbol{e}}_i$ denotes the unit vector in direction $i$ and 
$\otimes_{i=1}^d \Delta_i = \Delta_d \cdots \Delta_2 \Delta_1$.
 In general, taking a sample $\Delta G_\bell(\omega)$ will require a deterministic solution of the PDE at $2^d$ different grids, see~\figref{fig:MIMCgrid}. For example, to take a single sample of $\Delta G_{(1,2)}(\omega)$, we must solve the PDE four times, using four different values for the discretization parameters: $(1,2)$, $(0,2)$, $(1,1)$ and $(0,1)$. The multi-index difference is then computed as
\begin{align*}
\Delta G_{(1,2)}(\omega)&=\Delta_2(\Delta_1 G_{(1,2)}(\omega))\\
&=(G_{(1,2)}(\omega) - G_{(0,2)}(\omega)) - (G_{(1,1)}(\omega) - G_{(0,1)}(\omega))\\
&=G_{(1,2)}(\omega) - G_{(0,2)}(\omega) - G_{(1,1)}(\omega) + G_{(0,1)}(\omega).
\end{align*}
A key point is that these four solutions are based on the same realization of the random field $k(\bx,\omega)$, i.e., with the same sample $\omega$. Thus, the same random numbers $\xi_r(\omega)$, are used in its KL-expansion in~\eqref{eq:KL_trunc}. We therefore expect the quantity of interest on each of these grids to be close to each other, such \review{that the variance $\V[\Delta G_\bell(\omega)]$ of the multi-index differences will be small}. This is the rationale behind the Multi-Index Monte Carlo (MIMC) estimator
\begin{align}\label{eq:MIMC}
\mathcal{M}_L \coloneqq \sum_{\bell\in\I(L)} \mathcal{S}_{N_\bell}(\Delta G_\bell) = \sum_{\bell\in\I(L)} \frac{1}{N_\bell}\sum_{n=0}^{N_\bell-1} \Delta G_\bell(\omega_n),
\end{align}
where the set $\I(L)$ is conveniently called the \emph{index set}. The parameter $L\in\N_0$ controls the size of this index set, and $\I(L-1)\subset\I(L)$, $L=1,2\ldots$. The estimator is \review{asymptotically unbiased, i.e., the sequence $(\E[\mathcal{M}_L])_{L\geq1}$ converges to $\E[G]$}. \review{Let us denote the variance of the multi-index difference by $V_\bell\coloneqq\V[\Delta G_\bell]$. The variance of the estimator is then given by
\begin{equation}\label{eq:mi_variance}
\V[\mathcal{M}_L]=\sum_{\bell\in\I(L)}\frac{V_\bell}{N_\bell}\;\sr{=}\sum_{\bell\in\I(L)}\tilde{V}_\bell,
\end{equation}
where $\tilde{V}_\bell\coloneqq V_\bell/N_\bell$ is the contribution of index $\bell$ to the total variance of the estimator.}

\review{We assume that the index set $\calI(L)\subseteq \N_0^d$ is an \emph{admissible} (or \emph{downward closed}) {index set}, meaning that for all
\begin{align*}\label{PR_eq:DownwardClosed}
\quad\btau\leq\bell\in\calI(L)\Rightarrow\btau\in\calI(L),
\end{align*}
where $\btau\leq\bell$ means $\tau_j\leq\ell_j$ for all $j$, see~\cite{chkifa2014high}. Hence, for every index $\bell\ne(0,0,\ldots)$ in an admissible index set, all indices with smaller entries in at least one direction are also included in the set. Amongst others, this condition ensures that the index set does not contain gaps.}

Throughout this text, we will encounter two different types of index sets:
\vspace{0.5\baselineskip}
\begin{itemize}
\item Full Tensor (FT) index sets:
\end{itemize}
\begin{equation}\label{eq:FT}
\qquad\I(L) = \left\{\bell\in\N^d:\ell_i\leq L \textrm{ for all }  1\leq i\leq d\right\}\textrm{, and}
\end{equation}
\begin{itemize}
\item Total Degree (TD) index sets: 
\end{itemize}
\begin{equation}\label{eq:TD}
\qquad\I(L) = \left\{\bell\in\N^d:\sum_{i=1}^d\ell_i\leq L\right\}.
\end{equation}
The latter is inspired by the so-called Smolyak-construction in sparse grids~\cite{smolyak1960interpolation,gerstner2003dimension}. The FT index set would include all grids shown in~\figref{fig:MIMCgrid}, whereas the TD index set corresponds to the grids inside the upper left triangle. Note that the classical Multilevel Monte Carlo method only includes the grids on the main diagonal of~\figref{fig:MIMCgrid}.

In~\cite{haji2016multi}, it is shown that indices that are contained inside the weighted $d$-simplex, i.e., the weighted TD-type
\begin{equation*}
\I(L) = \left\{\bell\in\N^d:\sum_{i=1}^d\delta_i\ell_i\leq L\right\}, \quad \sum_{i=1}^d \delta_i = 1 \text{ and } 0<\delta_i\leq1,
\end{equation*}
\sr{form an optimal index set under certain conditions}.

\review{At this point, we should stress that the MIMC method is not limited to problems with a FV discretization on a \sr{unit cube}. Any discretization method on any domain $D$  that allows for a successive refinement in some direction can be used for solving the deterministic PDE underlying the parameterized model~\eqref{eq:SPDE}. Furthermore, the multi-index is not only restricted to the number of dimensions in the physical discretization, see e.g.,~\cite{robbe2017dimension}, where the multi-index controls the physical discretization as well as two KL expansions.}

The objective of estimator~\eqref{eq:MIMC} is to compute the expected value of the quantity of interest, $\E[G]$, to sufficient accuracy, for example by bounding the \emph{root mean square error} (RMSE) by a tolerance parameter $\epsilon>0$:
\begin{equation}\label{eq:RMSE}
\text{RMSE}=\sqrt{\E[(\mathcal{M}_L-\E[G])^2]}\leq\epsilon.
\end{equation}
The quantity under the square root is the MSE, which can be expanded as
\begin{equation}\label{eq:MSE}
\text{MSE}=\V[\mathcal{M}_L]+\textrm{Bias}(\mathcal{M}_L,G)^2,
\end{equation}
a classical result from statistics~\cite{wackerly2008mathematical}. The first term in~\eqref{eq:MSE} is the variance of the estimator, given by~\eqref{eq:mi_variance}, and represents the statistical part of the error. It can be reduced by taking more samples. The second term in~\eqref{eq:MSE} is the square of the bias of the estimator. It can be reduced by adding more indices to the index set. \review{If we want the $\mathrm{MSE}\leq\epsilon^2$, then it is sufficient to take $\V[\mathcal{M}_L]\leq\eta\epsilon^2$ and $\textrm{Bias}(\mathcal{M}_L,G)^2\leq(1-\eta)\epsilon^2$, with $\eta\in(0,1)$.}

Denote by $W_\bell$ the amount of work to compute a single sample of the difference $\Delta G_\bell$ of the quantity of interest at index $\bell$. The \sr{optimal} number of samples $N_\bell$ at each index $\bell$ can be computed by balancing the total amount of work 
\begin{align}\label{eq:MIMCcost}
\sr{W=\sum_{\bell\in\I(L)} N_\bell W_\bell}
\end{align}
over all indices in the index set, such that the statistical part of~\eqref{eq:MSE} is satisfied.  Then, the solution of the optimization problem
\begin{align*}
&\underset{N_\bell}{\textrm{min}} \;\sr{W} \\
&\textrm{s.t.} \sum_{\bell\in\I(L)}\frac{V_\bell}{N_\bell}\leq\review{\eta\epsilon^2}
\end{align*}
can be found as
\begin{equation}\label{eq:optimal_samples}
N_\bell = \frac{1}{\review{\eta\epsilon^2}}\sqrt{\frac{V_\bell}{W_\bell}}\sum_{\btau\in\I(L)}\sqrt{V_\btau W_\btau} \quad \textrm{for all }\bell\in\I(L) 
\end{equation}
using the method of the Lagrange multipliers. In practical computations, this value can be rounded up to the nearest largest integer $\lceil N_\bell\rceil$. \review{The variance $V_\bell$ in~\eqref{eq:optimal_samples} can be approximated by a sample variance, 
\begin{equation}\label{eq:sample_variance}
V_\bell\approx\frac{1}{N_\bell-1}\sum_{n=0}^{N_\bell-1}\left(\Delta G_\bell(\omega_n) - %\left(\frac{1}{N_\bell}\sum_{n=0}^{N_\bell-1}\Delta G_\bell(\omega_n)\right)
\mathcal{S}_{N_\bell}(\Delta G_\bell)\right)^2.
\end{equation}
The contribution of index $\bell$ to the total variance of the estimator~\eqref{eq:mi_variance} is thus approximated by
\begin{equation}\label{eq:mi_var_contr}
\tilde{V}_\bell=\frac{V_\bell}{N_\bell}\approx\frac{1}{N_\bell(N_\bell-1)}\sum_{n=0}^{N_\bell-1}\left(\Delta G_\bell(\omega_n) - %\left(\frac{1}{N_\bell}\sum_{n=0}^{N_\bell-1}\Delta G_\bell(\omega_n)\right)
\mathcal{S}_{N_\bell}(\Delta G_\bell)\right)^2.
\end{equation}}
The wall clock time can be used as a cost estimate for the true cost $W_\bell$.

%%%%%%%%%%%%%%%%%%%%%%%%%%%%%%%%%%%%%%%%%%%%%%%%%%%%%%%%%%%

\subsection{An Algorithm for MIMC Simulation}
\label{sec:MIMC_algorithm}

All elements are in place to formulate a complete algorithm for MIMC simulation (\algref{alg:MIMC}). \review{As input, the method requires a requested tolerance on the RMSE of the expected value of some quantity of interest. The outputs returned by the method are the value of the MIMC estimator and an error estimate on the computed result.} We clarify some of the essential components of the algorithm. 

The algorithm is adaptive in the index set parameter $L$. That is, we start from an index set $\{(0,\ldots,0)\}$ and add more indices to the set according to~\eqref{eq:FT} or~\eqref{eq:TD}, until the total error estimate is less than the requested accuracy $\epsilon$. At each new index, $\tilde{N}$\emph{warm-up samples} are taken to get an initial estimate for the variance contribution. Note that if this number of samples exceeds the optimal number of samples in~\eqref{eq:optimal_samples}, performance \sr{deterioration} may arise, see~\cite{pauli2015multilevel}. This often happens on the fine grids, where the required number of samples is small. We find in our numerical examples that $\tilde{N}=32$ is a good trade-off. \review{There are techniques to somewhat overcome this problem, such as regression on the variance model as suggested in~\cite{giles2015multilevel}, or \emph{continuation} Multilevel Monte Carlo~\cite{collier2014continuation}. It is the latter approach that we will use in our numerical experiments later.}

The bias is computed using the heuristic
\begin{align}\label{eq:alg_bias}
\mathrm{Bias}(\mathcal{M}_L,G)
=\left|\sum_{\bell\notin\I(L)}\E[\Delta G_\bell]\right|
\approx\hat{B}\coloneqq\left|\sum_{\bell\in\partial \I(L)}\E[\Delta G_\bell]\right|
\end{align}
where $\partial \I(L)=\I(L)\setminus\I(L-1)$ is the boundary of the index set, similar \sr{to}~\cite{gerstner2003dimension}. The mean $\E[\Delta G_\bell]$ can be approximated by a sample average. The approximation along the boundary is justified for cases where $\E[\Delta G_\bell]$ decays sufficiently fast with respect to $L$. In~\cite{haji2016multi}, for example, the analysis assumes that the decay is at least exponentially fast \review{with respect to $\bell$. We stress that~\eqref{eq:alg_bias} is a heuristic, and it might fail, even with sufficient decay.} 

The algorithm continues by adding samples at each index in the index set according to~\eqref{eq:optimal_samples}. Next, an estimate for the variance of the estimator is computed. When this estimate is larger than the allowed accuracy $\epsilon^2/2$, we double the number of samples at the index where the ratio of variance contribution \review{$\tilde{V}_\bell= V_\bell/N_\bell$} and cost $W_\bell$ is largest. \review{In the next iteration, formula~\eqref{eq:optimal_samples} is reevaluated and additional samples are taken accordingly. 

% algorithm: MIMC
\begin{algorithm}[H]
	\centering
	\caption{MIMC}
	\label{alg:MIMC}
	\begin{algorithmic}
		\REQUIRE{tolerance $\epsilon$ on RMSE}
		\STATE{$L=0$}
		\STATE{$\I(L)=\{(0,\ldots,0)\}$}
		\STATE{$\I(L-1)=\varnothing$}
		\STATE{$\mathrm{error}=\infty$}
		%\STATE{$\eta=1/2$}
		\REPEAT
			\FOR{$\bell\in\I(L)\textbackslash\I(L-1)$}
				\STATE{take $\tilde{N}$ MC warm-up samples at index $\bell$}
				\STATE{use the sample variance~\eqref{eq:sample_variance} as an estimate for the variance $V_\bell$}
				\STATE{compute the contribution $\tilde{V}_\bell$ to the total variance of the estimator using~\eqref{eq:mi_var_contr}}
			\ENDFOR
			\STATE{compute an estimate $\hat{V}$ for the variance of the estimator using~\eqref{eq:mi_variance}}
			\sr{\REPEAT
				%\STATE{compute an estimate $\hat{B}$ for the bias using~\eqref{eq:alg_bias}}
				%\IF{$\hat{B}^2<\epsilon^2/2$}
				%	\STATE{$\eta=1-\hat{B}^2/\epsilon^2$}
				%\ENDIF
				\FOR{$\bell\in\I(L)$}
					\STATE{compute the optimal number of samples $N_\bell$ at index $\bell$ using~\eqref{eq:optimal_samples}}
					\STATE{take additional MC samples at index $\bell$, to have at least $\lceil N_{\bell}\rceil$ MC samples}
					\STATE{use the sample variance~\eqref{eq:sample_variance} as an estimate for the variance $V_\bell$}
					\STATE{compute the contribution $\tilde{V}_\bell$ to the total variance of the estimator using~\eqref{eq:mi_var_contr}}
				\ENDFOR
				\STATE{compute an estimate $\hat{V}$ for the variance of the estimator using~\eqref{eq:mi_variance}}
				\IF{$\hat{V}>\eta\epsilon^2$}
				\STATE{find the index $\btau\in\I(L)$ with largest ratio $\tilde{V}_\btau/W_\btau$}
					\STATE{double the number of MC samples at index $\btau$}
					\STATE{use the sample variance~\eqref{eq:sample_variance} at index $\btau$ as an estimate for the variance $V_\btau$}
					\STATE{compute the contribution $\tilde{V}_\btau$ to the total variance of the estimator using~\eqref{eq:mi_var_contr}}
				\ENDIF
			\UNTIL{$\hat{V}\leq\epsilon^2/2$}}
			\IF{$L\geq2$}
				\STATE{compute an estimate $\hat{B}$ for the bias using~\eqref{eq:alg_bias}}
				\STATE{\sr{$\mathrm{error}=\sqrt{\hat{V}+\hat{B}^2}$}}
			\ENDIF
			\STATE{$L\coloneqq L+1$}
		\UNTIL{$\mathrm{error}<\epsilon$}
		\STATE{evaluate the MIMC estimator $\mathcal{M}_L$ using~\eqref{eq:MIMC}}
		\RETURN $\mathcal{M}_L$, error
	\end{algorithmic}
\end{algorithm}

That way, our estimator is guaranteed to have a variance smaller than or equal to a fraction \sr{$1/2$} of the MSE budget.} %In principle, this technical detail is not needed, but then the variance constraint \sr{would only be satisfied in a probabilistic manner, see~\cite{collier2014continuation}.}}

\sr{Note that \algref{alg:MIMC} is presented for a fixed $\eta=1/2$. In our implementation, we adapted the error splitting parameter when the square of the bias is smaller than $\epsilon^2/2$, where we further restrict \sr{$\eta\in[1/2,1)$}. We have that $\eta=1-\hat{B}^2/\epsilon^2$. Thus, the remaining portion of the MSE budget is used to relax the requirement on the variance of the estimator.}

% algorithm: MIQMC
\begin{algorithm}[H]
	\centering
	\caption{MIQMC}
	\label{alg:MIQMC}
	\begin{algorithmic}
		%\REQUIRE{tolerance $\epsilon$ on RMSE}
		\STATE{$L=0$}
		\STATE{$\I(L)=\{(0,\ldots,0)\}$}
		\STATE{$\I(L-1)=\varnothing$}
		\STATE{$\mathrm{error}=\infty$}
		%\STATE{$\eta=1/2$}
		\REPEAT
			\FOR{$\bell\in\I(L)\textbackslash\I(L-1)$}
				%\IF{$L\leq2$}
				\STATE{take $\tilde{N}^\star$ QMC warm-up samples at index $\bell$ for each random shift $\boldsymbol{\Xi}_{k,\bell}$} %and compute sample variance $\hat{V}_{\bell}$}
				\STATE{compute the variance contribution $\tilde{V}_\bell^\star$ of the difference $\Delta G_\bell$ using~\eqref{eq:QMC_variance_estimate}}
				%\ELSE
				%\STATE{fit $\hat{V}_{\bell}$ from the available indices in $\I(L-1)$}
				%\ENDIF
			\ENDFOR
			\STATE{compute an estimate $\hat{V}$ for the variance of the estimator using~\eqref{eq:MIQMC_variance}}
			\sr{\REPEAT
				\STATE{find the index $\btau\in\I(L)$ with largest ratio $\tilde{V}_\btau^\star/W_\btau$}
				\STATE{double the number of QMC samples at index $\btau$ for each random shift $\boldsymbol{\Xi}_{k,\btau}$}
				\STATE{compute the variance contribution $\tilde{V}_\btau^\star$ of the difference $\Delta G_\btau$ using~\eqref{eq:QMC_variance_estimate}}
				%\STATE{compute an estimate $\hat{B}$ for the bias using~\eqref{eq:alg_bias}}
				%\IF{$\hat{B}^2<\epsilon^2/2$}
				%	\STATE{$\eta=1-\hat{B}^2/\epsilon^2$}
				%\ENDIF
				\STATE{compute an estimate $\hat{V}$ for the variance of the estimator using~\eqref{eq:MIQMC_variance}}
			\UNTIL{$\hat{V}<\epsilon^2/2$}}
			\IF{$L\geq2$}
				\STATE{compute an estimate $\hat{B}$ for the bias using~\eqref{eq:alg_bias}}
				\STATE{\sr{$\mathrm{error}=\sqrt{\hat{V}+\hat{B}^2}$}}
			\ENDIF
			\STATE{$L\coloneqq L+1$}
		\UNTIL{$\mathrm{error}<\epsilon$}
		\STATE{evaluate the MIQMC estimator $\mathcal{M}^\star_L$ using~\eqref{eq:MIQMC}}
		\RETURN $\mathcal{M}^\star_L$, error
	\end{algorithmic}
\end{algorithm}

%%%%%%%%%%%%%%%%%%%%%%%%%%%%%%%%%%%%%%%%%%%%%%%%%%%%%%%%%%%
% 4 - MULTI-INDEX QUASI-MONTE CARLO
%%%%%%%%%%%%%%%%%%%%%%%%%%%%%%%%%%%%%%%%%%%%%%%%%%%%%%%%%%%
\section{Multi-Index Quasi-Monte Carlo Simulation}
\label{sec:MIQMC}

In this section, we derive the Multi-Index Quasi-Monte Carlo (MIQMC) estimator. We start with a short introduction on Quasi-Monte Carlo (QMC) methods, before combining such methods with the MIMC estimator from the previous section. Finally, we will discuss an algorithm for MIQMC simulation.

%%%%%%%%%%%%%%%%%%%%%%%%%%%%%%%%%%%%%%%%%%%%%%%%%%%%%%%%%%%
\subsection{Quasi-Monte Carlo Quadrature}
\label{sec:QMC}

The QMC method is a method to approximate high-dimensional integrals
\begin{equation*}
I_s(f) = \int_{[0,1]^s} f(\boldsymbol{y}) \, \mathrm{d} \boldsymbol{y}
\end{equation*}
over the unit cube $[0,1]^s$ by an equal-weight cubature rule
\begin{equation}\label{eq:QMC}
\mathcal{S}^\star_N(f) \coloneqq \frac{1}{N} \sum_{n=0}^{N-1} f(\boldsymbol{t}_n).
\end{equation}
\review{Note that we will use a $^\star$ to denote the QMC counterparts of the MC methods.} Formula~\eqref{eq:QMC} is seemingly identical to the Monte Carlo estimator $\mathcal{S}_N(f)$ in~\eqref{eq:MC}. However, instead of $\boldsymbol{t}_n\in [0,1]^s$ being i.i.d.\ uniform random numbers, the cubature points $\boldsymbol{t}_n$ are chosen deterministically to be better than random. \review{``Better'', in this setting, means ``more uniformly distributed'', a property that is measured by the \emph{discrepancy}~\cite{dick2013high}}. Some common techniques for generating these points are \emph{rank-1 lattice rules}~\cite{dick2013high} and \emph{digital nets}~\cite{dick2010digital}. Rather than the usual $\mathcal{O}(1/\sqrt{N})$ convergence behavior for Monte Carlo methods, QMC methods can, under certain conditions, achieve an integration error $\mathcal{O}(N^{-\alpha})$ with $\alpha>1/2$, see~\cite{kuo2005lifting,dick2013high}. In our work, we will use the \emph{rank-1 lattice rule} approach. An $N$-point \emph{rank-1 lattice rule} in $s$ dimensions is a QMC method with cubature points
\begin{equation} \label{eq:rank1latticerules}
\boldsymbol{t}_n = \left\{\frac{n\boldsymbol{z}}{N}\right\},\quad n=0,\ldots,N-1,
\end{equation}
where $\boldsymbol{z}\in\mathbb{Z}^s$ is an $s$-dimensional \emph{generating vector}, and $\{\cdot\}$ denotes the fractional part, i.e., $\{\bx\}=\bx-\lfloor\bx\rfloor$.

%The classical Monte Carlo method comes with a probabilistic error bound of the form $\sigma(f)/\sqrt{N}$, where $\sigma^2(f)\coloneqq I_s(f^2)-(I_s(f))^2$ is the variance of $f$. This error estimate comes practically for free, since $\sigma(f)$ can be estimated using the same samples used to approximate $I_s(f)$. 
Unfortunately, QMC methods do not provide \sr{an error bound derived from~\eqref{eq:mi_var_contr}}, since the points are chosen deterministically. However, this feature can be recovered by using \emph{random shifts}: each point in the lattice rule is shifted by a vector \review{$\boldsymbol{\Xi}\in[0,1]^s$}:
\begin{align*}
\boldsymbol{t}_n' = \left\{\frac{n\boldsymbol{z}}{N}+\boldsymbol{\Xi}\right\}=\{\boldsymbol{t}_n+\boldsymbol{\Xi}\},\quad n=0,\ldots,N-1\sr{.}
\end{align*}
We will denote the corresponding randomly shifted lattice rule as $\mathcal{S}^\star_N(f;\boldsymbol{\Xi})$. A probabilistic error estimate for the QMC method can be obtained by choosing $K$ i.i.d.\ shifts \review{$\boldsymbol{\Xi}_0,\ldots,$ $\boldsymbol{\Xi}_{K-1}$}. The approximation for the integral now becomes
\begin{align*}
\mathcal{S}^\star_{N,K}(f) &\coloneqq \frac{1}{K} \sum_{k=0}^{K-1}  \mathcal{S}^\star_N(f;\boldsymbol{\Xi}_k) \\
&= \frac{1}{K} \sum_{k=0}^{K-1} \frac{1}{N} \sum_{n=0}^{N-1} f(\{\boldsymbol{t}_n+\boldsymbol{\Xi}_k\}).
\end{align*}
\review{Since the $\mathcal{S}^\star_N(f;\boldsymbol{\Xi}_0),\ldots,\mathcal{S}^\star_N(f;\boldsymbol{\Xi}_{K-1})$ are i.i.d.\ random variables, the (sample) variance of $\mathcal{S}^\star_{N,K}(f)$, 
\begin{align}\label{eq:QMC_variance_estimate}
\V[\mathcal{S}^\star_{N,K}(f)]
\approx \frac{1}{K-1}\sum_{k=0}^{K-1}\left(\mathcal{S}^\star_N(f;\boldsymbol{\Xi}_k)-
%\left(\frac{1}{K}\sum_{k=0}^{K-1}\mathcal{S}^\star_N(f;\boldsymbol{\Xi}_k)\right)
\mathcal{S}^\star_{N,K}(f)\right)^2,
\end{align}
can be used to construct a confidence interval for $\mathcal{S}_{N,K}^\star(f)$ in the usual way, see~\cite{dick2013high}.}

The integral we consider here is the expectation of the quantity of interest, $\E[G]$.
Since the lognormal random field associated with~\eqref{eq:SPDE} is represented by an infinite number of $\mathcal{N}(0,1)$-distributed random numbers in the KL-expansion, see~\eqref{eq:KL}, we actually have to consider an integral over $\R^\infty$:
\begin{align*}
\E[G] = \E[\mathcal{G}(p(\bx,\omega))]
&=
\int_{\R^\infty} \mathcal{G}(p(\bx,\xi_1, \xi_2 \ldots)) \, \mathrm{d}\Phi(\boldsymbol{\xi})
\\
&=
\review{\int_{[0,1]^\infty} \mathcal{G}(p(\bx,\Phi^{-1}(y_1),\Phi^{-1}(y_2) \ldots)) \, \mathrm{d}\boldsymbol{y}}
\\
&\approx
\int_{[0,1]^s} \mathcal{G}(p(\bx,\Phi^{-1}(y_1),\ldots,\Phi^{-1}(y_s), 0, \ldots)) \, \mathrm{d}\boldsymbol{y} 
\\
&\approx 
\mathcal{S}^\star_N(G\circ\Phi^{-1}) = \frac{1}{N} \sum_{n=0}^{N-1} G(\Phi^{-1}(\boldsymbol{t}_n)),
\end{align*}
where $\Phi$ and $\Phi^{-1}$ are the cumulative normal density and its inverse respectively. We apply this change of variables component-wise, i.e.,
\begin{equation*}
\boldsymbol{\xi} = \Phi^{-1}(\by) = (\Phi^{-1}(y_1),\Phi^{-1}(y_2),\ldots) \in\R^\N \text{ and } \by\in(0,1)^\N.
\end{equation*}
The setting of approximating the expected value by applying a linear functional to the solution of the lognormal diffusion problem under consideration has been analyzed using randomly shifted lattice rules for a single level of discretization and for the multilevel algorithm, see, e.g.,~\cite{nichols2014fast,graham2014quasi,kuo2015multilevel,kuo2016application}.
In such a case it can be shown that the integrand belongs to a certain weighted Sobolev space with so-called \review{\emph{product and order dependent} (POD for short)} weights, where the weights denote the importance of different sets of variables.
A generating vector $\boldsymbol{z}$ for the lattice rule can then be constructed using a \emph{component-by-component} (CBC) algorithm with cost $\mathcal{O}(sN\log{N}+s^2N)$, see~\cite{nichols2014fast,kuo2016application,nuyens2006fast} for details.
Software accompanying \cite{kuo2016application} for constructing such rules is available on the internet \cite{QMC4PDE}.
The convergence \sr{rate depends} on the decay of the eigenvalues, but is limited to $O(N^{-1})$ because of the use of randomly shifted lattice rules on a non-periodic smooth function.
The convergence is however independent of the truncation dimension $s$ of the random field due to the \review{POD-weighted} Sobolev space.

%%%%%%%%%%%%%%%%%%%%%%%%%%%%%%%%%%%%%%%%%%%%%%%%%%%%%%%%%%%
\subsection{The MIQMC Estimator}
\label{sec:MIQMC_derivation}

In the remainder of this section, we will derive the Multi-Index Quasi-Monte Carlo (MIQMC) estimator. The idea of the MIQMC method is to replace the simple MC estimator for the differences $\Delta G_\bell$ in~\eqref{eq:MIMC} by the QMC method from~\secref{sec:QMC}. Due to the bias constraint we would like to satisfy in our algorithm, we need the estimator for the differences to be unbiased. This is satisfied for the randomly shifted rank-1 lattice rules presented above. The MIQMC estimator can be expressed as
\begin{equation}
\mathcal{M}^\star_L \coloneqq \sum_{\bell\in\I(L)} \mathcal{S}^\star_{N_\bell,K}(\Delta G_\bell).
\end{equation}
Fully expanded, the MIQMC estimator for $\E[G]$ based on rank-1 lattice rules reads
\begin{equation}\label{eq:MIQMC}
\mathcal{M}^\star_L=\sum_{\bell\in\I(L)}\frac{1}{K}\sum_{k=0}^{K-1} \frac{1}{N_\bell} \sum_{n=0}^{N_\bell-1} \Delta G_\bell(\Phi^{-1}(\{\bt_n+\boldsymbol{\Xi}_{k,\bell}\})),
\end{equation}
with $\Phi^{-1}$ the inverse cumulative normal. Note that we have now written an explicit dependence of the differences $\Delta G_\bell$ on the vector $\boldsymbol{\xi}=\Phi^{-1}(\{\bt_n+\boldsymbol{\Xi}_{k,\bell}\})\in\R^s$ in the KL-expansion~\eqref{eq:KL_trunc}. The MIQMC estimator is still an \review{asymptotically unbiased} estimator, and its variance is given by
\begin{align}\label{eq:MIQMC_variance}
\mathbb{V}[\mathcal{M}^\star_L]
&= \mathbb{V}\left[\sum_{\bell\in\I(L)}\mathcal{S}^\star_{N_\bell,K}(\Delta G_\bell)\right] \nonumber\\
&= \sum_{\bell\in\I(L)} \mathbb{V}\left[\mathcal{S}^\star_{N_\bell,K}(\Delta G_\bell)\right],
\end{align}
because of the i.i.d.\ random shifts $\boldsymbol{\Xi}_{k,\bell}$. \sr{The total work of the estimator is}
\begin{equation}\label{eq:MIQMCcost}
\sr{W^\star=K\sum_{\bell\in\I(L)}N_\bell W_\bell.}
\end{equation}

%%%%%%%%%%%%%%%%%%%%%%%%%%%%%%%%%%%%%%%%%%%%%%%%%%%%%%%%%%%
\subsection{An Algorithm for MIQMC Simulation}
\label{sec:MIQMC_algorithm}

We present an algorithm for MIQMC simulation in~\algref{alg:MIQMC}. Some of the remarks given in~\secref{sec:MIMC_algorithm} also apply here.

Contrary to~\algref{alg:MIMC}, there is no analytic expression for the required number of samples at each index, similar to~\eqref{eq:optimal_samples}. Instead, we will base our method on the simple yet effective algorithm given in~\cite{giles2009multilevel}: starting from an initial number of samples $\tilde{N}^\star$, we double the number of samples at the index with the largest ratio of variance contribution and cost. \review{The way this variance contribution is estimated, is the main difference with the MIMC algorithm in~\algref{alg:MIMC}.} \review{For MIQMC, the contribution $\tilde{V}_\bell^\star\coloneqq\V[\mathcal{S}^\star_{N,K}(\Delta G_\bell)]$ to the variance of the estimator, $\mathbb{V}[\mathcal{M}^\star_L]$, is computed by~\eqref{eq:QMC_variance_estimate}. This requires $K$ independent random shift $\boldsymbol{\Xi}_{k,\bell}$ at each index $\bell$, where all shifts are mutually independent.} The number of shifts $K$ needs to be chosen carefully. If $K$ is too small, the variance estimation can be poor and the algorithm may terminate too early. If $K$ is too large, it may kill the performance of the MIQMC estimator. Furthermore, the choice of $K$ also influences the choice of the number of warm-up samples $\tilde{N}^\star$. We numerically found that in our examples, presented below, any $8<K<32$ is acceptable, and in our experiments we will choose $K=16$ with $\tilde{N}^\star=4$.

%%%%%%%%%%%%%%%%%%%%%%%%%%%%%%%%%%%%%%%%%%%%%%%%%%%%%%%%%%%
% 4 - NUMERICAL RESULTS
%%%%%%%%%%%%%%%%%%%%%%%%%%%%%%%%%%%%%%%%%%%%%%%%%%%%%%%%%%%
\section{Numerical Results}
\label{sec:numerical_results}

\review{We investigate the performance of our MIQMC algorithm on the 3D flow problem introduced in~\secref{sec:problem_formulation}. We consider three different sets of parameters for the covariance function of the underlying Gaussian random field, with various degree of smoothness, and two different quantities of interest. We compare with standard MIMC simulation and the multilevel counterparts: MLMC and MLQMC~\cite{robbe2016thesis}. We show numerically that, for certain choices for the parameters in the covariance function and certain choices for the quantity of interest, such that the integrand is smooth, it is possible to obtain an estimator with a cost inversely proportional to the requested tolerance $\epsilon$ on the RMSE, which is the best possible result for randomly shifted \sr{lattice} rules in this setting~\cite{kuo2016application}.}

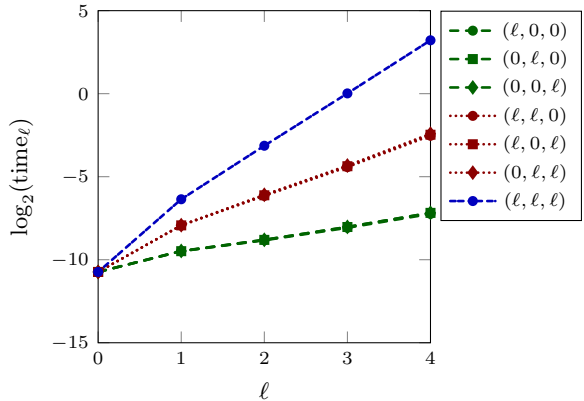
\begin{figure}[t]
\centering
	\tikzexternalenable
	\tikzsetnextfilename{G1F1_analyse_W}%
	\input{figures/G1F1_analyse_W.tikz}%
	\tikzexternaldisable

\caption{\review{Average run time to compute one realization of the multi-index difference $\Delta G_\bell$ of the quantity of interest \texttt{G1} in the three-dimensional flow problem. The notation $(\ell,0,0)$ means that we investigate how the run time behaves when only the first dimension is refined, and similar for all the other (mixed) dimensions.}}
\label{fig:analyse_cost}
\end{figure}

\review{We consider a domain $D=[0,1]^3$ and impose a grid hierarchy as explained in~\secref{sec:MIMC}. We choose $m_{i,0}=4$ and $M_i=2$ for all $1 \leq i \leq 3$. Each realization of the PDE is discretized using a cell-centered FV approach, and the resulting sparse system is solved using a preconditioned conjugate gradient method with an algebraic multigrid preconditioner~\cite{boyle2009hsl}. The average running time to compute a realization of the multi-index difference $\Delta G_\bell$ at each index $\bell$ shows an isotropic structure, i.e., the rates are the same in every dimension, as can be deduced from~\figref{fig:analyse_cost}. Also, the mixed dimension rates are the products of the respective single-dimension rates. All simulations are performed on a 2.8GHz Ivy Bridge processor with 64GB of RAM. For the implementation of MLQMC and MIQMC based on rank-1 lattice rules, we pick a standard generating vector $\boldsymbol{z}$ from~\cite{nuyens2010magic}, and choose $K=16$ random shifts.}

\setlength{\figureheight}{6cm}
\setlength{\figurewidth}{6cm}

Three different sets of parameter values for the Mat\'ern covariance function are provided, denoted as \texttt{F1} to \texttt{F3} , see~\tabref{tab:cov_parameters}. As the correlation length and smoothness decrease, we require more terms in the KL expansion of the underlying Gaussian random field. We used the criterion $\theta_s/\sr{\theta_1}\leq10^{-3}$ to determine the number of terms $s$, also listed in the table. The generated random fields have a zero-mean ($\mu(\bx)=0$), and we choose $p=1$ for the $\ell_p$-norm for all sets of parameters. Hence, we can use the analytic expressions for eigenvalues and eigenfunctions for \texttt{F3}.

\review{As a first example, consider the parameterized PDE~\eqref{eq:SPDE} with only Dirichlet boundary conditions, i.e., $p_D(\bx)=0$ and $\Gamma_D=\Gamma$. The quantity of interest is a point evaluation of the pressure at the middle of the domain, $\bx=(0.5,0.5,0.5)$. This problem will be denoted as \texttt{G1}.}

\begin{table}
\caption{\review{The three different sets of parameters used in the Mat\'ern covariance function (compare with the realizations in \figref{fig:realisations} and the eigenvalue decay in \figref{fig:decayofeigenvalues}).}}
\label{tab:cov_parameters}
\centering\review{
\begin{tabular}{r{3.8cm}p{2.5cm}p{2.5cm}p{2.5cm}}\toprule
& \texttt{F1} & \texttt{F2} & \texttt{F3} \\ \midrule
correlation length $\lambda$ & 1 & 0.3 & 0.075 \\
variance $\sigma^2$ & 1 & 1 & 1 \\
smoothness $\nu$ & 2.5 & 1 & 0.5 \\
number of KL terms $s$ & 12 & 201 & 3500 \\ \bottomrule
\end{tabular}}
\end{table}

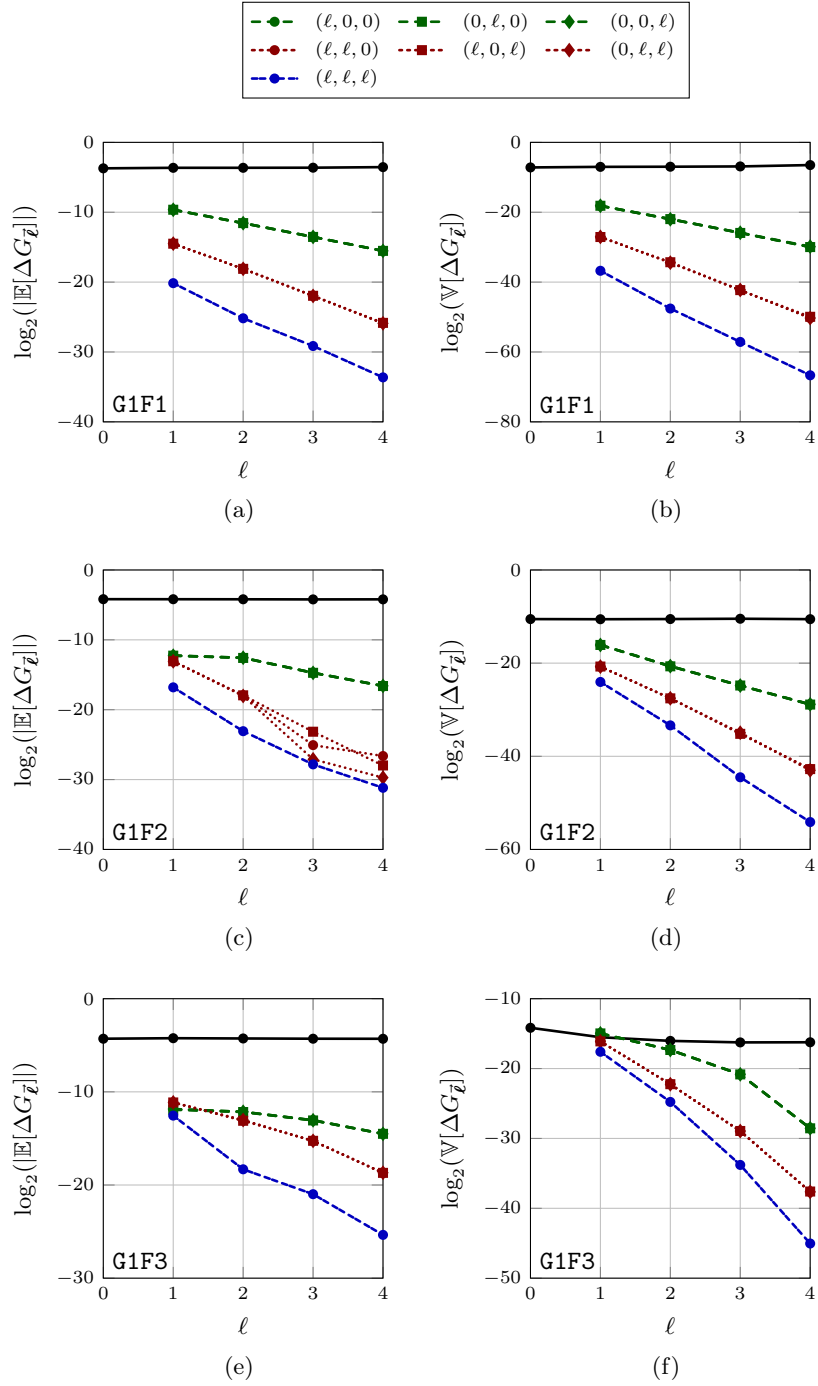
\begin{figure}[p]
\centering
\setlength{\figureheight}{2cm}
\setlength{\figurewidth}{10cm}
\hspace{1.25cm}%
	\tikzexternalenable
	\tikzsetnextfilename{legend_analyse}%
	\input{figures/legend_analyse.tikz}%
	\tikzexternaldisable
\\
\setlength{\figureheight}{5.3cm}
\setlength{\figurewidth}{5.3cm}
\subfloat[][{\makebox[-1cm]{}}]{%
	\tikzexternalenable
	\tikzsetnextfilename{G1F1_analyse_E}%
	\input{figures/G1F1_analyse_E.tikz}%
	\tikzexternaldisable
}\hspace{0.2cm}
\subfloat[][{\makebox[-1cm]{}}]{%
	\tikzexternalenable
	\tikzsetnextfilename{G1F1_analyse_V}%
	\input{figures/G1F1_analyse_V.tikz}%
	\tikzexternaldisable
}\\
\subfloat[][{\makebox[-1cm]{}}]{%
	\tikzexternalenable
	\tikzsetnextfilename{G1F2_analyse_E}%
	\input{figures/G1F2_analyse_E.tikz}%
	\tikzexternaldisable
}\hspace{0.2cm}
\subfloat[][{\makebox[-1cm]{}}]{%
	\tikzexternalenable
	\tikzsetnextfilename{G1F2_analyse_V}%
	\input{figures/G1F2_analyse_V.tikz}%
	\tikzexternaldisable
}\\
\subfloat[][{\makebox[-1cm]{}}]{%
	\tikzexternalenable
	\tikzsetnextfilename{G1F3_analyse_E}%
	\input{figures/G1F3_analyse_E.tikz}%
	\tikzexternaldisable
}\hspace{0.2cm}
\subfloat[][{\makebox[-1cm]{}}]{%
	\tikzexternalenable
	\tikzsetnextfilename{G1F3_analyse_V}%
	\input{figures/G1F3_analyse_V.tikz}%
	\tikzexternaldisable
}
\caption{\review{Behavior of the estimated mean and variance of the multi-index differences $\Delta G_\bell$ for the \sr{first} quantity of interest \sr{\texttt{G1}}. The notation $(\ell,0,0)$ means that we investigate how these quantities behave when only the first dimension is refined, and similar for all the other (mixed) dimensions. For reference, the full black line corresponds to the approximation of the quantity of interest $G_\bell$ at $(\ell,\ell,\ell)$.}}
\label{fig:analyse_G1}
\end{figure}

\review{We analyze the behavior of the mean $|\E[\Delta G_\bell]|$ and variance $\V[\Delta G_\bell]$ of the multi-index differences in all (mixed) directions of the problem. We clearly see isotropy and a product structure from~\figref{fig:analyse_G1}.}

\begin{figure}[p]
\centering
\setlength{\figureheight}{2cm}
\setlength{\figurewidth}{10cm}
	\tikzexternalenable
	\tikzsetnextfilename{legend_times_G1}%
	\input{figures/legend_times_G1.tikz}%
	\tikzexternaldisable
\\
\setlength{\figureheight}{6cm}
\setlength{\figurewidth}{6cm}
\subfloat[][{\makebox[0cm]{}}]{%
	\tikzexternalenable
	\tikzsetnextfilename{cost_G1F1}%
	\input{figures/cost_G1F1.tikz}%
	\tikzexternaldisable
}\hspace{2.0cm}
\subfloat[][{\makebox[0cm]{}}]{%
	\tikzexternalenable
	\tikzsetnextfilename{times_G1F1}%
	\input{figures/times_G1F1.tikz}%
	\tikzexternaldisable
}\\
\subfloat[][{\makebox[0cm]{}}]{%
	\tikzexternalenable
	\tikzsetnextfilename{cost_G1F2}%
	\input{figures/cost_G1F2.tikz}%
	\tikzexternaldisable
}\hspace{2.0cm}
\subfloat[][{\makebox[0cm]{}}]{%
	\tikzexternalenable
	\tikzsetnextfilename{times_G1F2}%
	\input{figures/times_G1F2.tikz}%
	\tikzexternaldisable
}\\
\subfloat[][{\makebox[0cm]{}}]{%
	\tikzexternalenable
	\tikzsetnextfilename{cost_G1F3}%
	\input{figures/cost_G1F3.tikz}%
	\tikzexternaldisable
}\hspace{2.0cm}
\subfloat[][{\makebox[0cm]{}}]{%
	\tikzexternalenable
	\tikzsetnextfilename{times_G1F3}%
	\input{figures/times_G1F3.tikz}%
	\tikzexternaldisable
}\\
\caption{\review{Performance comparison of Multi-Index (both FT and TD) and Multilevel (Quasi) Monte Carlo under different test conditions for \texttt{G1}.}}
\label{fig:G1}
\end{figure}
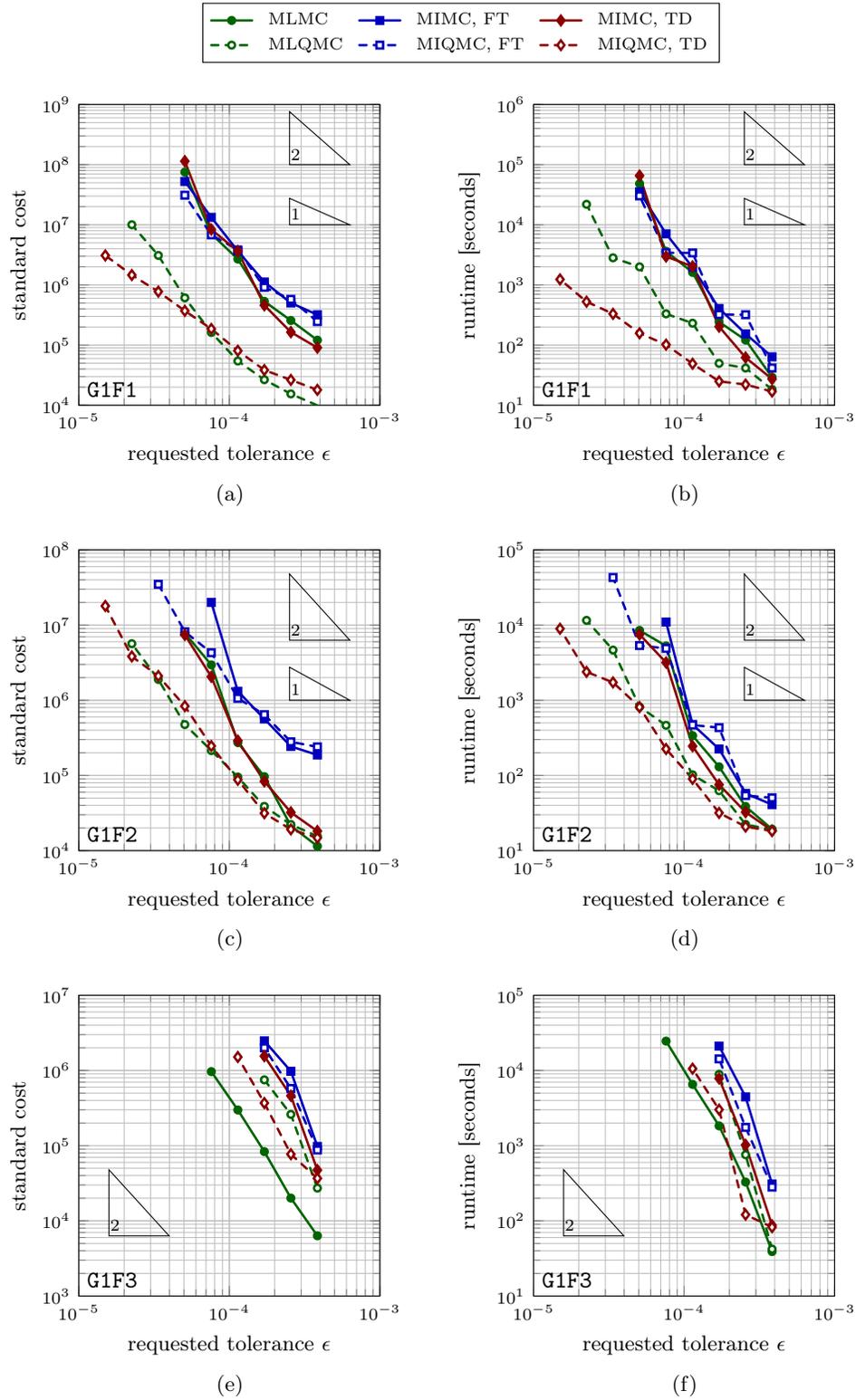

\review{Next, we compare our MIQMC estimator with both the MLQMC estimator from~\cite{robbe2016thesis,kuo2015multilevel}, and the variants based on plain Monte Carlo sampling: MIMC (with both FT and TD index sets) and MLMC. Note that for the multilevel methods, we use the implementation of the multi-index method with a single index representing the refinement in all levels at the same time. We measure the total simulation time\sr{, and the total amount of work using~\eqref{eq:MIMCcost} and~\eqref{eq:MIQMCcost},} where $W_\bell$ is computed using a regression on the actual run times from~\figref{fig:analyse_cost}. We ran the six different algorithms for a sequence of decreasing tolerances $\epsilon$ and present the results in~\figref{fig:G1}. For the smooth field \texttt{F1}, we clearly see the benefit of both QMC methods. Our MIQMC algorithm with TD index sets reaches an accuracy $\epsilon = \mathcal{O}(1/\mathrm{time}^r)$ with $r=0.92$. This is nearly optimal, since we are working with rank-1 lattice rules. Thus, the best rate we hope to achieve is $r=1$, i.e., a cost inversely proportional to the desired accuracy. Compare this to the methods based on MC-sampling, that have $r=0.5$. Furthermore, the MIMC methods that use FT index sets have a suboptimal performance. This was already observed in~\cite{haji2016multi}, and there seems to be no improvement when switching to QMC-sampling.}

\review{When the smoothness of the Mat\'ern kernel of the underlying Gaussian field decreases, the achieved rate $r$ also decreases. For \texttt{F2}, we find numerically that $r=0.71$ for MIQMC with TD index sets. This can also be seen in~\figref{fig:G1}, where for the very rough field \texttt{F3}, the benefits of the QMC method have disappeared. All methods have the same asymptotic convergence rate $r=0.5$, and the classical MLMC is actually the best method.}

\review{As a second example, consider the parameterized PDE~\eqref{eq:SPDE} with \emph{flow cell} geometry, i.e.,
\begin{align*}
\begin{cases}
p((x_1,x_2,x_3),\cdot)=1 &\mathrm{\;on\;} \left.\Gamma \right|_{x_1=0}, \\
p((x_1,x_2,x_3),\cdot)=0 &\mathrm{\;on\;} \left.\Gamma \right|_{x_1=1} \quad\mathrm{ and} \\
-k \nabla p\cdot \boldsymbol{n}=0&\mathrm{\;elsewhere.}
\end{cases}
\end{align*}
The quantity of interest we consider here is the effective permeability through the side of the domain at $x_1=1$, i.e.,
\begin{equation}\label{eq:Q2}
	\mathcal{G}=-\int_0^1\int_0^1k\left.\frac{\partial p}{\partial x_1}\right|_{x_1=1}\mathrm{d}x_2\, \mathrm{d}x_3,
\end{equation}
see~\cite{graham2011quasi,cliffe2011multilevel}.
We approximate the derivative in~\eqref{eq:Q2} by a second-order finite difference, and the integral by the two-dimensional trapezoidal rule.}

\begin{figure}[p]
\centering
\setlength{\figureheight}{2cm}
\setlength{\figurewidth}{10cm}
\hspace{1.25cm}%
	\tikzexternalenable
	\tikzsetnextfilename{legend_analyse}%
	\input{figures/legend_analyse.tikz}%
	\tikzexternaldisable
\\
\setlength{\figureheight}{5.3cm}
\setlength{\figurewidth}{5.3cm}
\subfloat[][{\makebox[-1cm]{}}]{%
	\tikzexternalenable
	\tikzsetnextfilename{G2F1_analyse_E}%
	\input{figures/G2F1_analyse_E.tikz}%
	\tikzexternaldisable
}\hspace{0.2cm}
\subfloat[][{\makebox[-1cm]{}}]{%
	\tikzexternalenable
	\tikzsetnextfilename{G2F1_analyse_V}%
	\input{figures/G2F1_analyse_V.tikz}%
	\tikzexternaldisable
}\\
\subfloat[][{\makebox[-1cm]{}}]{%
	\tikzexternalenable
	\tikzsetnextfilename{G2F2_analyse_E}%
	\input{figures/G2F2_analyse_E.tikz}%
	\tikzexternaldisable
}\hspace{0.2cm}
\subfloat[][{\makebox[-1cm]{}}]{%
	\tikzexternalenable
	\tikzsetnextfilename{G2F2_analyse_V}%
	\input{figures/G2F2_analyse_V.tikz}%
	\tikzexternaldisable
}\\
\subfloat[][{\makebox[-1cm]{}}]{%
	\tikzexternalenable
	\tikzsetnextfilename{G2F3_analyse_E}%
	\input{figures/G2F3_analyse_E.tikz}%
	\tikzexternaldisable
}\hspace{0.2cm}
\subfloat[][{\makebox[-1cm]{}}]{%
	\tikzexternalenable
	\tikzsetnextfilename{G2F3_analyse_V}%
	\input{figures/G2F3_analyse_V.tikz}%
	\tikzexternaldisable
}
\caption{\review{Behavior of the estimated mean and variance of the multi-index differences $\Delta G_\bell$ for the second quantity of interest \texttt{G2}. The notation $(\ell,0,0)$ means that we investigate how these quantities behave when only the first dimension is refined, and similar for all the other (mixed) dimensions. For reference, the full black line corresponds to the approximation of the quantity of interest $G_\bell$ at $(\ell,\ell,\ell)$.}}
\label{fig:analyse_G2}
\end{figure}
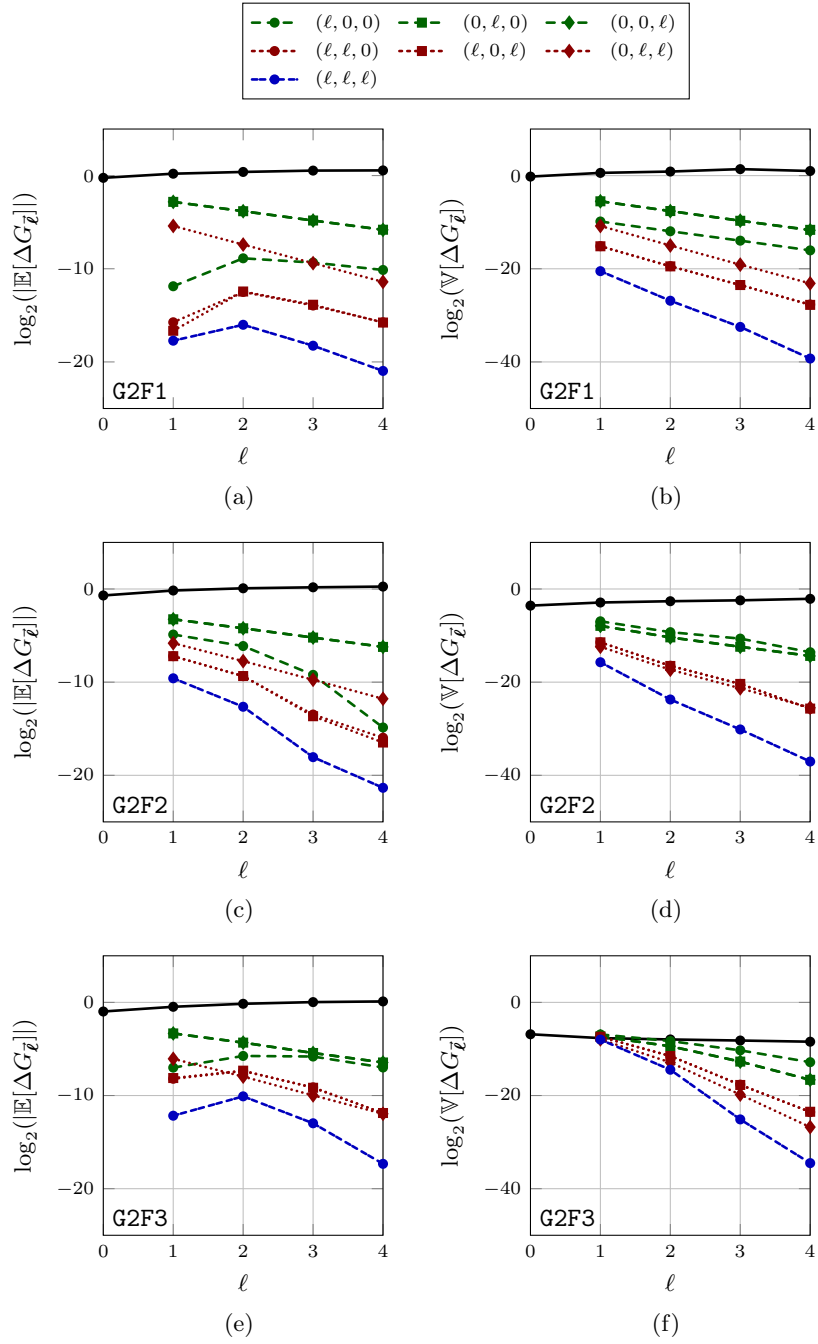

\review{Again, we analyze the behavior of the mean $|\E[\Delta G_\bell]|$ and variance $\V[\Delta G_\bell]$ of the multi-index differences in all (mixed) directions of this second problem. Now, there is a clear anisotropy when refining the differences as can be seen in~\figref{fig:analyse_G2}. It turns out that refining in the $x_1$-direction is much more advantageous than refining in the $x_2$- or $x_3$-direction. This is evident when we consider the \sr{asymmetry of the quantity of interest \texttt{G2}, where the flux in the $x_1$-direction is considered.} Observe that the convergence rates of the mean $|\E[\Delta G_\bell]|$ and the variance $\V[\Delta G_\bell]$ are much smaller compared to the first test problem, \texttt{G1}.}

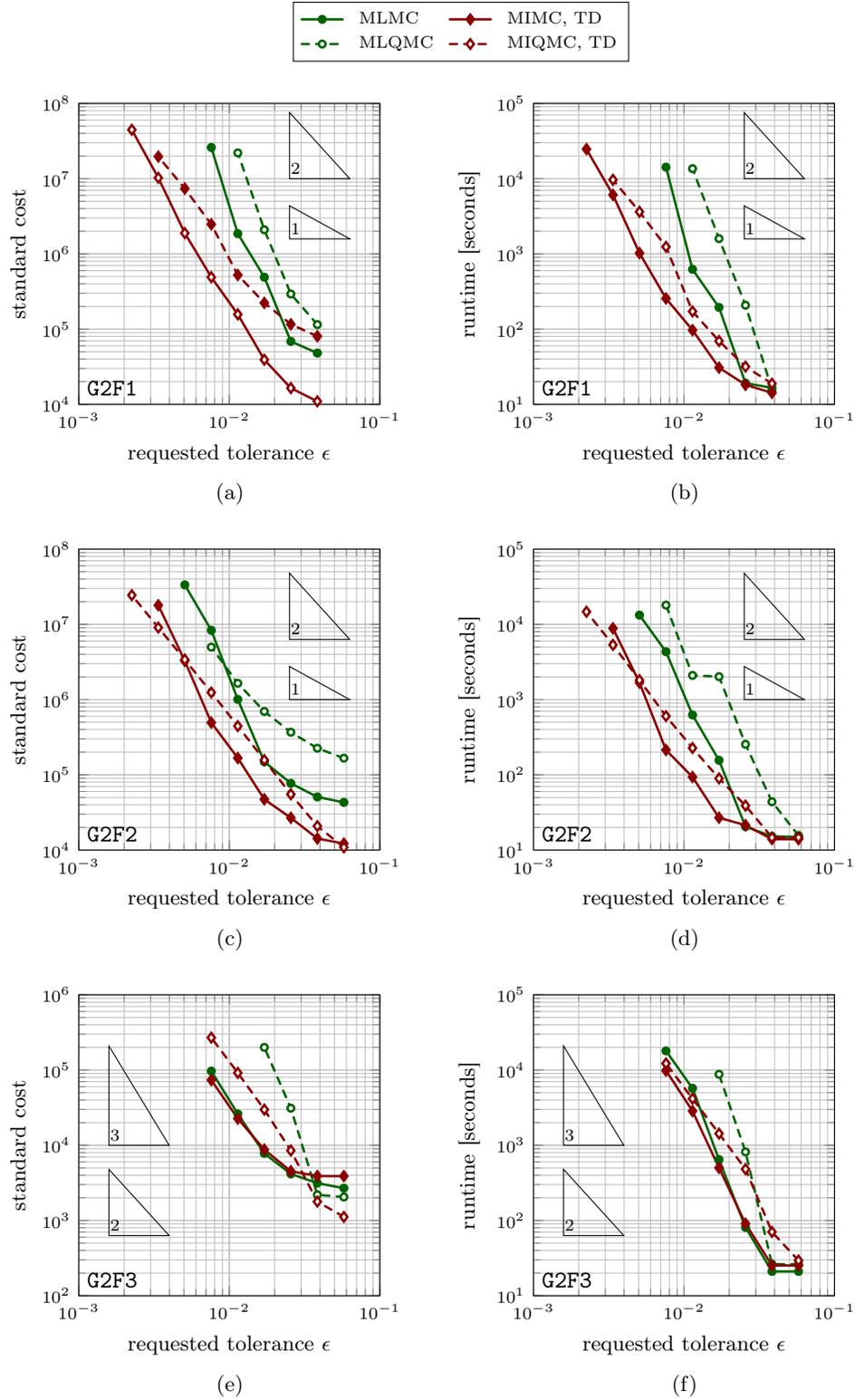
\begin{figure}[p]
\centering
\setlength{\figureheight}{2cm}
\setlength{\figurewidth}{10cm}
	\tikzexternalenable
	\tikzsetnextfilename{legend_times_G2}%
	\input{figures/legend_times_G2.tikz}%
	\tikzexternaldisable
\\
\setlength{\figureheight}{6cm}
\setlength{\figurewidth}{6cm}
\subfloat[][{\makebox[0cm]{}}]{%
	\tikzexternalenable
	\tikzsetnextfilename{cost_G2F1}%
	\input{figures/cost_G2F1.tikz}%
	\tikzexternaldisable
}\hspace{2.0cm}
\subfloat[][{\makebox[0cm]{}}]{%
	\tikzexternalenable
	\tikzsetnextfilename{times_G2F1}%
	\input{figures/times_G2F1.tikz}%
	\tikzexternaldisable
}\\
\subfloat[][{\makebox[0cm]{}}]{%
	\tikzexternalenable
	\tikzsetnextfilename{cost_G2F2}%
	\input{figures/cost_G2F2.tikz}%
	\tikzexternaldisable
}\hspace{2.0cm}
\subfloat[][{\makebox[0cm]{}}]{%
	\tikzexternalenable
	\tikzsetnextfilename{times_G2F2}%
	\input{figures/times_G2F2.tikz}%
	\tikzexternaldisable
}\\
\subfloat[][{\makebox[0cm]{}}]{%
	\tikzexternalenable
	\tikzsetnextfilename{cost_G2F3}%
	\input{figures/cost_G2F3.tikz}%
	\tikzexternaldisable
}\hspace{2.0cm}
\subfloat[][{\makebox[0cm]{}}]{%
	\tikzexternalenable
	\tikzsetnextfilename{times_G2F3}%
	\input{figures/times_G2F3.tikz}%
	\tikzexternaldisable
}\\
\caption{\review{Performance comparison of Multi-Index (both FT and TD) and Multilevel (Quasi) Monte Carlo under different test conditions for \texttt{G2}. Notice that the achieved accuracies are much \sr{higher} compared to the first test case \texttt{G1}.}}
\label{fig:G2}
\end{figure}

\review{As before, we run all methods for a sequence of decreasing tolerances $\epsilon$ and compare the performance. This is illustrated in~\figref{fig:G2}. However, we do not plot the results for the FT index set, since it behaves quite \sr{badly, similar to} our previous example. For the smooth case, \texttt{F1}, the benefit of multi-index methods over multilevel methods is clearly visible again. The multilevel methods have an asymptotic convergence rate $r=0.32$, thus, a $\mathrm{cost}=\mathcal{O}(\epsilon^{-3})$. This is due to the slow convergence rate of the variance of the multilevel differences. By \sr{also considering the other meshes} included in the multi-index telescoping sum in~\eqref{eq:MIMC} or~\eqref{eq:MIQMC}, we are again able to recover the order-2 $\epsilon$-convergence rate: $\mathrm{cost}=\mathcal{O}(\epsilon^{-2})$. For the less smooth case, \texttt{F2}, the results indicate that there might be some benefit in switching to QMC. However, the gain is only apparent for very small tolerances $\epsilon$. In the non-smooth case \texttt{F3}, we have again that all methods have the same asymptotic $\mathrm{cost}=\mathcal{O}(\epsilon^{-3})$. In this case, the extension to MIMC does not help.}

\review{We conclude from our experiments that the MIQMC estimator with TD index sets is able to recover the best possible $\epsilon$-convergence rate of order 1 when using randomly shifted lattice rules for smooth problems. For problems with less smoothness, and hence, a slower convergence of the variance of the multi-index differences, the convergence rate \sr{deteriorates}, but in our experiments the asymptotic $\epsilon$-complexity is never worse than the classical MLMC method.}

%%%%%%%%%%%%%%%%%%%%%%%%%%%%%%%%%%%%%%%%%%%%%%%%%%%%%%%%%%%
% 6 - CONCLUSIONS AND FURTHER WORK
%%%%%%%%%%%%%%%%%%%%%%%%%%%%%%%%%%%%%%%%%%%%%%%%%%%%%%%%%%%
\section{Conclusions and Further Work}

We have proposed a new Multi-Index Quasi-Monte Carlo algorithm for the solution of mathematical models in the form of partial differential equations with random coefficients. The MIQMC method combines the Multi-Index Monte Carlo method developed in~\cite{haji2016multi} with QMC methods to \sr{obtain faster convergence} of the multi-index differences. Motivated by problems in subsurface flow, we applied our method to an elliptic PDE in three dimensions with a diffusion coefficient given by a lognormal random field with underlying Mat\'ern covariance function. For problems with a small correlation length, several thousand uncertainties are required to accurately model the underlying random field. 

Our numerical results show that the MIQMC method performs remarkably well in the case of smooth problems. We are able to recover a cost $\mathcal{O}(1/\epsilon)$, associated with QMC methods, where $\epsilon$ is the requested tolerance on the estimator. \review{For problems with less smoothness, the benefit of QMC is less pronounced, as expected.} Because of the faster QMC convergence, we are able to reduce the simulation time from \sr{several hours with MLMC to only a couple of minutes with MIQMC}, for problems with $2$ million degrees of freedom and up to 3500 uncertainties.

Future work may focus on either further reducing the variance of the multi-index differences, using for example interlaced polynomial lattice rules~\cite{dick2013higher} or higher order digital nets~\cite{dick2010digital}. Provided enough smoothness in the problem, we expect the estimator to achieve a complexity $\mathcal{O}(1/\epsilon^{p})$, with $p<2$. One could also improve the multi-index method by introducing an adaptive strategy for choosing the indices. This is already done in~\cite{robbe2017dimension} for the Multi-Index Monte Carlo method, based on a greedy algorithm outlined in~\cite{haji2016misc,gerstner2003dimension,nobile2015convergence}. Using this approach, a quasi-optimal index set is constructed, without prior knowledge of the underlying problem. For this adaptive method, we expect similar gains as for adaptive sparse grids~\cite{gerstner2003dimension}.

Finally, we note that \review{the analysis of our MIQMC estimator should be similar to what was done in~\cite{haji2016multi}, except for the convergence rate of the sampling method}. Such analysis is outside the scope of the present paper.

\section*{Acknowledgments}

The authors would like to thank the referees for their valuable input. Their detailed comments and remarks helped to significantly improve the manuscript.

\bibliographystyle{abbrv}
\bibliography{references}

\end{document}

%% file: paper_template.tex
\usepackage[T1]{fontenc}
\usepackage{amssymb, amsmath, mathtools}
\usepackage{graphicx}
\usepackage{ pgfplots}
	\pgfplotsset{compat=newest}
	\newlength\figureheight
	\newlength\figurewidth
\usepackage{pgfplotstable}
\usepackage{tikz}
	\usetikzlibrary{external}
	\usetikzlibrary{backgrounds}
  	\usetikzlibrary{shapes.geometric}
\usepackage{mathrsfs}
\usepackage{hyperref}
\usepackage{algorithm}
\usepackage{caption}
\usepackage[caption=false]{subfig}
\usepackage{cite}
\usepackage{algorithmic}
	
\usepackage{pifont}
\usepackage{multirow}
\usepackage{booktabs}
\usepackage[top=1in, bottom=1.25in, left=1.25in, right=1.25in]{geometry}

\newcommand{\figref}[1]{Figure~\ref{#1}}
\newcommand{\tabref}[1]{Table~\ref{#1}}
\newcommand{\algref}[1]{Algorithm~\ref{#1}}

\newcommand{\secref}[1]{\S\ref{#1}}

\newcommand{\R}{\mathbb{R}}
\newcommand{\N}{\mathbb{N}}
\newcommand{\I}{\mathcal{I}}
\newcommand{\calI}{\mathcal{I}}
\newcommand{\bx}{\boldsymbol{x}}
\newcommand{\by}{\boldsymbol{y}}
\newcommand{\bt}{\boldsymbol{t}}
\newcommand{\bmu}{\boldsymbol{\mu}}
\newcommand{\btau}{{\vec{\boldsymbol{\tau}}}}
\newcommand{\bell}{{\vec{\boldsymbol{\ell}}}}

\newcommand{\E}{\mathbb{E}}
\newcommand{\V}{\mathbb{V}}

\newcommand{\TheTitle}[1]{\def\TheTitle{#1}}
\newcommand{\TheAuthors}[1]{\def\TheAuthors{#1}}

\definecolor{green}{RGB}{0,150,0}
\definecolor{darkgrey}{rgb}{0.8,0.8,0.8}
\definecolor{darkred}{RGB}{139,0,0}
\definecolor{darkgreen}{RGB}{0,100,0}
\definecolor{darkmagenta}{RGB}{180,0,180}
\definecolor{darkblue}{RGB}{0,0,190}
\definecolor{ratecolor1}{RGB}{0,100,0}
\definecolor{ratecolor2}{RGB}{139,0,0}
\definecolor{ratecolor3}{RGB}{0,0,190}

\newcommand{\review}[1]{{\color{black}{#1}}}
\newcommand{\sr}[1]{#1}
\newcolumntype{p}[1]{>{\centering\let\newline\\\arraybackslash\hspace{0pt}}m{#1}}
\newcolumntype{r}[1]{>{\raggedright\let\newline\\\arraybackslash\hspace{0pt}}m{#1}}

\pgfdeclareplotmark{(ml)}
  {\draw[ratecolor1,fill] (0pt, 0pt) circle [radius=1.4pt];
  \draw[white,fill] (0pt, 0pt) circle [radius=0.7pt];}
\pgfdeclareplotmark{(td)}
  {%\node [draw,scale=0.5,diamond,fill=ratecolor2] at (-0.05,0) {}; 
  \node [draw,scale=0.4,diamond,fill=white] at (-0.07,0) {};}
\pgfdeclareplotmark{(ft)}
  {\draw[ratecolor3,fill] (0pt, 0pt) diamond [radius=1.4pt];
  \draw[white,fill] (0pt, 0pt) diamond [radius=0.7pt];}

\newcommand{\LogLogSlopeTriangle}[5]
{
    \pgfplotsextra
    {
        \pgfkeysgetvalue{/pgfplots/xmin}{\xmin}
        \pgfkeysgetvalue{/pgfplots/xmax}{\xmax}
        \pgfkeysgetvalue{/pgfplots/ymin}{\ymin}
        \pgfkeysgetvalue{/pgfplots/ymax}{\ymax}

        \pgfmathsetmacro{\xArel}{#1}
        \pgfmathsetmacro{\yArel}{#3}
        \pgfmathsetmacro{\xBrel}{#1-#2}
        \pgfmathsetmacro{\yBrel}{\yArel}
        \pgfmathsetmacro{\xCrel}{\xBrel}

        \pgfmathsetmacro{\lnxB}{\xmin*(1-(#1-#2))+\xmax*(#1-#2)} % in [xmin,xmax]
        \pgfmathsetmacro{\lnxA}{\xmin*(1-#1)+\xmax*#1} % in [xmin,xmax]
        \pgfmathsetmacro{\lnyA}{\ymin*(1-#3)+\ymax*#3} % in [ymin,ymax]
        \pgfmathsetmacro{\lnyC}{\lnyA+1.1*#4*(\lnxA-\lnxB)}
        \pgfmathsetmacro{\yCrel}{\lnyC-\ymin)/(\ymax-\ymin)}
        
        \coordinate (A) at (rel axis cs:\xArel,\yArel);
        \coordinate (B) at (rel axis cs:\xBrel,\yBrel);
        \coordinate (C) at (rel axis cs:\xCrel,\yCrel);

        \draw[#5]   (A)--node[pos=0.9,yshift=1ex] {\scriptsize #4}
                    (B)--
                    (C)-- 
                    cycle;
    }
}

\makeatletter
\newenvironment{customlegend}[1][]{
    \begingroup
    \pgfplots@init@cleared@structures
    \pgfplotsset{#1}
}{
    \pgfplots@createlegend
    \endgroup
}
\def\addlegendimage{\pgfplots@addlegendimage}
\makeatother

%% file: figures/colorbar_3.tikz
\begin{tikzpicture}
\begin{axis}[
    hide axis,
    scale only axis,
    height=0pt,
    width=0pt,
    colormap/jet,
    colorbar,
    point meta min=-3,
    point meta max=3,
    colorbar style={
        height=\figureheight,
        width=\figurewidth,
        ytick={-3,-1,1,3},
        yticklabel style={text width=width("$-3$"),align=right},
        axis line style={draw=none},
        ytick style={draw=none},
        every y tick label/.append style={font=\tiny, xshift=-1.5ex}
    }
    ]
    \addplot [draw=none] coordinates {(0,0)};
\end{axis}
\end{tikzpicture}

%% file: figures/decay.tikz
\begin{tikzpicture}
	\begin{axis}[
		width=\figurewidth,
		height=\figureheight,
		xmode=log,
		xmin=1,
		xmax=1000000,
		xminorticks=true,
		xlabel={\scriptsize $r$},
		every x tick label/.append style={font=\tiny},
		ymode=log,
		ymin=1e-7,
		ymax=1,
		yminorticks=true,
		ytick={1e0,1e-1,1e-2,1e-3,1e-4,1e-5,1e-6,1e-7},
		every y tick label/.append style={font=\tiny},
		ylabel={\scriptsize  $\theta_r$},
		axis background/.style={fill=white},
		legend style={legend cell align=left,font=\scriptsize}
	]
	
	\pgfplotstableread{data/eigenvalues/lambda1sigma1nu2pt5.txt} \data;
	\addplot [
		color=ratecolor1,
		solid,
		line width=1pt
	]
	table[
		x expr=\thisrow{x},
		y expr=\thisrow{y}
	] {\data};
	\addlegendentry{$\{1,1,2.5\}$};
	
	\pgfplotstableread{data/eigenvalues/lambda0pt3sigma1nu1.txt} \data;
	\addplot [
		color=ratecolor2,
		dashed,
		line cap=round,
		line width=1pt
	]
	table[
		x expr=\thisrow{x},
		y expr=\thisrow{y}
	] {\data};
	\addlegendentry{$\{0.3,1,1\}$};
	
	\pgfplotstableread{data/eigenvalues/lambda0pt075sigma1nu0pt5.txt} \data;
	\addplot [
		color=ratecolor3,
		dashed,
		dash pattern=on 5pt off 3pt,
		line width=1pt,
		line cap=round
	]
	table[
		x expr=\thisrow{x},
		y expr=\thisrow{y}
	] {\data};
	\addlegendentry{$\{0.075,1,0.5\}$};
	
	\addplot[color=black,dashed,line cap=round,line width=0.75pt,domain=1:1e6] {25*x^(-3/2)};
	\addplot[color=black,dashed,line cap=round,line width=0.75pt,domain=1:1e6] {50*x^(-2.1)};
	\addplot[color=black,dashed,line cap=round,line width=0.75pt,domain=1:1e6] {75*x^(-7/2-0.2)};
		
	\pgfplotsset{
    after end axis/.code={
        \node[black,above] at (axis cs:350,0.75e-6){\scriptsize $r^{-8/3}$};
        \node[black,above] at (axis cs:1.25e4,0.75e-6){\scriptsize $r^{-5/3}$};
        \node[black,above] at (axis cs:2.5e5,0.75e-6){\scriptsize $r^{-4/3}$};  
    }
}
	\end{axis}
\end{tikzpicture}

%% file: figures/MIMCgrid.tikz
\begin{tikzpicture}[scale=0.25,framed]
	% bottom line
	\draw[xstep=4cm,ystep=1cm] (0,0) grid (8,8);
		\foreach \x in {2,6}
			\foreach \y in {0.5,1.5,2.5,3.5,4.5,5.5,6.5,7.5}
				\node[shape=circle,fill=black,scale=0.3] (\x-\y) at (\x,\y){};
	\draw[thick,latex-] (9,4) -- (11,4);
	\draw[xstep=2cm,ystep=1cm] (12,0) grid (20,8);
		\foreach \x in {1,3,5,7}
			\foreach \y in {0.5,1.5,2.5,3.5,4.5,5.5,6.5,7.5}
				\node[shape=circle,fill=black,scale=0.3] (\x-\y) at (12+\x,\y){};
	\draw[thick,latex-] (21,4) -- (23,4);
	\draw[step=1cm] (24,0) grid (32,8);
		\foreach \x in {0.5,1.5,2.5,3.5,4.5,5.5,6.5,7.5}
			\foreach \y in {0.5,1.5,2.5,3.5,4.5,5.5,6.5,7.5}
				\node[shape=circle,fill=black,scale=0.3] (\x-\y) at (24+\x,\y){};
	\draw[thick,latex-] (4,11) -- (4,9);
	\draw[thick,latex-] (16,11) -- (16,9);
	\draw[thick,latex-] (28,11) -- (28,9);
	\draw[thick,latex-] (21,11) -- (23,9);
	\draw[thick,latex-] (9,11) -- (11,9);
	% middle line
	\draw[xstep=4,ystep=2cm] (0,12) grid (8,20);
		\foreach \x in {2,6}
			\foreach \y in {1,3,5,7}
				\node[shape=circle,fill=black,scale=0.3] (\x-\y) at (\x,12+\y){};
	\draw[thick,latex-] (9,16) -- (11,16);
	\draw[step=2cm] (12,12) grid (20,20);
		\foreach \x in {1,3,5,7}
			\foreach \y in {1,3,5,7}
				\node[shape=circle,fill=black,scale=0.3] (\x-\y) at (12+\x,12+\y){};
	\draw[thick,latex-] (21,16) -- (23,16);
	\draw[xstep=1cm,ystep=2cm] (24,12) grid (32,20);
		\foreach \x in {0.5,1.5,2.5,3.5,4.5,5.5,6.5,7.5}
			\foreach \y in {1,3,5,7}
				\node[shape=circle,fill=black,scale=0.3] (\x-\y) at (24+\x,12+\y){};
	\draw[thick,latex-] (9,23) -- (11,21);
	\draw[thick,latex-] (21,23) -- (23,21);
	\draw[thick,latex-] (4,23) -- (4,21);
	\draw[thick,latex-] (16,23) -- (16,21);
	\draw[thick,latex-] (28,23) -- (28,21);
	% top line
	\draw[step=4cm] (0,24) grid (8,32);
		\foreach \x in {2,6}
			\foreach \y in {2,6}
				\node[shape=circle,fill=black,scale=0.3] (\x-\y) at (\x,24+\y){};
	\draw[thick,latex-] (9,28) -- (11,28);
	\draw[ystep=4,xstep=2cm] (12,24) grid (20,32);
		\foreach \x in {1,3,5,7}
			\foreach \y in {2,6}
				\node[shape=circle,fill=black,scale=0.3] (\x-\y) at (12+\x,24+\y){};
	\draw[thick,latex-] (21,28) -- (23,28);
	\draw[xstep=1cm,ystep=4cm] (24,24) grid (32,32);
		\foreach \x in {0.5,1.5,2.5,3.5,4.5,5.5,6.5,7.5}
			\foreach \y in {2,6}
				\node[shape=circle,fill=black,scale=0.3] (\x-\y) at (24+\x,24+\y){};
	% outline
	\foreach \y in {0.5,1.5,2.5,3.5,4.5,5.5,6.5,7.5}
		\node[shape=circle,fill=black,scale=0.3] (\y) at (-2,\y){};
	\foreach \y in {1,3,5,7}
		\node[shape=circle,fill=black,scale=0.3] (\y) at (-2,12+\y){};
	\foreach \y in {2,6}
		\node[shape=circle,fill=black,scale=0.3] (\y) at (-2,24+\y){};
	\foreach \x in {0.5,1.5,2.5,3.5,4.5,5.5,6.5,7.5}
		\node[shape=circle,fill=black,scale=0.3] (\x) at (24+\x,34){};
	\foreach \x in {1,3,5,7}
		\node[shape=circle,fill=black,scale=0.3] (\x) at (12+\x,34){};
	\foreach \x in {2,6}
		\node[shape=circle,fill=black,scale=0.3] (\x) at (\x,34){};
	% text
	\node at (4,36) {$\ell_1=0$};
	\node at (16,36) {$\ell_1=1$};
	\node at (28,36) {$\ell_1=2$};
	\node[rotate=90] at (-4,4) {$\ell_2=2$};
	\node[rotate=90] at (-4,16) {$\ell_2=1$};
	\node[rotate=90] at (-4,28) {$\ell_2=0$};
	\node[] at (35,-3) {};
\end{tikzpicture}

%% file: figures/G1F1_analyse_W.tikz
\begin{tikzpicture}
\def\maxL{4}
\begin{axis}[
	height=\figureheight,
	width=\figurewidth,
	xmin=0,
	xmax=\maxL,
	xminorticks=true,
	xlabel={\small $\ell$},
	every x tick label/.append style={font=\scriptsize},
	ymin=-15,
	ymax=5,
	yminorticks=true,
	ylabel={\small $\log_2(\text{time}_\ell)$},
	every y tick label/.append style={font=\scriptsize},
	legend style ={ at={(1.03,1)}, 
        anchor=north west, draw=black, style={font=\scriptsize},
        fill=white,align=left}
]

%
% SINGLES
%

	\pgfplotstableread{data/G1F1/Winn.txt} \data;
	\addplot [
		color=ratecolor1,
		dashed,
		line width=1pt,
		line cap=round,
		mark=*,
		mark size = 1.4,
		mark options={solid,ratecolor1}
	]
	table[
		x expr=\thisrow{x},
		y expr=\thisrow{y}
	] {\data};
	\addlegendentry{$(\ell,0,0)$};
	
	\pgfplotstableread{data/G1F1/Wnin.txt} \data;
	\addplot [
		color=ratecolor1,
		dashed,
		line width=1pt,
		line cap=round,
		mark=square*,
		mark size = 1.4,
		mark options={solid,ratecolor1}
	]
	table[
		x expr=\thisrow{x},
		y expr=\thisrow{y}
	] {\data};
	\addlegendentry{$(0,\ell,0)$};
	
	\pgfplotstableread{data/G1F1/Wnni.txt} \data;
	\addplot [
		color=ratecolor1,
		dashed,
		line width=1pt,
		line cap=round,
		mark=diamond*,
		mark options={solid,ratecolor1}
	]
	table[
		x expr=\thisrow{x},
		y expr=\thisrow{y}
	] {\data};
	\addlegendentry{$(0,0,\ell)$};
	
%
% DOUBLES
%

	\pgfplotstableread{data/G1F1/Wiin.txt} \data;
	\addplot [
		color=ratecolor2,
		dotted,
		line width=1pt,
		line cap=round,
		dash pattern=on 0pt off 2\pgflinewidth,
		mark=*,
		mark size = 1.4,
		mark options={solid,ratecolor2}
	]
	table[
		x expr=\thisrow{x},
		y expr=\thisrow{y}
	] {\data};
	\addlegendentry{$(\ell,\ell,0)$};
	
	\pgfplotstableread{data/G1F1/Wini.txt} \data;
	\addplot [
		color=ratecolor2,
		dotted,
		line width=1pt,
		line cap=round,
		dash pattern=on 0pt off 2\pgflinewidth,
		mark=square*,
		mark size = 1.4,
		mark options={solid,ratecolor2}
	]
	table[
		x expr=\thisrow{x},
		y expr=\thisrow{y}
	] {\data};
	\addlegendentry{$(\ell,0,\ell)$};
	
	\pgfplotstableread{data/G1F1/Wnii.txt} \data;
	\addplot [
		color=ratecolor2,
		dotted,
		line width=1pt,
		line cap=round,
		dash pattern=on 0pt off 2\pgflinewidth,
		mark=diamond*,
		mark options={solid,ratecolor2}
	]
	table[
		x expr=\thisrow{x},
		y expr=\thisrow{y}
	] {\data};
	\addlegendentry{$(0,\ell,\ell)$};

%
% TRIPPLE
%
	\pgfplotstableread{data/G1F1/Wiii.txt} \data;
	\addplot [
		color=ratecolor3,
		dash pattern=on 2pt off 1pt on 2pt off 2pt,
		line width=1pt,
		line cap=round,
		mark=*,
		mark size = 1.4,
		mark options={solid,ratecolor3}
	]
	table[
		x expr=\thisrow{x},
		y expr=\thisrow{y}
	] {\data};
	\addlegendentry{$(\ell,\ell,\ell)$};
	
\end{axis}
\end{tikzpicture}

%% file: figures/legend_analyse.tikz
\begin{tikzpicture}
    \begin{customlegend}[
    	legend columns=3,
		legend style={column sep=1ex,font=\scriptsize},
		legend cell align=left,
		legend entries={$(\ell,0,0)$ \\ $(0,\ell,0)$ \\ $(0,0,\ell)$ \\ $(\ell,\ell,0)$ \\ $(\ell,0,\ell)$ \\ $(0,\ell,\ell)$ \\ $(\ell,\ell,\ell)$\\}
		]
		\addlegendimage{ratecolor1,dashed,line cap=round,line width=1pt,mark=*,mark size=1.4pt,mark options={solid}}
		\addlegendimage{ratecolor1,dashed,line cap=round,line width=1pt,mark=square*,mark size=1.4pt,mark options={solid}}
		\addlegendimage{ratecolor1,dashed,line cap=round,line width=1pt,mark=diamond*,mark options={solid}}
		\addlegendimage{ratecolor2,dotted,line cap=round,line width=1pt,mark=*,mark size=1.4pt,mark options={solid}}
		\addlegendimage{ratecolor2,dotted,line cap=round,line width=1pt,mark=square*,mark size=1.4pt,mark options={solid}}
		\addlegendimage{ratecolor2,dotted,line cap=round,line width=1pt,mark=diamond*,mark options={solid}}
		\addlegendimage{ratecolor3,dash pattern=on 2pt off 1pt on 2pt off 2pt,line width=1pt,mark=*,mark size=1.4pt,mark options={solid},line cap=round}
   \end{customlegend}
\end{tikzpicture}

%% file: figures/G1F1_analyse_E.tikz
\begin{tikzpicture}
\def\maxL{4}
\begin{axis}[
	height=\figureheight,
	width=\figurewidth,
	xmin=0,
	xmax=\maxL,
	xminorticks=true,
	xlabel={\small $\ell$},
	every x tick label/.append style={font=\scriptsize},
	ymin=-40,
	ymax=0,
	yminorticks=true,
	ylabel={\small $\log_2(|\E[\Delta G_\bell]|)$},
	every y tick label/.append style={font=\scriptsize},
	legend style ={ at={(1.03,1)}, 
        anchor=north west, draw=black, style={font=\scriptsize},
        fill=white,align=left},
   ymajorgrids,
	xmajorgrids
]

%
% ZERO
%

	\pgfplotstableread{data/G1F1/Ennn.txt} \data;
	\addplot [
		color=black,
		solid,
		line width=1pt,
		line cap=round,
		mark=*,
		mark size = 1.4,
		mark options={solid}
	]
	table[
		x expr=\thisrow{x},
		y expr=\thisrow{y}
	] {\data};
	%\addlegendentry{$(0,0,0)$};

%
% SINGLES
%

	\pgfplotstableread{data/G1F1/Einn.txt} \data;
	\addplot [
		color=ratecolor1,
		dashed,
		line width=1pt,
		line cap=round,
		mark=*,
		mark size = 1.4,
		mark options={solid,ratecolor1}
	]
	table[
		x expr=\thisrow{x},
		y expr=\thisrow{y}
	] {\data};
	%\addlegendentry{$(\ell,0,0)$};
	
	\pgfplotstableread{data/G1F1/Enin.txt} \data;
	\addplot [
		color=ratecolor1,
		dashed,
		line width=1pt,
		line cap=round,
		mark=square*,
		mark size = 1.4,
		mark options={solid,ratecolor1}
	]
	table[
		x expr=\thisrow{x},
		y expr=\thisrow{y}
	] {\data};
	%\addlegendentry{$(0,\ell,0)$};
	
	\pgfplotstableread{data/G1F1/Enni.txt} \data;
	\addplot [
		color=ratecolor1,
		dashed,
		line width=1pt,
		line cap=round,
		mark=diamond*,
		mark options={solid,ratecolor1}
	]
	table[
		x expr=\thisrow{x},
		y expr=\thisrow{y}
	] {\data};
	%\addlegendentry{$(0,0,\ell)$};
	
%
% DOUBLES
%

	\pgfplotstableread{data/G1F1/Eiin.txt} \data;
	\addplot [
		color=ratecolor2,
		dotted,
		line width=1pt,
		line cap=round,
		dash pattern=on 0pt off 2\pgflinewidth,
		mark=*,
		mark size = 1.4,
		mark options={solid,ratecolor2}
	]
	table[
		x expr=\thisrow{x},
		y expr=\thisrow{y}
	] {\data};
	%\addlegendentry{$(\ell,\ell,0)$};
	
	\pgfplotstableread{data/G1F1/Eini.txt} \data;
	\addplot [
		color=ratecolor2,
		dotted,
		line width=1pt,
		line cap=round,
		dash pattern=on 0pt off 2\pgflinewidth,
		mark=square*,
		mark size = 1.4,
		mark options={solid,ratecolor2}
	]
	table[
		x expr=\thisrow{x},
		y expr=\thisrow{y}
	] {\data};
	%\addlegendentry{$(\ell,0,\ell)$};
	
	\pgfplotstableread{data/G1F1/Enii.txt} \data;
	\addplot [
		color=ratecolor2,
		dotted,
		line width=1pt,
		line cap=round,
		dash pattern=on 0pt off 2\pgflinewidth,
		mark=diamond*,
		mark options={solid,ratecolor2}
	]
	table[
		x expr=\thisrow{x},
		y expr=\thisrow{y}
	] {\data};
	%\addlegendentry{$(0,\ell,\ell)$};

%
% TRIPPLE
%

	\pgfplotstableread{data/G1F1/Eiii.txt} \data;
	\addplot [
		color=ratecolor3,
		dash pattern=on 2pt off 1pt on 2pt off 2pt,
		line width=1pt,
		line cap=round,
		mark=*,
		mark size = 1.4,
		mark options={solid,ratecolor3}
	]
	table[
		x expr=\thisrow{x},
		y expr=\thisrow{y}
	] {\data};
	%\addlegendentry{$(\ell,\ell,\ell)$};
	
		\pgfplotsset{
    	after end axis/.code={
        	\node[above right] at (0.03,0.03){\texttt{G1F1}};
    	}
	}
	
\end{axis}
\end{tikzpicture}

%% file: figures/G1F1_analyse_V.tikz
\begin{tikzpicture}
\def\maxL{4}
\begin{axis}[
	height=\figureheight,
	width=\figurewidth,
	xmin=0,
	xmax=\maxL,
	xminorticks=true,
	xlabel={\small $\ell$},
	every x tick label/.append style={font=\scriptsize},
	ymin=-80,
	ymax=0,
	yminorticks=true,
	ylabel={\small $\log_2(\V[\Delta G_\bell])$},
	every y tick label/.append style={font=\scriptsize},
	%legend style ={ at={(1.06,1)}, 
   %     anchor=north west, draw=black, style={font=\scriptsize},
   %     fill=white,align=left},
   ymajorgrids,
	xmajorgrids
]

%
% ZERO
%

	\pgfplotstableread{data/G1F1/Vnnn.txt} \data;
	\addplot [
		color=black,
		solid,
		line width=1pt,
		line cap=round,
		mark=*,
		mark size = 1.4,
		mark options={solid},
		forget plot
	]
	table[
		x expr=\thisrow{x},
		y expr=\thisrow{y}
	] {\data};
	%\addlegendentry{$(0,0,0)$};

%
% SINGLES
%

	\pgfplotstableread{data/G1F1/Vinn.txt} \data;
	\addplot [
		color=ratecolor1,
		dashed,
		line width=1pt,
		line cap=round,
		mark=*,
		mark size = 1.4,
		mark options={solid,ratecolor1}
	]
	table[
		x expr=\thisrow{x},
		y expr=\thisrow{y}
	] {\data};
	%\addlegendentry{$(\ell,0,0)$};
	
	\pgfplotstableread{data/G1F1/Vnin.txt} \data;
	\addplot [
		color=ratecolor1,
		dashed,
		line width=1pt,
		line cap=round,
		mark=square*,
		mark size = 1.4,
		mark options={solid,ratecolor1}
	]
	table[
		x expr=\thisrow{x},
		y expr=\thisrow{y}
	] {\data};
	%\addlegendentry{$(0,\ell,0)$};
	
	\pgfplotstableread{data/G1F1/Vnni.txt} \data;
	\addplot [
		color=ratecolor1,
		dashed,
		line width=1pt,
		line cap=round,
		mark=diamond*,
		mark options={solid,ratecolor1}
	]
	table[
		x expr=\thisrow{x},
		y expr=\thisrow{y}
	] {\data};
	%\addlegendentry{$(0,0,\ell)$};
	
%
% DOUBLES
%

	\pgfplotstableread{data/G1F1/Viin.txt} \data;
	\addplot [
		color=ratecolor2,
		dotted,
		line width=1pt,
		line cap=round,
		dash pattern=on 0pt off 2\pgflinewidth,
		mark=*,
		mark size = 1.4,
		mark options={solid,ratecolor2}
	]
	table[
		x expr=\thisrow{x},
		y expr=\thisrow{y}
	] {\data};
	%\addlegendentry{$(\ell,\ell,0)$};
	
	\pgfplotstableread{data/G1F1/Vini.txt} \data;
	\addplot [
		color=ratecolor2,
		dotted,
		line width=1pt,
		line cap=round,
		dash pattern=on 0pt off 2\pgflinewidth,
		mark=square*,
		mark size = 1.4,
		mark options={solid,ratecolor2}
	]
	table[
		x expr=\thisrow{x},
		y expr=\thisrow{y}
	] {\data};
	%\addlegendentry{$(\ell,0,\ell)$};
	
	\pgfplotstableread{data/G1F1/Vnii.txt} \data;
	\addplot [
		color=ratecolor2,
		dotted,
		line width=1pt,
		line cap=round,
		dash pattern=on 0pt off 2\pgflinewidth,
		mark=diamond*,
		mark options={solid,ratecolor2}
	]
	table[
		x expr=\thisrow{x},
		y expr=\thisrow{y}
	] {\data};
	%\addlegendentry{$(0,\ell,\ell)$};

%
% TRIPPLE
%

	\pgfplotstableread{data/G1F1/Viii.txt} \data;
	\addplot [
		color=ratecolor3,
		dash pattern=on 2pt off 1pt on 2pt off 2pt,
		line width=1pt,
		line cap=round,
		mark=*,
		mark size = 1.4,
		mark options={solid,ratecolor3}
	]
	table[
		x expr=\thisrow{x},
		y expr=\thisrow{y}
	] {\data};
	%\addlegendentry{$(\ell,\ell,\ell)$};
	
	\pgfplotsset{
    	after end axis/.code={
        	\node[above right] at (0.03,0.03){\texttt{G1F1}};
    	}
	}

\end{axis}
\end{tikzpicture}

%% file: figures/G1F2_analyse_E.tikz
\begin{tikzpicture}
\def\maxL{4}
\begin{axis}[
	height=\figureheight,
	width=\figurewidth,
	xmin=0,
	xmax=\maxL,
	xminorticks=true,
	xlabel={\small $\ell$},
	every x tick label/.append style={font=\scriptsize},
	ymin=-40,
	ymax=0,
	yminorticks=true,
	ylabel={\small $\log_2(|\E[\Delta G_\bell]|)$},
	every y tick label/.append style={font=\scriptsize},
	legend style ={ at={(1.03,1)}, 
        anchor=north west, draw=black, style={font=\scriptsize},
        fill=white,align=left},
   ymajorgrids,
	xmajorgrids
]

%
% ZERO
%

	\pgfplotstableread{data/G1F2/Ennn.txt} \data;
	\addplot [
		color=black,
		solid,
		line width=1pt,
		line cap=round,
		mark=*,
		mark size = 1.4,
		mark options={solid}
	]
	table[
		x expr=\thisrow{x},
		y expr=\thisrow{y}
	] {\data};
	%\addlegendentry{$(0,0,0)$};

%
% SINGLES
%

	\pgfplotstableread{data/G1F2/Einn.txt} \data;
	\addplot [
		color=ratecolor1,
		dashed,
		line width=1pt,
		line cap=round,
		mark=*,
		mark size = 1.4,
		mark options={solid,ratecolor1}
	]
	table[
		x expr=\thisrow{x},
		y expr=\thisrow{y}
	] {\data};
	%\addlegendentry{$(\ell,0,0)$};
	
	\pgfplotstableread{data/G1F2/Enin.txt} \data;
	\addplot [
		color=ratecolor1,
		dashed,
		line width=1pt,
		line cap=round,
		mark=square*,
		mark size = 1.4,
		mark options={solid,ratecolor1}
	]
	table[
		x expr=\thisrow{x},
		y expr=\thisrow{y}
	] {\data};
	%\addlegendentry{$(0,\ell,0)$};
	
	\pgfplotstableread{data/G1F2/Enni.txt} \data;
	\addplot [
		color=ratecolor1,
		dashed,
		line width=1pt,
		line cap=round,
		mark=diamond*,
		mark options={solid,ratecolor1}
	]
	table[
		x expr=\thisrow{x},
		y expr=\thisrow{y}
	] {\data};
	%\addlegendentry{$(0,0,\ell)$};
	
%
% DOUBLES
%

	\pgfplotstableread{data/G1F2/Eiin.txt} \data;
	\addplot [
		color=ratecolor2,
		dotted,
		line width=1pt,
		line cap=round,
		dash pattern=on 0pt off 2\pgflinewidth,
		mark=*,
		mark size = 1.4,
		mark options={solid,ratecolor2}
	]
	table[
		x expr=\thisrow{x},
		y expr=\thisrow{y}
	] {\data};
	%\addlegendentry{$(\ell,\ell,0)$};
	
	\pgfplotstableread{data/G1F2/Eini.txt} \data;
	\addplot [
		color=ratecolor2,
		dotted,
		line width=1pt,
		line cap=round,
		dash pattern=on 0pt off 2\pgflinewidth,
		mark=square*,
		mark size = 1.4,
		mark options={solid,ratecolor2}
	]
	table[
		x expr=\thisrow{x},
		y expr=\thisrow{y}
	] {\data};
	%\addlegendentry{$(\ell,0,\ell)$};
	
	\pgfplotstableread{data/G1F2/Enii.txt} \data;
	\addplot [
		color=ratecolor2,
		dotted,
		line width=1pt,
		line cap=round,
		dash pattern=on 0pt off 2\pgflinewidth,
		mark=diamond*,
		mark options={solid,ratecolor2}
	]
	table[
		x expr=\thisrow{x},
		y expr=\thisrow{y}
	] {\data};
	%\addlegendentry{$(0,\ell,\ell)$};

%
% TRIPPLE
%

	\pgfplotstableread{data/G1F2/Eiii.txt} \data;
	\addplot [
		color=ratecolor3,
		dash pattern=on 2pt off 1pt on 2pt off 2pt,
		line width=1pt,
		line cap=round,
		mark=*,
		mark size = 1.4,
		mark options={solid,ratecolor3}
	]
	table[
		x expr=\thisrow{x},
		y expr=\thisrow{y}
	] {\data};
	%\addlegendentry{$(\ell,\ell,\ell)$};
	
		\pgfplotsset{
    	after end axis/.code={
        	\node[above right] at (0.03,0.03){\texttt{G1F2}};
    	}
	}
	
\end{axis}
\end{tikzpicture}

%% file: figures/G1F2_analyse_V.tikz
\begin{tikzpicture}
\def\maxL{4}
\begin{axis}[
	height=\figureheight,
	width=\figurewidth,
	xmin=0,
	xmax=\maxL,
	xminorticks=true,
	xlabel={\small $\ell$},
	every x tick label/.append style={font=\scriptsize},
	ymin=-60,
	ymax=0,
	yminorticks=true,
	ylabel={\small $\log_2(\V[\Delta G_\bell])$},
	every y tick label/.append style={font=\scriptsize},
	legend style ={ at={(1.03,1)}, 
        anchor=north west, draw=black, style={font=\scriptsize},
        fill=white,align=left},
   ymajorgrids,
	xmajorgrids
]

%
% ZERO
%

	\pgfplotstableread{data/G1F2/Vnnn.txt} \data;
	\addplot [
		color=black,
		solid,
		line width=1pt,
		line cap=round,
		mark=*,
		mark size = 1.4,
		mark options={solid}
	]
	table[
		x expr=\thisrow{x},
		y expr=\thisrow{y}
	] {\data};
	%\addlegendentry{$(0,0,0)$};

%
% SINGLES
%

	\pgfplotstableread{data/G1F2/Vinn.txt} \data;
	\addplot [
		color=ratecolor1,
		dashed,
		line width=1pt,
		line cap=round,
		mark=*,
		mark size = 1.4,
		mark options={solid,ratecolor1}
	]
	table[
		x expr=\thisrow{x},
		y expr=\thisrow{y}
	] {\data};
	%\addlegendentry{$(\ell,0,0)$};
	
	\pgfplotstableread{data/G1F2/Vnin.txt} \data;
	\addplot [
		color=ratecolor1,
		dashed,
		line width=1pt,
		line cap=round,
		mark=square*,
		mark size = 1.4,
		mark options={solid,ratecolor1}
	]
	table[
		x expr=\thisrow{x},
		y expr=\thisrow{y}
	] {\data};
	%\addlegendentry{$(0,\ell,0)$};
	
	\pgfplotstableread{data/G1F2/Vnni.txt} \data;
	\addplot [
		color=ratecolor1,
		dashed,
		line width=1pt,
		line cap=round,
		mark=diamond*,
		mark options={solid,ratecolor1}
	]
	table[
		x expr=\thisrow{x},
		y expr=\thisrow{y}
	] {\data};
	%\addlegendentry{$(0,0,\ell)$};
	
%
% DOUBLES
%

	\pgfplotstableread{data/G1F2/Viin.txt} \data;
	\addplot [
		color=ratecolor2,
		dotted,
		line width=1pt,
		line cap=round,
		dash pattern=on 0pt off 2\pgflinewidth,
		mark=*,
		mark size = 1.4,
		mark options={solid,ratecolor2}
	]
	table[
		x expr=\thisrow{x},
		y expr=\thisrow{y}
	] {\data};
	%\addlegendentry{$(\ell,\ell,0)$};
	
	\pgfplotstableread{data/G1F2/Vini.txt} \data;
	\addplot [
		color=ratecolor2,
		dotted,
		line width=1pt,
		line cap=round,
		dash pattern=on 0pt off 2\pgflinewidth,
		mark=square*,
		mark size = 1.4,
		mark options={solid,ratecolor2}
	]
	table[
		x expr=\thisrow{x},
		y expr=\thisrow{y}
	] {\data};
	%\addlegendentry{$(\ell,0,\ell)$};
	
	\pgfplotstableread{data/G1F2/Vnii.txt} \data;
	\addplot [
		color=ratecolor2,
		dotted,
		line width=1pt,
		line cap=round,
		dash pattern=on 0pt off 2\pgflinewidth,
		mark=diamond*,
		mark options={solid,ratecolor2}
	]
	table[
		x expr=\thisrow{x},
		y expr=\thisrow{y}
	] {\data};
	%\addlegendentry{$(0,\ell,\ell)$};

%
% TRIPPLE
%

	\pgfplotstableread{data/G1F2/Viii.txt} \data;
	\addplot [
		color=ratecolor3,
		dash pattern=on 2pt off 1pt on 2pt off 2pt,
		line width=1pt,
		line cap=round,
		mark=*,
		mark size = 1.4,
		mark options={solid,ratecolor3}
	]
	table[
		x expr=\thisrow{x},
		y expr=\thisrow{y}
	] {\data};
	%\addlegendentry{$(\ell,\ell,\ell)$};
	
		\pgfplotsset{
    	after end axis/.code={
        	\node[above right] at (0.03,0.03){\texttt{G1F2}};
    	}
	}
	
\end{axis}
\end{tikzpicture}

%% file: figures/G1F3_analyse_E.tikz
\begin{tikzpicture}
\def\maxL{4}
\begin{axis}[
	height=\figureheight,
	width=\figurewidth,
	xmin=0,
	xmax=\maxL,
	xminorticks=true,
	xlabel={\small $\ell$},
	every x tick label/.append style={font=\scriptsize},
	ymin=-30,
	ymax=0,
	yminorticks=true,
	ylabel={\small $\log_2(|\E[\Delta G_\bell]|)$},
	every y tick label/.append style={font=\scriptsize},
	legend style ={ at={(1.03,1)}, 
        anchor=north west, draw=black, style={font=\scriptsize},
        fill=white,align=left},
   ymajorgrids,
	xmajorgrids
]

%
% ZERO
%

	\pgfplotstableread{data/G1F3/Ennn.txt} \data;
	\addplot [
		color=black,
		solid,
		line width=1pt,
		line cap=round,
		mark=*,
		mark size = 1.4,
		mark options={solid}
	]
	table[
		x expr=\thisrow{x},
		y expr=\thisrow{y}
	] {\data};
	%\addlegendentry{$(0,0,0)$};

%
% SINGLES
%

	\pgfplotstableread{data/G1F3/Einn.txt} \data;
	\addplot [
		color=ratecolor1,
		dashed,
		line width=1pt,
		line cap=round,
		mark=*,
		mark size = 1.4,
		mark options={solid,ratecolor1}
	]
	table[
		x expr=\thisrow{x},
		y expr=\thisrow{y}
	] {\data};
	%\addlegendentry{$(\ell,0,0)$};
	
	\pgfplotstableread{data/G1F3/Enin.txt} \data;
	\addplot [
		color=ratecolor1,
		dashed,
		line width=1pt,
		line cap=round,
		mark=square*,
		mark size = 1.4,
		mark options={solid,ratecolor1}
	]
	table[
		x expr=\thisrow{x},
		y expr=\thisrow{y}
	] {\data};
	%\addlegendentry{$(0,\ell,0)$};
	
	\pgfplotstableread{data/G1F3/Enni.txt} \data;
	\addplot [
		color=ratecolor1,
		dashed,
		line width=1pt,
		line cap=round,
		mark=diamond*,
		mark options={solid,ratecolor1}
	]
	table[
		x expr=\thisrow{x},
		y expr=\thisrow{y}
	] {\data};
	%\addlegendentry{$(0,0,\ell)$};
	
%
% DOUBLES
%

	\pgfplotstableread{data/G1F3/Eiin.txt} \data;
	\addplot [
		color=ratecolor2,
		dotted,
		line width=1pt,
		line cap=round,
		dash pattern=on 0pt off 2\pgflinewidth,
		mark=*,
		mark size = 1.4,
		mark options={solid,ratecolor2}
	]
	table[
		x expr=\thisrow{x},
		y expr=\thisrow{y}
	] {\data};
	%\addlegendentry{$(\ell,\ell,0)$};
	
	\pgfplotstableread{data/G1F3/Eini.txt} \data;
	\addplot [
		color=ratecolor2,
		dotted,
		line width=1pt,
		line cap=round,
		dash pattern=on 0pt off 2\pgflinewidth,
		mark=square*,
		mark size = 1.4,
		mark options={solid,ratecolor2}
	]
	table[
		x expr=\thisrow{x},
		y expr=\thisrow{y}
	] {\data};
	%\addlegendentry{$(\ell,0,\ell)$};
	
	\pgfplotstableread{data/G1F3/Enii.txt} \data;
	\addplot [
		color=ratecolor2,
		dotted,
		line width=1pt,
		line cap=round,
		dash pattern=on 0pt off 2\pgflinewidth,
		mark=diamond*,
		mark options={solid,ratecolor2}
	]
	table[
		x expr=\thisrow{x},
		y expr=\thisrow{y}
	] {\data};
	%\addlegendentry{$(0,\ell,\ell)$};

%
% TRIPPLE
%

	\pgfplotstableread{data/G1F3/Eiii.txt} \data;
	\addplot [
		color=ratecolor3,
		dash pattern=on 2pt off 1pt on 2pt off 2pt,
		line width=1pt,
		line cap=round,
		mark=*,
		mark size = 1.4,
		mark options={solid,ratecolor3}
	]
	table[
		x expr=\thisrow{x},
		y expr=\thisrow{y}
	] {\data};
	%\addlegendentry{$(\ell,\ell,\ell)$};
	
		\pgfplotsset{
    	after end axis/.code={
        	\node[above right] at (0.03,0.03){\texttt{G1F3}};
    	}
	}
	
\end{axis}
\end{tikzpicture}

%% file: figures/G1F3_analyse_V.tikz
\begin{tikzpicture}
\def\maxL{4}
\begin{axis}[
	height=\figureheight,
	width=\figurewidth,
	xmin=0,
	xmax=\maxL,
	xminorticks=true,
	xlabel={\small $\ell$},
	every x tick label/.append style={font=\scriptsize},
	ymin=-50,
	ymax=-10,
	yminorticks=true,
	ylabel={\small $\log_2(\V[\Delta G_\bell])$},
	every y tick label/.append style={font=\scriptsize},
	legend style ={ at={(1.03,1)}, 
        anchor=north west, draw=black, style={font=\scriptsize},
        fill=white,align=left},
   ymajorgrids,
	xmajorgrids
]

%
% ZERO
%

	\pgfplotstableread{data/G1F3/Vnnn.txt} \data;
	\addplot [
		color=black,
		solid,
		line width=1pt,
		line cap=round,
		mark=*,
		mark size = 1.4,
		mark options={solid}
	]
	table[
		x expr=\thisrow{x},
		y expr=\thisrow{y}
	] {\data};
	%\addlegendentry{$(0,0,0)$};

%
% SINGLES
%

	\pgfplotstableread{data/G1F3/Vinn.txt} \data;
	\addplot [
		color=ratecolor1,
		dashed,
		line width=1pt,
		line cap=round,
		mark=*,
		mark size = 1.4,
		mark options={solid,ratecolor1}
	]
	table[
		x expr=\thisrow{x},
		y expr=\thisrow{y}
	] {\data};
	%\addlegendentry{$(\ell,0,0)$};
	
	\pgfplotstableread{data/G1F3/Vnin.txt} \data;
	\addplot [
		color=ratecolor1,
		dashed,
		line width=1pt,
		line cap=round,
		mark=square*,
		mark size = 1.4,
		mark options={solid,ratecolor1}
	]
	table[
		x expr=\thisrow{x},
		y expr=\thisrow{y}
	] {\data};
	%\addlegendentry{$(0,\ell,0)$};
	
	\pgfplotstableread{data/G1F3/Vnni.txt} \data;
	\addplot [
		color=ratecolor1,
		dashed,
		line width=1pt,
		line cap=round,
		mark=diamond*,
		mark options={solid,ratecolor1}
	]
	table[
		x expr=\thisrow{x},
		y expr=\thisrow{y}
	] {\data};
	%\addlegendentry{$(0,0,\ell)$};
	
%
% DOUBLES
%

	\pgfplotstableread{data/G1F3/Viin.txt} \data;
	\addplot [
		color=ratecolor2,
		dotted,
		line width=1pt,
		line cap=round,
		dash pattern=on 0pt off 2\pgflinewidth,
		mark=*,
		mark size = 1.4,
		mark options={solid,ratecolor2}
	]
	table[
		x expr=\thisrow{x},
		y expr=\thisrow{y}
	] {\data};
	%\addlegendentry{$(\ell,\ell,0)$};
	
	\pgfplotstableread{data/G1F3/Vini.txt} \data;
	\addplot [
		color=ratecolor2,
		dotted,
		line width=1pt,
		line cap=round,
		dash pattern=on 0pt off 2\pgflinewidth,
		mark=square*,
		mark size = 1.4,
		mark options={solid,ratecolor2}
	]
	table[
		x expr=\thisrow{x},
		y expr=\thisrow{y}
	] {\data};
	%\addlegendentry{$(\ell,0,\ell)$};
	
	\pgfplotstableread{data/G1F3/Vnii.txt} \data;
	\addplot [
		color=ratecolor2,
		dotted,
		line width=1pt,
		line cap=round,
		dash pattern=on 0pt off 2\pgflinewidth,
		mark=diamond*,
		mark options={solid,ratecolor2}
	]
	table[
		x expr=\thisrow{x},
		y expr=\thisrow{y}
	] {\data};
	%\addlegendentry{$(0,\ell,\ell)$};

%
% TRIPPLE
%

	\pgfplotstableread{data/G1F3/Viii.txt} \data;
	\addplot [
		color=ratecolor3,
		dash pattern=on 2pt off 1pt on 2pt off 2pt,
		line width=1pt,
		line cap=round,
		mark=*,
		mark size = 1.4,
		mark options={solid,ratecolor3}
	]
	table[
		x expr=\thisrow{x},
		y expr=\thisrow{y}
	] {\data};
	%\addlegendentry{$(\ell,\ell,\ell)$};
	
		\pgfplotsset{
    	after end axis/.code={
        	\node[above right] at (0.03,0.03){\texttt{G1F3}};
    	}
	}
	
\end{axis}
\end{tikzpicture}

%% file: figures/legend_times_G1.tikz
\begin{tikzpicture}
    \begin{customlegend}[
    	legend columns=3,
		legend style={column sep=1ex,font=\scriptsize},
		legend cell align=left,
		legend entries={MLMC\\MIMC, FT\\MIMC, TD\\MLQMC\\MIQMC, FT\\MIQMC, TD\\}
		]
		\addlegendimage{ratecolor1,solid,line width=1pt,mark=*,mark size=1.4pt,mark options={solid},line cap=round}
		\addlegendimage{ratecolor3,solid,line width=1pt,mark=square*,mark size=1.4pt,mark options={solid},line cap=round}
		\addlegendimage{ratecolor2,solid,line width=1pt,mark=diamond*,mark options={solid},line cap=round}
   	\addlegendimage{ratecolor1,dashed,line cap=round,line width=1pt,mark=*,mark size=1.4pt,mark options={solid,fill=white}}
   	\addlegendimage{ratecolor3,dashed,line cap=round,line width=1pt,mark=square*,mark size=1.4pt,mark options={solid,fill=white}}
   	\addlegendimage{ratecolor2,dashed,line cap=round,line width=1pt,mark=diamond*,mark options={solid,fill=white}}
   \end{customlegend}
\end{tikzpicture}

%% file: figures/cost_G1F1.tikz
	\begin{tikzpicture}[trim axis left, trim axis right]
	\begin{axis}[
		width=\figurewidth,
		height=\figureheight,
		xmode=log,
		xmin=1e-05,
		xmax=0.001,
		xminorticks=true,
		xlabel={\small requested tolerance $\epsilon$},
		xmajorgrids,
		xminorgrids,
		every x tick label/.append style={font=\scriptsize},
		ymode=log,
		ymin=10000, 
		ymax=1000000000,
		yminorticks=true,
		ylabel={\small standard cost},
		ytick={100,1000,10000,100000,1000000,10000000,100000000,1000000000},
		ymajorgrids,
		yminorgrids,
		every y tick label/.append style={font=\scriptsize},
	]
	
	\addplot [
		color=ratecolor1,
		solid,
		line width=1pt,
		mark=*,
		mark size = 1.4,
		mark options={solid,ratecolor1}
	]
 	table[
		x = x,
		y = y
	] {data/G1F1MLMC/cost_formatted.txt};
	%\addlegendentry{MLMC};

	\addplot [
		color=ratecolor1,
		dashed,
		line cap=round,
		line width=1pt,
		mark=*,
		mark size = 1.4,
		mark options={solid,fill=white}
	]
 	table[
		x = x,
		y = y
	] {data/G1F1MLQM/cost_formatted.txt};
	%\addlegendentry{MLQMC};
	
\addplot [		
		color=ratecolor3,
		solid,
		mark=square*,
		mark size=1.4pt,
		line width = 1pt,
    	mark options={solid,ratecolor3}
	]
 	table[
		x = x,
		y = y
	] {data/G1F1FTMC/cost_formatted.txt};
	%\addlegendentry{MIMC, FT};
	
	\addplot [
		color=ratecolor3,
		dashed,
		line cap=round,
		mark=square*,
		mark size=1.4pt,
		line width = 1pt,
    	mark options={solid,fill=white}
	]
 	table[
		x = x,
		y = y
	] {data/G1F1FTQM/cost_formatted.txt};
	%\addlegendentry{MIQMC, FT};

	\addplot [
		color=ratecolor2,
		solid,
		mark=diamond*,
		line width=1pt,
		mark options={solid,ratecolor2}
	]
 	table[
		x = x,
		y = y
	] {data/G1F1TDMC/cost_formatted.txt};
	%\addlegendentry{MIMC, TD};

	\addplot [
		color=ratecolor2,
		dashed,
		line cap=round,
		mark=diamond*,
		line width=1pt,
		mark options={solid,fill=white}
	]
 	table[
		x = x,
		y = y
	] {data/G1F1TDQM/cost_formatted.txt};
	%\addlegendentry{MIQMC, TD};
	
	\LogLogSlopeTriangle{0.9}{0.2}{0.8}{2}{black};
	\LogLogSlopeTriangle{0.9}{0.2}{0.6}{1}{black};
	
	\pgfplotsset{
    	after end axis/.code={
        	\node[above right] at (axis cs:1e-5,10000){\texttt{G1F1}};
    	}
	}
	
	\end{axis}
\end{tikzpicture}

%% file: figures/times_G1F1.tikz
	\begin{tikzpicture}[trim axis left, trim axis right]
	\begin{axis}[
		width=\figurewidth,
		height=\figureheight,
		xmode=log,
		xmin=1e-05,
		xmax=0.001,
		xminorticks=true,
		xlabel={\small requested tolerance $\epsilon$},
		xmajorgrids,
		xminorgrids,
		every x tick label/.append style={font=\scriptsize},
		ymode=log,
		ymin=10,
		ymax=1000000,
		yminorticks=true,
		ylabel={\small runtime [seconds]},
		ytick={10,100,1000,10000,100000,1000000},
		ymajorgrids,
		yminorgrids,
		every y tick label/.append style={font=\scriptsize},
	]

	\addplot [
		color=ratecolor1,
		solid,
		line width=1pt,
		mark=*,
		mark size = 1.4,
		mark options={solid,ratecolor1}
	]
 	table[
		x = x,
		y = y
	] {data/G1F1MLMC/times_formatted.txt};
	%\addlegendentry{MLMC};

	\addplot [
		color=ratecolor1,
		dashed,
		line cap=round,
		line width=1pt,
		mark=*,
		mark size = 1.4,
		mark options={solid,fill=white}
	]
 	table[
		x = x,
		y = y
	] {data/G1F1MLQM/times_formatted.txt};
	%\addlegendentry{MLQMC};
	
	\addplot [		
		color=ratecolor3,
		solid,
		mark=square*,
		mark size=1.4pt,
		line width = 1pt,
    	mark options={solid,ratecolor3}
	]
 	table[
		x = x,
		y = y
	] {data/G1F1FTMC/times_formatted.txt};
	%\addlegendentry{MIMC, FT};
	
	\addplot [
		color=ratecolor3,
		dashed,
		line cap=round,
		mark=square*,
		mark size=1.4pt,
		line width = 1pt,
    	mark options={solid,fill=white}
	]
 	table[
		x = x,
		y = y
	] {data/G1F1FTQM/times_formatted.txt};
	%\addlegendentry{MIQMC, FT};
	
	\addplot [
		color=ratecolor2,
		solid,
		mark=diamond*,
		line width=1pt,
		mark options={solid,ratecolor2}
	]
 	table[
		x = x,
		y = y
	] {data/G1F1TDMC/times_formatted.txt};
	%\addlegendentry{MIMC, TD};

	\addplot [
		color=ratecolor2,
		dashed,
		line cap=round,
		mark=diamond*,
		line width=1pt,
		mark options={solid,fill=white}
	]
 	table[
		x = x,
		y = y
	] {data/G1F1TDQM/times_formatted.txt};
	%\addlegendentry{MIQMC, TD};
	
	\LogLogSlopeTriangle{0.9}{0.2}{0.8}{2}{black};
	\LogLogSlopeTriangle{0.9}{0.2}{0.6}{1}{black};

		\pgfplotsset{
    	after end axis/.code={
        	\node[above right] at (axis cs:1e-5,10){\texttt{G1F1}};
    	}
	}
	
	\end{axis}
\end{tikzpicture}

%% file: figures/cost_G1F2.tikz
	\begin{tikzpicture}[trim axis left, trim axis right]
	\begin{axis}[
		width=\figurewidth,
		height=\figureheight,
		xmode=log,
		xmin=1e-05,
		xmax=0.001,
		xminorticks=true,
		xlabel={\small requested tolerance $\epsilon$},
		xmajorgrids,
		xminorgrids,
		every x tick label/.append style={font=\scriptsize},
		ymode=log,
		ymin=10000, 
		ymax=100000000,
		yminorticks=true,
		ylabel={\small standard cost},
		ytick={100,1000,10000,100000,1000000,10000000,100000000,1000000000},
		ymajorgrids,
		yminorgrids,
		every y tick label/.append style={font=\scriptsize},
	]
	
	\addplot [
		color=ratecolor1,
		solid,
		line width=1pt,
		mark=*,
		mark size = 1.4,
		mark options={solid,ratecolor1}
	]
 	table[
		x = x,
		y = y
	] {data/G1F2MLMC/cost_formatted.txt};
	%\addlegendentry{MLMC};

	\addplot [
		color=ratecolor1,
		dashed,
		line cap=round,
		line width=1pt,
		mark=*,
		mark size = 1.4,
		mark options={solid,fill=white}
	]
 	table[
		x = x,
		y = y
	] {data/G1F2MLQM/cost_formatted.txt};
	%\addlegendentry{MLQMC};
	
\addplot [		
		color=ratecolor3,
		solid,
		mark=square*,
		mark size=1.4pt,
		line width = 1pt,
    	mark options={solid,ratecolor3}
	]
 	table[
		x = x,
		y = y
	] {data/G1F2FTMC/cost_formatted.txt};
	%\addlegendentry{MIMC, FT};
	
	\addplot [
		color=ratecolor3,
		dashed,
		line cap=round,
		mark=square*,
		mark size=1.4pt,
		line width = 1pt,
    	mark options={solid,fill=white}
	]
 	table[
		x = x,
		y = y
	] {data/G1F2FTQM/cost_formatted.txt};
	%\addlegendentry{MIQMC, FT};

	\addplot [
		color=ratecolor2,
		solid,
		mark=diamond*,
		line width=1pt,
		mark options={solid,ratecolor2}
	]
 	table[
		x = x,
		y = y
	] {data/G1F2TDMC/cost_formatted.txt};
	%\addlegendentry{MIMC, TD};

	\addplot [
		color=ratecolor2,
		dashed,
		line cap=round,
		mark=diamond*,
		line width=1pt,
		mark options={solid,fill=white}
	]
 	table[
		x = x,
		y = y
	] {data/G1F2TDQM/cost_formatted.txt};
	%\addlegendentry{MIQMC, TD};
	
	\LogLogSlopeTriangle{0.9}{0.2}{0.7}{2}{black};
	\LogLogSlopeTriangle{0.9}{0.2}{0.5}{1}{black};
	
	\pgfplotsset{
    	after end axis/.code={
        	\node[above right] at (axis cs:1e-5,10000){\texttt{G1F2}};
    	}
	}
	
	\end{axis}
\end{tikzpicture}

%% file: figures/times_G1F2.tikz
	\begin{tikzpicture}[trim axis left, trim axis right]
	\begin{axis}[
		width=\figurewidth,
		height=\figureheight,
		xmode=log,
		xmin=1e-05,
		xmax=0.001,
		xminorticks=true,
		xlabel={\small requested tolerance $\epsilon$},
		xmajorgrids,
		xminorgrids,
		every x tick label/.append style={font=\scriptsize},
		ymode=log,
		ymin=10,
		ymax=100000,
		yminorticks=true,
		ylabel={\small runtime [seconds]},
		ytick={10,100,1000,10000,100000,1000000},
		ymajorgrids,
		yminorgrids,
		every y tick label/.append style={font=\scriptsize},
	]

	\addplot [
		color=ratecolor1,
		solid,
		line width=1pt,
		mark=*,
		mark size = 1.4,
		mark options={solid,ratecolor1}
	]
 	table[
		x = x,
		y = y
	] {data/G1F2MLMC/times_formatted.txt};
	%\addlegendentry{MLMC};

	\addplot [
		color=ratecolor1,
		dashed,
		line cap=round,
		line width=1pt,
		mark=*,
		mark size = 1.4,
		mark options={solid,fill=white}
	]
 	table[
		x = x,
		y = y
	] {data/G1F2MLQM/times_formatted.txt};
	%\addlegendentry{MLQMC};
	
\addplot [		
		color=ratecolor3,
		solid,
		mark=square*,
		mark size=1.4pt,
		line width = 1pt,
    	mark options={solid,ratecolor3}
	]
 	table[
		x = x,
		y = y
	] {data/G1F2FTMC/times_formatted.txt};
	%\addlegendentry{MIMC, FT};
	
	\addplot [
		color=ratecolor3,
		dashed,
		line cap=round,
		mark=square*,
		mark size=1.4pt,
		line width = 1pt,
    	mark options={solid,fill=white}
	]
 	table[
		x = x,
		y = y
	] {data/G1F2FTQM/times_formatted.txt};
	%\addlegendentry{MIQMC, FT};

	\addplot [
		color=ratecolor2,
		solid,
		mark=diamond*,
		line width=1pt,
		mark options={solid,ratecolor2}
	]
 	table[
		x = x,
		y = y
	] {data/G1F2TDMC/times_formatted.txt};
	%\addlegendentry{MIMC};

	\addplot [
		color=ratecolor2,
		dashed,
		line cap=round,
		mark=diamond*,
		line width=1pt,
		mark options={solid,fill=white}
	]
 	table[
		x = x,
		y = y
	] {data/G1F2TDQM/times_formatted.txt};
	%\addlegendentry{MIQMC};
	
	\LogLogSlopeTriangle{0.9}{0.2}{0.7}{2}{black};
	\LogLogSlopeTriangle{0.9}{0.2}{0.5}{1}{black};

		\pgfplotsset{
    	after end axis/.code={
        	\node[above right] at (axis cs:1e-5,10){\texttt{G1F2}};
    	}
	}

	\end{axis}
\end{tikzpicture}

%% file: figures/cost_G1F3.tikz
	\begin{tikzpicture}[trim axis left, trim axis right]
	\begin{axis}[
		width=\figurewidth,
		height=\figureheight,
		xmode=log,
		xmin=1e-03,
		xmax=0.1,
		xminorticks=true,
		xlabel={\small requested tolerance $\epsilon$},
		xmajorgrids,
		xminorgrids,
		every x tick label/.append style={font=\scriptsize},
		ymode=log,
		ymin=1000, 
		ymax=10000000,
		yminorticks=true,
		ylabel={\small standard cost},
		ytick={100,1000,10000,100000,1000000,10000000},
		xticklabels={$10^{-6}$,$10^{-5}$,$10^{-4}$,$10^{-3}$},
		ymajorgrids,
		yminorgrids,
		every y tick label/.append style={font=\scriptsize},
	]

	\addplot [
		color=ratecolor1,
		solid,
		line width=1pt,
		mark=*,
		mark size = 1.4,
		mark options={solid,ratecolor1}
	]
 	table[
		x = x,
		y = y
	] {data/G1F3MLMC/cost_formatted.txt};
	%\addlegendentry{MLMC};

	\addplot [
		color=ratecolor1,
		dashed,
		line cap=round,
		line width=1pt,
		mark=*,
		mark size = 1.4,
		mark options={solid,fill=white}
	]
 	table[
		x = x,
		y = y
	] {data/G1F3MLQM/cost_formatted.txt};
	%\addlegendentry{MLQMC};
	
\addplot [		
		color=ratecolor3,
		solid,
		mark=square*,
		mark size=1.4pt,
		line width = 1pt,
    	mark options={solid,ratecolor3}
	]
 	table[
		x = x,
		y = y
	] {data/G1F3FTMC/cost_formatted.txt};
	%\addlegendentry{MIMC, FT};
	
	\addplot [
		color=ratecolor3,
		dashed,
		line cap=round,
		mark=square*,
		mark size=1.4pt,
		line width = 1pt,
    	mark options={solid,fill=white}
	]
 	table[
		x = x,
		y = y
	] {data/G1F3FTQM/cost_formatted.txt};
	%\addlegendentry{MIQMC, FT};

	\addplot [
		color=ratecolor2,
		solid,
		mark=diamond*,
		line width=1pt,
		mark options={solid,ratecolor2}
	]
 	table[
		x = x,
		y = y
	] {data/G1F3TDMC/cost_formatted.txt};
	%\addlegendentry{MIMC, TD};

	\addplot [
		color=ratecolor2,
		dashed,
		line cap=round,
		mark=diamond*,
		line width=1pt,
		mark options={solid,fill=white}
	]
 	table[
		x = x,
		y = y
	] {data/G1F3TDQM/cost_formatted.txt};
	%\addlegendentry{MIQMC, TD};
	
	%\LogLogSlopeTriangle{0.3}{0.2}{0.4}{2}{black};
	\LogLogSlopeTriangle{0.3}{0.2}{0.2}{2}{black};
	
	\pgfplotsset{
    	after end axis/.code={
        	\node[above right] at (axis cs:1e-3,1000){\texttt{G1F3}};
    	}
	}
	
	\end{axis}
\end{tikzpicture}

%% file: figures/times_G1F3.tikz
	\begin{tikzpicture}[trim axis left, trim axis right]
	\begin{axis}[
		width=\figurewidth,
		height=\figureheight,
		xmode=log,
		xmin=1e-05,
		xmax=0.001,
		xminorticks=true,
		xlabel={\small requested tolerance $\epsilon$},
		xmajorgrids,
		xminorgrids,
		every x tick label/.append style={font=\scriptsize},
		ymode=log,
		ymin=10,
		ymax=100000,
		yminorticks=true,
		ylabel={\small runtime [seconds]},
		ytick={10,100,1000,10000,100000,1000000},
		ymajorgrids,
		yminorgrids,
		every y tick label/.append style={font=\scriptsize},
	]

	\addplot [
		color=ratecolor1,
		solid,
		line width=1pt,
		mark=*,
		mark size = 1.4,
		mark options={solid,ratecolor1}
	]
 	table[
		x = x,
		y = y
	] {data/G1F3MLMC/times_formatted.txt};
	%\addlegendentry{MLMC};

	\addplot [
		color=ratecolor1,
		dashed,
		line cap=round,
		line width=1pt,
		mark=*,
		mark size = 1.4,
		mark options={solid,fill=white}
	]
 	table[
		x = x,
		y = y
	] {data/G1F3MLQM/times_formatted.txt};
	%\addlegendentry{MLQMC};
	
\addplot [		
		color=ratecolor3,
		solid,
		mark=square*,
		mark size=1.4pt,
		line width = 1pt,
    	mark options={solid,ratecolor3}
	]
 	table[
		x = x,
		y = y
	] {data/G1F3FTMC/times_formatted.txt};
	%\addlegendentry{MIMC, FT};
	
	\addplot [
		color=ratecolor3,
		dashed,
		line cap=round,
		mark=square*,
		mark size=1.4pt,
		line width = 1pt,
    	mark options={solid,fill=white}
	]
 	table[
		x = x,
		y = y
	] {data/G1F3FTQM/times_formatted.txt};
	%\addlegendentry{MIQMC, FT};

	\addplot [
		color=ratecolor2,
		solid,
		mark=diamond*,
		line width=1pt,
		mark options={solid,ratecolor2}
	]
 	table[
		x = x,
		y = y
	] {data/G1F3TDMC/times_formatted.txt};
	%\addlegendentry{MIMC};

	\addplot [
		color=ratecolor2,
		dashed,
		line cap=round,
		mark=diamond*,
		line width=1pt,
		mark options={solid,fill=white}
	]
 	table[
		x = x,
		y = y
	] {data/G1F3TDQM/times_formatted.txt};
	%\addlegendentry{MIQMC};
	
	%\LogLogSlopeTriangle{0.3}{0.2}{0.4}{2}{black};
	\LogLogSlopeTriangle{0.3}{0.2}{0.2}{2}{black};

		\pgfplotsset{
    	after end axis/.code={
        	\node[above right] at (axis cs:1e-5,10){\texttt{G1F3}};
    	}
	}
	
	\end{axis}
\end{tikzpicture}

%% file: figures/G2F1_analyse_E.tikz
\begin{tikzpicture}
\def\maxL{4}
\begin{axis}[
	height=\figureheight,
	width=\figurewidth,
	xmin=0,
	xmax=\maxL,
	xminorticks=true,
	xlabel={\small $\ell$},
	every x tick label/.append style={font=\scriptsize},
	ymin=-25,
	ymax=5,
	yminorticks=true,
	ylabel={\small $\log_2(|\E[\Delta G_\bell]|)$},
	every y tick label/.append style={font=\scriptsize},
	legend style ={ at={(1.03,1)}, 
        anchor=north west, draw=black, style={font=\scriptsize},
        fill=white,align=left}
]

%
% ZERO
%

	\pgfplotstableread{data/G2F1/Ennn.txt} \data;
	\addplot [
		color=black,
		solid,
		line width=1pt,
		line cap=round,
		mark=*,
		mark size = 1.4,
		mark options={solid}
	]
	table[
		x expr=\thisrow{x},
		y expr=\thisrow{y}
	] {\data};
	%\addlegendentry{$(0,0,0)$};

%
% SINGLES
%

	\pgfplotstableread{data/G2F1/Einn.txt} \data;
	\addplot [
		color=ratecolor1,
		dashed,
		line width=1pt,
		line cap=round,
		mark=*,
		mark size = 1.4,
		mark options={solid,ratecolor1}
	]
	table[
		x expr=\thisrow{x},
		y expr=\thisrow{y}
	] {\data};
	%\addlegendentry{$(\ell,0,0)$};
	
	\pgfplotstableread{data/G2F1/Enin.txt} \data;
	\addplot [
		color=ratecolor1,
		dashed,
		line width=1pt,
		line cap=round,
		mark=square*,
		mark size = 1.4,
		mark options={solid,ratecolor1}
	]
	table[
		x expr=\thisrow{x},
		y expr=\thisrow{y}
	] {\data};
	%\addlegendentry{$(0,\ell,0)$};
	
	\pgfplotstableread{data/G2F1/Enni.txt} \data;
	\addplot [
		color=ratecolor1,
		dashed,
		line width=1pt,
		line cap=round,
		mark=diamond*,
		mark options={solid,ratecolor1}
	]
	table[
		x expr=\thisrow{x},
		y expr=\thisrow{y}
	] {\data};
	%\addlegendentry{$(0,0,\ell)$};
	
%
% DOUBLES
%

	\pgfplotstableread{data/G2F1/Eiin.txt} \data;
	\addplot [
		color=ratecolor2,
		dotted,
		line width=1pt,
		line cap=round,
		dash pattern=on 0pt off 2\pgflinewidth,
		mark=*,
		mark size = 1.4,
		mark options={solid,ratecolor2}
	]
	table[
		x expr=\thisrow{x},
		y expr=\thisrow{y}
	] {\data};
	%\addlegendentry{$(\ell,\ell,0)$};
	
	\pgfplotstableread{data/G2F1/Eini.txt} \data;
	\addplot [
		color=ratecolor2,
		dotted,
		line width=1pt,
		line cap=round,
		dash pattern=on 0pt off 2\pgflinewidth,
		mark=square*,
		mark size = 1.4,
		mark options={solid,ratecolor2}
	]
	table[
		x expr=\thisrow{x},
		y expr=\thisrow{y}
	] {\data};
	%\addlegendentry{$(\ell,0,\ell)$};
	
	\pgfplotstableread{data/G2F1/Enii.txt} \data;
	\addplot [
		color=ratecolor2,
		dotted,
		line width=1pt,
		line cap=round,
		dash pattern=on 0pt off 2\pgflinewidth,
		mark=diamond*,
		mark options={solid,ratecolor2}
	]
	table[
		x expr=\thisrow{x},
		y expr=\thisrow{y}
	] {\data};
	%\addlegendentry{$(0,\ell,\ell)$};

%
% TRIPPLE
%

	\pgfplotstableread{data/G2F1/Eiii.txt} \data;
	\addplot [
		color=ratecolor3,
		dash pattern=on 2pt off 1pt on 2pt off 2pt,
		line width=1pt,
		line cap=round,
		mark=*,
		mark size = 1.4,
		mark options={solid,ratecolor3}
	]
	table[
		x expr=\thisrow{x},
		y expr=\thisrow{y}
	] {\data};
	%\addlegendentry{$(\ell,\ell,\ell)$};
	
			\pgfplotsset{
    	after end axis/.code={
        	\node[above right] at (0.03,0.03){\texttt{G2F1}};
    	}
	}
	
\end{axis}
\end{tikzpicture}

%% file: figures/G2F1_analyse_V.tikz
\begin{tikzpicture}
\def\maxL{4}
\begin{axis}[
	height=\figureheight,
	width=\figurewidth,
	xmin=0,
	xmax=\maxL,
	xminorticks=true,
	xlabel={\small $\ell$},
	every x tick label/.append style={font=\scriptsize},
	ymin=-50,
	ymax=10,
	yminorticks=true,
	ylabel={\small $\log_2(\V[\Delta G_\bell])$},
	every y tick label/.append style={font=\scriptsize},
	%legend style ={ at={(1.06,1)}, 
    %    anchor=north west, draw=black, style={font=\scriptsize},
     %   fill=white,align=left},
   ymajorgrids,
	xmajorgrids
]

%
% ZERO
%

	\pgfplotstableread{data/G2F1/Vnnn.txt} \data;
	\addplot [
		color=black,
		solid,
		line width=1pt,
		line cap=round,
		mark=*,
		mark size = 1.4,
		mark options={solid},
		forget plot
	]
	table[
		x expr=\thisrow{x},
		y expr=\thisrow{y}
	] {\data};
	%\addlegendentry{$(0,0,0)$};

%
% SINGLES
%

	\pgfplotstableread{data/G2F1/Vinn.txt} \data;
	\addplot [
		color=ratecolor1,
		dashed,
		line width=1pt,
		line cap=round,
		mark=*,
		mark size = 1.4,
		mark options={solid,ratecolor1}
	]
	table[
		x expr=\thisrow{x},
		y expr=\thisrow{y}
	] {\data};
	%\addlegendentry{$(\ell,0,0)$};
	
	\pgfplotstableread{data/G2F1/Vnin.txt} \data;
	\addplot [
		color=ratecolor1,
		dashed,
		line width=1pt,
		line cap=round,
		mark=square*,
		mark size = 1.4,
		mark options={solid,ratecolor1}
	]
	table[
		x expr=\thisrow{x},
		y expr=\thisrow{y}
	] {\data};
	%\addlegendentry{$(0,\ell,0)$};
	
	\pgfplotstableread{data/G2F1/Vnni.txt} \data;
	\addplot [
		color=ratecolor1,
		dashed,
		line width=1pt,
		line cap=round,
		mark=diamond*,
		mark options={solid,ratecolor1}
	]
	table[
		x expr=\thisrow{x},
		y expr=\thisrow{y}
	] {\data};
	%\addlegendentry{$(0,0,\ell)$};
	
%
% DOUBLES
%

	\pgfplotstableread{data/G2F1/Viin.txt} \data;
	\addplot [
		color=ratecolor2,
		dotted,
		line width=1pt,
		line cap=round,
		dash pattern=on 0pt off 2\pgflinewidth,
		mark=*,
		mark size = 1.4,
		mark options={solid,ratecolor2}
	]
	table[
		x expr=\thisrow{x},
		y expr=\thisrow{y}
	] {\data};
	%\addlegendentry{$(\ell,\ell,0)$};
	
	\pgfplotstableread{data/G2F1/Vini.txt} \data;
	\addplot [
		color=ratecolor2,
		dotted,
		line width=1pt,
		line cap=round,
		dash pattern=on 0pt off 2\pgflinewidth,
		mark=square*,
		mark size = 1.4,
		mark options={solid,ratecolor2}
	]
	table[
		x expr=\thisrow{x},
		y expr=\thisrow{y}
	] {\data};
	%\addlegendentry{$(\ell,0,\ell)$};
	
	\pgfplotstableread{data/G2F1/Vnii.txt} \data;
	\addplot [
		color=ratecolor2,
		dotted,
		line width=1pt,
		line cap=round,
		dash pattern=on 0pt off 2\pgflinewidth,
		mark=diamond*,
		mark options={solid,ratecolor2}
	]
	table[
		x expr=\thisrow{x},
		y expr=\thisrow{y}
	] {\data};
	%\addlegendentry{$(0,\ell,\ell)$};

%
% TRIPPLE
%

	\pgfplotstableread{data/G2F1/Viii.txt} \data;
	\addplot [
		color=ratecolor3,
		dash pattern=on 2pt off 1pt on 2pt off 2pt,
		line width=1pt,
		line cap=round,
		mark=*,
		mark size = 1.4,
		mark options={solid,ratecolor3}
	]
	table[
		x expr=\thisrow{x},
		y expr=\thisrow{y}
	] {\data};
	%\addlegendentry{$(\ell,\ell,\ell)$};
	
			\pgfplotsset{
    	after end axis/.code={
        	\node[above right] at (0.03,0.03){\texttt{G2F1}};
    	}
	}
	
\end{axis}
\end{tikzpicture}

%% file: figures/G2F2_analyse_E.tikz
\begin{tikzpicture}
\def\maxL{4}
\begin{axis}[
	height=\figureheight,
	width=\figurewidth,
	xmin=0,
	xmax=\maxL,
	xminorticks=true,
	xlabel={\small $\ell$},
	every x tick label/.append style={font=\scriptsize},
	ymin=-25,
	ymax=5,
	yminorticks=true,
	ylabel={\small $\log_2(|\E[\Delta G_\bell]|)$},
	every y tick label/.append style={font=\scriptsize},
	legend style ={ at={(1.03,1)}, 
        anchor=north west, draw=black, style={font=\scriptsize},
        fill=white,align=left},
   ymajorgrids,
	xmajorgrids
]

%
% ZERO
%

	\pgfplotstableread{data/G2F2/Ennn.txt} \data;
	\addplot [
		color=black,
		solid,
		line width=1pt,
		line cap=round,
		mark=*,
		mark size = 1.4,
		mark options={solid}
	]
	table[
		x expr=\thisrow{x},
		y expr=\thisrow{y}
	] {\data};
	%\addlegendentry{$(0,0,0)$};

%
% SINGLES
%

	\pgfplotstableread{data/G2F2/Einn.txt} \data;
	\addplot [
		color=ratecolor1,
		dashed,
		line width=1pt,
		line cap=round,
		mark=*,
		mark size = 1.4,
		mark options={solid,ratecolor1}
	]
	table[
		x expr=\thisrow{x},
		y expr=\thisrow{y}
	] {\data};
	%\addlegendentry{$(\ell,0,0)$};
	
	\pgfplotstableread{data/G2F2/Enin.txt} \data;
	\addplot [
		color=ratecolor1,
		dashed,
		line width=1pt,
		line cap=round,
		mark=square*,
		mark size = 1.4,
		mark options={solid,ratecolor1}
	]
	table[
		x expr=\thisrow{x},
		y expr=\thisrow{y}
	] {\data};
	%\addlegendentry{$(0,\ell,0)$};
	
	\pgfplotstableread{data/G2F2/Enni.txt} \data;
	\addplot [
		color=ratecolor1,
		dashed,
		line width=1pt,
		line cap=round,
		mark=diamond*,
		mark options={solid,ratecolor1}
	]
	table[
		x expr=\thisrow{x},
		y expr=\thisrow{y}
	] {\data};
	%\addlegendentry{$(0,0,\ell)$};
	
%
% DOUBLES
%

	\pgfplotstableread{data/G2F2/Eiin.txt} \data;
	\addplot [
		color=ratecolor2,
		dotted,
		line width=1pt,
		line cap=round,
		dash pattern=on 0pt off 2\pgflinewidth,
		mark=*,
		mark size = 1.4,
		mark options={solid,ratecolor2}
	]
	table[
		x expr=\thisrow{x},
		y expr=\thisrow{y}
	] {\data};
	%\addlegendentry{$(\ell,\ell,0)$};
	
	\pgfplotstableread{data/G2F2/Eini.txt} \data;
	\addplot [
		color=ratecolor2,
		dotted,
		line width=1pt,
		line cap=round,
		dash pattern=on 0pt off 2\pgflinewidth,
		mark=square*,
		mark size = 1.4,
		mark options={solid,ratecolor2}
	]
	table[
		x expr=\thisrow{x},
		y expr=\thisrow{y}
	] {\data};
	%\addlegendentry{$(\ell,0,\ell)$};
	
	\pgfplotstableread{data/G2F2/Enii.txt} \data;
	\addplot [
		color=ratecolor2,
		dotted,
		line width=1pt,
		line cap=round,
		dash pattern=on 0pt off 2\pgflinewidth,
		mark=diamond*,
		mark options={solid,ratecolor2}
	]
	table[
		x expr=\thisrow{x},
		y expr=\thisrow{y}
	] {\data};
	%\addlegendentry{$(0,\ell,\ell)$};

%
% TRIPPLE
%

	\pgfplotstableread{data/G2F2/Eiii.txt} \data;
	\addplot [
		color=ratecolor3,
		dash pattern=on 2pt off 1pt on 2pt off 2pt,
		line width=1pt,
		line cap=round,
		mark=*,
		mark size = 1.4,
		mark options={solid,ratecolor3}
	]
	table[
		x expr=\thisrow{x},
		y expr=\thisrow{y}
	] {\data};
	%\addlegendentry{$(\ell,\ell,\ell)$};
	
			\pgfplotsset{
    	after end axis/.code={
        	\node[above right] at (0.03,0.03){\texttt{G2F2}};
    	}
	}
	
\end{axis}
\end{tikzpicture}

%% file: figures/G2F2_analyse_V.tikz
\begin{tikzpicture}
\def\maxL{4}
\begin{axis}[
	height=\figureheight,
	width=\figurewidth,
	xmin=0,
	xmax=\maxL,
	xminorticks=true,
	xlabel={\small $\ell$},
	every x tick label/.append style={font=\scriptsize},
	ymin=-50,
	ymax=10,
	yminorticks=true,
	ylabel={\small $\log_2(\V[\Delta G_\bell])$},
	every y tick label/.append style={font=\scriptsize},
	legend style ={ at={(1.06,1)}, 
        anchor=north west, draw=black, style={font=\scriptsize},
        fill=white,align=left},
   ymajorgrids,
	xmajorgrids
]

%
% ZERO
%

	\pgfplotstableread{data/G2F2/Vnnn.txt} \data;
	\addplot [
		color=black,
		solid,
		line width=1pt,
		line cap=round,
		mark=*,
		mark size = 1.4,
		mark options={solid},
		forget plot
	]
	table[
		x expr=\thisrow{x},
		y expr=\thisrow{y}
	] {\data};
	%\addlegendentry{$(0,0,0)$};

%
% SINGLES
%

	\pgfplotstableread{data/G2F2/Vinn.txt} \data;
	\addplot [
		color=ratecolor1,
		dashed,
		line width=1pt,
		line cap=round,
		mark=*,
		mark size = 1.4,
		mark options={solid,ratecolor1}
	]
	table[
		x expr=\thisrow{x},
		y expr=\thisrow{y}
	] {\data};
	%\addlegendentry{$(\ell,0,0)$};
	
	\pgfplotstableread{data/G2F2/Vnin.txt} \data;
	\addplot [
		color=ratecolor1,
		dashed,
		line width=1pt,
		line cap=round,
		mark=square*,
		mark size = 1.4,
		mark options={solid,ratecolor1}
	]
	table[
		x expr=\thisrow{x},
		y expr=\thisrow{y}
	] {\data};
	%\addlegendentry{$(0,\ell,0)$};
	
	\pgfplotstableread{data/G2F2/Vnni.txt} \data;
	\addplot [
		color=ratecolor1,
		dashed,
		line width=1pt,
		line cap=round,
		mark=diamond*,
		mark options={solid,ratecolor1}
	]
	table[
		x expr=\thisrow{x},
		y expr=\thisrow{y}
	] {\data};
	%\addlegendentry{$(0,0,\ell)$};
	
%
% DOUBLES
%

	\pgfplotstableread{data/G2F2/Viin.txt} \data;
	\addplot [
		color=ratecolor2,
		dotted,
		line width=1pt,
		line cap=round,
		dash pattern=on 0pt off 2\pgflinewidth,
		mark=*,
		mark size = 1.4,
		mark options={solid,ratecolor2}
	]
	table[
		x expr=\thisrow{x},
		y expr=\thisrow{y}
	] {\data};
	%\addlegendentry{$(\ell,\ell,0)$};
	
	\pgfplotstableread{data/G2F2/Vini.txt} \data;
	\addplot [
		color=ratecolor2,
		dotted,
		line width=1pt,
		line cap=round,
		dash pattern=on 0pt off 2\pgflinewidth,
		mark=square*,
		mark size = 1.4,
		mark options={solid,ratecolor2}
	]
	table[
		x expr=\thisrow{x},
		y expr=\thisrow{y}
	] {\data};
	%\addlegendentry{$(\ell,0,\ell)$};
	
	\pgfplotstableread{data/G2F2/Vnii.txt} \data;
	\addplot [
		color=ratecolor2,
		dotted,
		line width=1pt,
		line cap=round,
		dash pattern=on 0pt off 2\pgflinewidth,
		mark=diamond*,
		mark options={solid,ratecolor2}
	]
	table[
		x expr=\thisrow{x},
		y expr=\thisrow{y}
	] {\data};
	%\addlegendentry{$(0,\ell,\ell)$};

%
% TRIPPLE
%

	\pgfplotstableread{data/G2F2/Viii.txt} \data;
	\addplot [
		color=ratecolor3,
		dash pattern=on 2pt off 1pt on 2pt off 2pt,
		line width=1pt,
		line cap=round,
		mark=*,
		mark size = 1.4,
		mark options={solid,ratecolor3}
	]
	table[
		x expr=\thisrow{x},
		y expr=\thisrow{y}
	] {\data};
	%\addlegendentry{$(\ell,\ell,\ell)$};
	
			\pgfplotsset{
    	after end axis/.code={
        	\node[above right] at (0.03,0.03){\texttt{G2F2}};
    	}
	}
	
\end{axis}
\end{tikzpicture}

%% file: figures/G2F3_analyse_E.tikz
\begin{tikzpicture}
\def\maxL{4}
\begin{axis}[
	height=\figureheight,
	width=\figurewidth,
	xmin=0,
	xmax=\maxL,
	xminorticks=true,
	xlabel={\small $\ell$},
	every x tick label/.append style={font=\scriptsize},
	ymin=-25,
	ymax=5,
	yminorticks=true,
	ylabel={\small $\log_2(|\E[\Delta G_\bell]|)$},
	every y tick label/.append style={font=\scriptsize},
	legend style ={ at={(1.03,1)}, 
        anchor=north west, draw=black, style={font=\scriptsize},
        fill=white,align=left},
   ymajorgrids,
	xmajorgrids
]

%
% ZERO
%

	\pgfplotstableread{data/G2F3/Ennn.txt} \data;
	\addplot [
		color=black,
		solid,
		line width=1pt,
		line cap=round,
		mark=*,
		mark size = 1.4,
		mark options={solid}
	]
	table[
		x expr=\thisrow{x},
		y expr=\thisrow{y}
	] {\data};
	%\addlegendentry{$(0,0,0)$};

%
% SINGLES
%

	\pgfplotstableread{data/G2F3/Einn.txt} \data;
	\addplot [
		color=ratecolor1,
		dashed,
		line width=1pt,
		line cap=round,
		mark=*,
		mark size = 1.4,
		mark options={solid,ratecolor1}
	]
	table[
		x expr=\thisrow{x},
		y expr=\thisrow{y}
	] {\data};
	%\addlegendentry{$(\ell,0,0)$};
	
	\pgfplotstableread{data/G2F3/Enin.txt} \data;
	\addplot [
		color=ratecolor1,
		dashed,
		line width=1pt,
		line cap=round,
		mark=square*,
		mark size = 1.4,
		mark options={solid,ratecolor1}
	]
	table[
		x expr=\thisrow{x},
		y expr=\thisrow{y}
	] {\data};
	%\addlegendentry{$(0,\ell,0)$};
	
	\pgfplotstableread{data/G2F3/Enni.txt} \data;
	\addplot [
		color=ratecolor1,
		dashed,
		line width=1pt,
		line cap=round,
		mark=diamond*,
		mark options={solid,ratecolor1}
	]
	table[
		x expr=\thisrow{x},
		y expr=\thisrow{y}
	] {\data};
	%\addlegendentry{$(0,0,\ell)$};
	
%
% DOUBLES
%

	\pgfplotstableread{data/G2F3/Eiin.txt} \data;
	\addplot [
		color=ratecolor2,
		dotted,
		line width=1pt,
		line cap=round,
		dash pattern=on 0pt off 2\pgflinewidth,
		mark=*,
		mark size = 1.4,
		mark options={solid,ratecolor2}
	]
	table[
		x expr=\thisrow{x},
		y expr=\thisrow{y}
	] {\data};
	%\addlegendentry{$(\ell,\ell,0)$};
	
	\pgfplotstableread{data/G2F3/Eini.txt} \data;
	\addplot [
		color=ratecolor2,
		dotted,
		line width=1pt,
		line cap=round,
		dash pattern=on 0pt off 2\pgflinewidth,
		mark=square*,
		mark size = 1.4,
		mark options={solid,ratecolor2}
	]
	table[
		x expr=\thisrow{x},
		y expr=\thisrow{y}
	] {\data};
	%\addlegendentry{$(\ell,0,\ell)$};
	
	\pgfplotstableread{data/G2F3/Enii.txt} \data;
	\addplot [
		color=ratecolor2,
		dotted,
		line width=1pt,
		line cap=round,
		dash pattern=on 0pt off 2\pgflinewidth,
		mark=diamond*,
		mark options={solid,ratecolor2}
	]
	table[
		x expr=\thisrow{x},
		y expr=\thisrow{y}
	] {\data};
	%\addlegendentry{$(0,\ell,\ell)$};

%
% TRIPPLE
%

	\pgfplotstableread{data/G2F3/Eiii.txt} \data;
	\addplot [
		color=ratecolor3,
		dash pattern=on 2pt off 1pt on 2pt off 2pt,
		line width=1pt,
		line cap=round,
		mark=*,
		mark size = 1.4,
		mark options={solid,ratecolor3}
	]
	table[
		x expr=\thisrow{x},
		y expr=\thisrow{y}
	] {\data};
	%\addlegendentry{$(\ell,\ell,\ell)$};
	
			\pgfplotsset{
    	after end axis/.code={
        	\node[above right] at (0.03,0.03){\texttt{G2F3}};
    	}
	}
	
\end{axis}
\end{tikzpicture}

%% file: figures/G2F3_analyse_V.tikz
\begin{tikzpicture}
\def\maxL{4}
\begin{axis}[
	height=\figureheight,
	width=\figurewidth,
	xmin=0,
	xmax=\maxL,
	xminorticks=true,
	xlabel={\small $\ell$},
	every x tick label/.append style={font=\scriptsize},
	ymin=-50,
	ymax=10,
	yminorticks=true,
	ylabel={\small $\log_2(\V[\Delta G_\bell])$},
	every y tick label/.append style={font=\scriptsize},
	legend style ={ at={(1.06,1)}, 
        anchor=north west, draw=black, style={font=\scriptsize},
        fill=white,align=left},
   ymajorgrids,
	xmajorgrids
]

%
% ZERO
%

	\pgfplotstableread{data/G2F3/Vnnn.txt} \data;
	\addplot [
		color=black,
		solid,
		line width=1pt,
		line cap=round,
		mark=*,
		mark size = 1.4,
		mark options={solid},
		forget plot
	]
	table[
		x expr=\thisrow{x},
		y expr=\thisrow{y}
	] {\data};
	%\addlegendentry{$(0,0,0)$};

%
% SINGLES
%

	\pgfplotstableread{data/G2F3/Vinn.txt} \data;
	\addplot [
		color=ratecolor1,
		dashed,
		line width=1pt,
		line cap=round,
		mark=*,
		mark size = 1.4,
		mark options={solid,ratecolor1}
	]
	table[
		x expr=\thisrow{x},
		y expr=\thisrow{y}
	] {\data};
	%\addlegendentry{$(\ell,0,0)$};
	
	\pgfplotstableread{data/G2F3/Vnin.txt} \data;
	\addplot [
		color=ratecolor1,
		dashed,
		line width=1pt,
		line cap=round,
		mark=square*,
		mark size = 1.4,
		mark options={solid,ratecolor1}
	]
	table[
		x expr=\thisrow{x},
		y expr=\thisrow{y}
	] {\data};
	%\addlegendentry{$(0,\ell,0)$};
	
	\pgfplotstableread{data/G2F3/Vnni.txt} \data;
	\addplot [
		color=ratecolor1,
		dashed,
		line width=1pt,
		line cap=round,
		mark=diamond*,
		mark options={solid,ratecolor1}
	]
	table[
		x expr=\thisrow{x},
		y expr=\thisrow{y}
	] {\data};
	%\addlegendentry{$(0,0,\ell)$};
	
%
% DOUBLES
%

	\pgfplotstableread{data/G2F3/Viin.txt} \data;
	\addplot [
		color=ratecolor2,
		dotted,
		line width=1pt,
		line cap=round,
		dash pattern=on 0pt off 2\pgflinewidth,
		mark=*,
		mark size = 1.4,
		mark options={solid,ratecolor2}
	]
	table[
		x expr=\thisrow{x},
		y expr=\thisrow{y}
	] {\data};
	%\addlegendentry{$(\ell,\ell,0)$};
	
	\pgfplotstableread{data/G2F3/Vini.txt} \data;
	\addplot [
		color=ratecolor2,
		dotted,
		line width=1pt,
		line cap=round,
		dash pattern=on 0pt off 2\pgflinewidth,
		mark=square*,
		mark size = 1.4,
		mark options={solid,ratecolor2}
	]
	table[
		x expr=\thisrow{x},
		y expr=\thisrow{y}
	] {\data};
	%\addlegendentry{$(\ell,0,\ell)$};
	
	\pgfplotstableread{data/G2F3/Vnii.txt} \data;
	\addplot [
		color=ratecolor2,
		dotted,
		line width=1pt,
		line cap=round,
		dash pattern=on 0pt off 2\pgflinewidth,
		mark=diamond*,
		mark options={solid,ratecolor2}
	]
	table[
		x expr=\thisrow{x},
		y expr=\thisrow{y}
	] {\data};
	%\addlegendentry{$(0,\ell,\ell)$};

%
% TRIPPLE
%

	\pgfplotstableread{data/G2F3/Viii.txt} \data;
	\addplot [
		color=ratecolor3,
		dash pattern=on 2pt off 1pt on 2pt off 2pt,
		line width=1pt,
		line cap=round,
		mark=*,
		mark size = 1.4,
		mark options={solid,ratecolor3}
	]
	table[
		x expr=\thisrow{x},
		y expr=\thisrow{y}
	] {\data};
	%\addlegendentry{$(\ell,\ell,\ell)$};
	
			\pgfplotsset{
    	after end axis/.code={
        	\node[above right] at (0.03,0.03){\texttt{G2F3}};
    	}
	}
	
\end{axis}
\end{tikzpicture}

%% file: figures/legend_times_G2.tikz
\begin{tikzpicture}
    \begin{customlegend}[
    	legend columns=2,
		legend style={column sep=1ex,font=\scriptsize},
		legend cell align=left,
		legend entries={MLMC\\MIMC, TD\\MLQMC\\MIQMC, TD\\}
		]
		\addlegendimage{ratecolor1,line cap=round,solid,line width=1pt,mark=*,mark size=1.4pt,mark options={solid}}
		%\addlegendimage{ratecolor3,line cap=round,solid,line width=1pt,mark=square*,mark size=1.4pt,mark options={solid}}
		\addlegendimage{ratecolor2,line cap=round,solid,line width=1pt,mark=diamond*,mark options={solid}}
   	\addlegendimage{ratecolor1,dashed,line cap=round,line width=1pt,mark=*,mark size=1.4pt,mark options={solid,fill=white}}
   	%\addlegendimage{ratecolor3,dashed,line cap=round,line width=1pt,mark=square*,mark size=1.4pt,mark options={solid,fill=white}}
   	\addlegendimage{ratecolor2,dashed,line cap=round,line width=1pt,mark=diamond*,mark options={solid,fill=white}}
   \end{customlegend}
\end{tikzpicture}

%% file: figures/cost_G2F1.tikz
	\begin{tikzpicture}[trim axis left, trim axis right]
	\begin{axis}[
		width=\figurewidth,
		height=\figureheight,
		xmode=log,
		xmin=1e-03,
		xmax=0.1,
		xminorticks=true,
		xlabel={\small requested tolerance $\epsilon$},
		xmajorgrids,
		xminorgrids,
		every x tick label/.append style={font=\scriptsize},
		ymode=log,
		ymin=10000, 
		ymax=100000000,
		yminorticks=true,
		ylabel={\small standard cost},
		ytick={100,1000,10000,100000,1000000,10000000,100000000,1000000000},
		ymajorgrids,
		yminorgrids,
		every y tick label/.append style={font=\scriptsize},
	]
	
	\addplot [
		color=ratecolor1,
		solid,
		line width=1pt,
		mark=*,
		mark size = 1.4,
		mark options={solid,ratecolor1}
	]
 	table[
		x = x,
		y = y
	] {data/G2F1MLMC/cost_formatted.txt};
	%\addlegendentry{MLMC};

	\addplot [
		color=ratecolor1,
		dashed,
		line cap=round,
		line width=1pt,
		mark=*,
		mark size = 1.4,
		mark options={solid,fill=white}
	]
 	table[
		x = x,
		y = y
	] {data/G2F1MLQM/cost_formatted.txt};
	%\addlegendentry{MLQMC};
	
%	\addplot [
%		color=blue,
%		solid,
%		mark=square*,
%		mark size=0.75pt,
%		line width = 0.5pt,
%		y dir = both,
%		y explicit,
%		error bar style={line width=0.5pt},
%    	error mark options={line width=0.5pt,rotate=90}
%	]
% 	table[
%		x = x,
%		y = y,
%		y error plus = z,
%		y error minus = u
%	] {data/G2F1_FTMC.txt};
%	%\addlegendentry{MIMC, FT};
%	
%	\addplot [
%		color=blue,
%		dashed,
%		mark=square*,
%		mark size=0.75pt,
%		line width = 0.5pt,
%		y dir = both,
%		y explicit,
%		error bar style={line width=0.5pt},
%    	error mark options={line width=0.5pt,rotate=90}
%	]
% 	table[
%		x = x,
%		y = y,
%		y error plus = z,
%		y error minus = u
%	] {data/G2F1_FTQM.txt};
%	%\addlegendentry{MIQMC, FT};

	\addplot [
		color=ratecolor2,
		solid,
		mark=diamond*,
		line width=1pt,
		mark options={solid,fill=white}
	]
 	table[
		x = x,
		y = y
	] {data/G2F1TDMC/cost_formatted.txt};
	%\addlegendentry{MIMC, TD};

	\addplot [
		color=ratecolor2,
		dashed,
		line cap=round,
		mark=diamond*,
		line width=1pt,
		mark options={solid,ratecolor2}
	]
 	table[
		x = x,
		y = y
	] {data/G2F1TDQM/cost_formatted.txt};
	%\addlegendentry{MIQMC, TD};
	
	\LogLogSlopeTriangle{0.9}{0.2}{0.75}{2}{black};
	\LogLogSlopeTriangle{0.9}{0.2}{0.55}{1}{black};
	
	\pgfplotsset{
    	after end axis/.code={
        	\node[above right] at (axis cs:1e-3,10000){\texttt{G2F1}};
    	}
	}
	
	\end{axis}
\end{tikzpicture}

%% file: figures/times_G2F1.tikz
	\begin{tikzpicture}[trim axis left, trim axis right]
	\begin{axis}[
		width=\figurewidth,
		height=\figureheight,
		xmode=log,
		xmin=1e-03,
		xmax=0.1,
		xminorticks=true,
		xlabel={\small requested tolerance $\epsilon$},
		xmajorgrids,
		xminorgrids,
		every x tick label/.append style={font=\scriptsize},
		ymode=log,
		ymin=10,
		ymax=100000,
		yminorticks=true,
		ylabel={\small runtime [seconds]},
		ytick={10,100,1000,10000,100000,1000000},
		ymajorgrids,
		yminorgrids,
		every y tick label/.append style={font=\scriptsize},
	]
		
	\addplot [
		color=ratecolor1,
		solid,
		line width=1pt,
		mark=*,
		mark size = 1.4,
		mark options={solid,ratecolor1}
	]
 	table[
		x = x,
		y = y
	] {data/G2F1MLMC/times_formatted.txt};
	%\addlegendentry{MLMC};

	\addplot [
		color=ratecolor1,
		dashed,
		line cap=round,
		line width=1pt,
		mark=*,
		mark size = 1.4,
		mark options={solid,fill=white}
	]
 	table[
		x = x,
		y = y
	] {data/G2F1MLQM/times_formatted.txt};
	%\addlegendentry{MLQMC};
	
%	\addplot [
%		color=blue,
%		solid,
%		mark=square*,
%		mark size=0.75pt,
%		line width = 0.5pt,
%		y dir = both,
%		y explicit,
%		error bar style={line width=0.5pt},
%    	error mark options={line width=0.5pt,rotate=90}
%	]
% 	table[
%		x = x,
%		y = y,
%		y error plus = z,
%		y error minus = u
%	] {data/G2F1_FTMC.txt};
%	%\addlegendentry{MIMC, FT};
%	
%	\addplot [
%		color=blue,
%		dashed,
%		mark=square*,
%		mark size=0.75pt,
%		line width = 0.5pt,
%		y dir = both,
%		y explicit,
%		error bar style={line width=0.5pt},
%    	error mark options={line width=0.5pt,rotate=90}
%	]
% 	table[
%		x = x,
%		y = y,
%		y error plus = z,
%		y error minus = u
%	] {data/G2F1_FTQM.txt};
%	%\addlegendentry{MIQMC, FT};

	\addplot [
		color=ratecolor2,
		solid,
		mark=diamond*,
		line width=1pt,
		mark options={solid,ratecolor2}
	]
 	table[
		x = x,
		y = y
	] {data/G2F1TDMC/times_formatted.txt};
	%\addlegendentry{MIMC, TD};

	\addplot [
		color=ratecolor2,
		dashed,
		line cap=round,
		mark=diamond*,
		line width=1pt,
		mark options={solid,fill=white}
	]
 	table[
		x = x,
		y = y
	] {data/G2F1TDQM/times_formatted.txt};
	%\addlegendentry{MIQMC, TD};
	
	\LogLogSlopeTriangle{0.9}{0.2}{0.75}{2}{black};
	\LogLogSlopeTriangle{0.9}{0.2}{0.55}{1}{black};

		\pgfplotsset{
    	after end axis/.code={
        	\node[above right] at (axis cs:1e-3,10){\texttt{G2F1}};
    	}
	}
	
	\end{axis}
\end{tikzpicture}

%% file: figures/cost_G2F2.tikz
	\begin{tikzpicture}[trim axis left, trim axis right]
	\begin{axis}[
		width=\figurewidth,
		height=\figureheight,
		xmode=log,
		xmin=1e-03,
		xmax=0.1,
		xminorticks=true,
		xlabel={\small requested tolerance $\epsilon$},
		xmajorgrids,
		xminorgrids,
		every x tick label/.append style={font=\scriptsize},
		ymode=log,
		ymin=10000, 
		ymax=100000000,
		yminorticks=true,
		ylabel={\small standard cost},
		ytick={100,1000,10000,100000,1000000,10000000,100000000,1000000000},
		ymajorgrids,
		yminorgrids,
		every y tick label/.append style={font=\scriptsize},
	]

	\addplot [
		color=ratecolor1,
		solid,
		line width=1pt,
		mark=*,
		mark size = 1.4,
		mark options={solid,ratecolor1}
	]
 	table[
		x = x,
		y = y
	] {data/G2F2MLMC/cost_formatted.txt};
	%\addlegendentry{MLMC};

	\addplot [
		color=ratecolor1,
		dashed,
		line cap=round,
		line width=1pt,
		mark=*,
		mark size = 1.4,
		mark options={solid,fill=white}
	]
 	table[
		x = x,
		y = y
	] {data/G2F2MLQM/cost_formatted.txt};
	%\addlegendentry{MLQMC};
	
%	\addplot [
%		color=blue,
%		solid,
%		mark=square*,
%		mark size=0.75pt,
%		line width = 0.5pt,
%		y dir = both,
%		y explicit,
%		error bar style={line width=0.5pt},
%    	error mark options={line width=0.5pt,rotate=90}
%	]
% 	table[
%		x = x,
%		y = y,
%		y error plus = z,
%		y error minus = u
%	] {data/G2F2_FTMC.txt};
%	%\addlegendentry{MIMC, FT};
%	
%	\addplot [
%		color=blue,
%		dashed,
%		mark=square*,
%		mark size=0.75pt,
%		line width = 0.5pt,
%		y dir = both,
%		y explicit,
%		error bar style={line width=0.5pt},
%    	error mark options={line width=0.5pt,rotate=90}
%	]
% 	table[
%		x = x,
%		y = y,
%		y error plus = z,
%		y error minus = u
%	] {data/G2F2_FTQM.txt};
%	%\addlegendentry{MIQMC, FT};

	\addplot [
		color=ratecolor2,
		solid,
		mark=diamond*,
		line width=1pt,
		mark options={solid,ratecolor2}
	]
 	table[
		x = x,
		y = y
	] {data/G2F2TDMC/cost_formatted.txt};
	%\addlegendentry{MIMC, TD};

	\addplot [
		color=ratecolor2,
		dashed,
		line cap=round,
		mark=diamond*,
		line width=1pt,
		mark options={solid,fill=white}
	]
 	table[
		x = x,
		y = y
	] {data/G2F2TDQM/cost_formatted.txt};
	%\addlegendentry{MIQMC, TD};
	
	\LogLogSlopeTriangle{0.9}{0.2}{0.7}{2}{black};
	\LogLogSlopeTriangle{0.9}{0.2}{0.5}{1}{black};
	
	\pgfplotsset{
    	after end axis/.code={
        	\node[above right] at (axis cs:1e-3,10000){\texttt{G2F2}};
    	}
	}
	
	\end{axis}
\end{tikzpicture}

%% file: figures/times_G2F2.tikz
	\begin{tikzpicture}[trim axis left, trim axis right]
	\begin{axis}[
		width=\figurewidth,
		height=\figureheight,
		xmode=log,
		xmin=1e-03,
		xmax=0.1,
		xminorticks=true,
		xlabel={\small requested tolerance $\epsilon$},
		xmajorgrids,
		xminorgrids,
		every x tick label/.append style={font=\scriptsize},
		ymode=log,
		ymin=10,
		ymax=100000,
		yminorticks=true,
		ylabel={\small runtime [seconds]},
		ytick={10,100,1000,10000,100000,1000000},
		ymajorgrids,
		yminorgrids,
		every y tick label/.append style={font=\scriptsize},
	]
	
	\addplot [
		color=ratecolor1,
		solid,
		line width=1pt,
		mark=*,
		mark size = 1.4,
		mark options={solid,ratecolor1}
	]
 	table[
		x = x,
		y = y
	] {data/G2F2MLMC/times_formatted.txt};
	%\addlegendentry{MLMC};

	\addplot [
		color=ratecolor1,
		dashed,
		line cap=round,
		line width=1pt,
		mark=*,
		mark size = 1.4,
		mark options={solid,fill=white}
	]
 	table[
		x = x,
		y = y
	] {data/G2F2MLQM/times_formatted.txt};
	%\addlegendentry{MLQMC};
	
%	\addplot [
%		color=blue,
%		solid,
%		mark=square*,
%		mark size=0.75pt,
%		line width = 0.5pt,
%		y dir = both,
%		y explicit,
%		error bar style={line width=0.5pt},
%    	error mark options={line width=0.5pt,rotate=90}
%	]
% 	table[
%		x = x,
%		y = y,
%		y error plus = z,
%		y error minus = u
%	] {data/G2F2_FTMC.txt};
%	%\addlegendentry{MIMC, FT};
%	
%	\addplot [
%		color=blue,
%		dashed,
%		mark=square*,
%		mark size=0.75pt,
%		line width = 0.5pt,
%		y dir = both,
%		y explicit,
%		error bar style={line width=0.5pt},
%    	error mark options={line width=0.5pt,rotate=90}
%	]
% 	table[
%		x = x,
%		y = y,
%		y error plus = z,
%		y error minus = u
%	] {data/G2F2_FTQM.txt};
%	%\addlegendentry{MIQMC, FT};

	\addplot [
		color=ratecolor2,
		solid,
		mark=diamond*,
		line width=1pt,
		mark options={solid,ratecolor2}
	]
 	table[
		x = x,
		y = y
	] {data/G2F2TDMC/times_formatted.txt};
	%\addlegendentry{MIMC};

	\addplot [
		color=ratecolor2,
		dashed,
		line cap=round,
		mark=diamond*,
		line width=1pt,
		mark options={solid,fill=white}
	]
 	table[
		x = x,
		y = y
	] {data/G2F2TDQM/times_formatted.txt};
	%\addlegendentry{MIQMC};
	
	\LogLogSlopeTriangle{0.9}{0.2}{0.7}{2}{black};
	\LogLogSlopeTriangle{0.9}{0.2}{0.5}{1}{black};

		\pgfplotsset{
    	after end axis/.code={
        	\node[above right] at (axis cs:1e-3,10){\texttt{G2F2}};
    	}
	}

	\end{axis}
\end{tikzpicture}

%% file: figures/cost_G2F3.tikz
	\begin{tikzpicture}[trim axis left, trim axis right]
	\begin{axis}[
		width=\figurewidth,
		height=\figureheight,
		xmode=log,
		xmin=1e-03,
		xmax=0.1,
		xminorticks=true,
		xlabel={\small requested tolerance $\epsilon$},
		xmajorgrids,
		xminorgrids,
		every x tick label/.append style={font=\scriptsize},
		ymode=log,
		ymin=100, 
		ymax=1000000,
		yminorticks=true,
		ylabel={\small standard cost},
		ytick={100,1000,10000,100000,1000000,10000000,100000000,1000000000},
		ymajorgrids,
		yminorgrids,
		every y tick label/.append style={font=\scriptsize},
	]

	\addplot [
		color=ratecolor1,
		solid,
		line width=1pt,
		mark=*,
		mark size = 1.4,
		mark options={solid,ratecolor1}
	]
 	table[
		x = x,
		y = y
	] {data/G2F3MLMC/cost_formatted.txt};
	%\addlegendentry{MLMC};

	\addplot [
		color=ratecolor1,
		dashed,
		line cap=round,
		line width=1pt,
		mark=*,
		mark size = 1.4,
		mark options={solid,fill=white}
	]
 	table[
		x = x,
		y = y
	] {data/G2F3MLQM/cost_formatted.txt};
	%\addlegendentry{MLQMC};
	
%	\addplot [
%		color=blue,
%		solid,
%		mark=square*,
%		mark size=0.75pt,
%		line width = 0.5pt,
%		y dir = both,
%		y explicit,
%		error bar style={line width=0.5pt},
%    	error mark options={line width=0.5pt,rotate=90}
%	]
% 	table[
%		x = x,
%		y = y,
%		y error plus = z,
%		y error minus = u
%	] {data/G2F3_FTMC.txt};
%	%\addlegendentry{MIMC, FT};
%	
%	\addplot [
%		color=blue,
%		dashed,
%		mark=square*,
%		mark size=0.75pt,
%		line width = 0.5pt,
%		y dir = both,
%		y explicit,
%		error bar style={line width=0.5pt},
%    	error mark options={line width=0.5pt,rotate=90}
%	]
% 	table[
%		x = x,
%		y = y,
%		y error plus = z,
%		y error minus = u
%	] {data/G2F3_FTQM.txt};
%	%\addlegendentry{MIQMC, FT};

	\addplot [
		color=ratecolor2,
		solid,
		mark=diamond*,
		line width=1pt,
		mark options={solid,ratecolor2}
	]
 	table[
		x = x,
		y = y
	] {data/G2F3TDMC/cost_formatted.txt};
	%\addlegendentry{MIMC, TD};

	\addplot [
		color=ratecolor2,
		dashed,
		line cap=round,
		mark=diamond*,
		line width=1pt,
		mark options={solid,fill=white}
	]
 	table[
		x = x,
		y = y
	] {data/G2F3TDQM/cost_formatted.txt};
	%\addlegendentry{MIQMC, TD};
	
	\LogLogSlopeTriangle{0.3}{0.2}{0.2}{2}{black};
	\LogLogSlopeTriangle{0.3}{0.2}{0.5}{3}{black};
	
	\pgfplotsset{
    	after end axis/.code={
        	\node[above right] at (axis cs:1e-3,100){\texttt{G2F3}};
    	}
	}

	\end{axis}
\end{tikzpicture}

%% file: figures/times_G2F3.tikz
	\begin{tikzpicture}[trim axis left, trim axis right]
	\begin{axis}[
		width=\figurewidth,
		height=\figureheight,
		xmode=log,
		xmin=1e-03,
		xmax=0.1,
		xminorticks=true,
		xlabel={\small requested tolerance $\epsilon$},
		xmajorgrids,
		xminorgrids,
		every x tick label/.append style={font=\scriptsize},
		ymode=log,
		ymin=10,
		ymax=100000,
		yminorticks=true,
		ylabel={\small runtime [seconds]},
		ytick={10,100,1000,10000,100000,1000000},
		ymajorgrids,
		yminorgrids,
		every y tick label/.append style={font=\scriptsize},
	]

	\addplot [
		color=ratecolor1,
		solid,
		line width=1pt,
		mark=*,
		mark size = 1.4,
		mark options={solid,ratecolor1}
	]
 	table[
		x = x,
		y = y
	] {data/G2F3MLMC/times_formatted.txt};
	%\addlegendentry{MLMC};

	\addplot [
		color=ratecolor1,
		dashed,
		line cap=round,
		line width=1pt,
		mark=*,
		mark size = 1.4,
		mark options={solid,fill=white}
	]
 	table[
		x = x,
		y = y
	] {data/G2F3MLQM/times_formatted.txt};
	%\addlegendentry{MLQMC};
	
%	\addplot [
%		color=blue,
%		solid,
%		mark=square*,
%		mark size=0.75pt,
%		line width = 0.5pt,
%		y dir = both,
%		y explicit,
%		error bar style={line width=0.5pt},
%    	error mark options={line width=0.5pt,rotate=90}
%	]
% 	table[
%		x = x,
%		y = y,
%		y error plus = z,
%		y error minus = u
%	] {data/G2F3_FTMC.txt};
%	%\addlegendentry{MIMC, FT};
%	
%	\addplot [
%		color=blue,
%		dashed,
%		mark=square*,
%		mark size=0.75pt,
%		line width = 0.5pt,
%		y dir = both,
%		y explicit,
%		error bar style={line width=0.5pt},
%    	error mark options={line width=0.5pt,rotate=90}
%	]
% 	table[
%		x = x,
%		y = y,
%		y error plus = z,
%		y error minus = u
%	] {data/G2F3_FTQM.txt};
%	%\addlegendentry{MIQMC, FT};

	\addplot [
		color=ratecolor2,
		solid,
		mark=diamond*,
		line width=1pt,
		mark options={solid,ratecolor2}
	]
 	table[
		x = x,
		y = y
	] {data/G2F3TDMC/times_formatted.txt};
	%\addlegendentry{MIMC};

	\addplot [
		color=ratecolor2,
		dashed,
		line cap=round,
		mark=diamond*,
		line width=1pt,
		mark options={solid,fill=white}
	]
 	table[
		x = x,
		y = y
	] {data/G2F3TDQM/times_formatted.txt};
	%\addlegendentry{MIQMC};
	
	%\LogLogSlopeTriangle{0.9}{0.2}{0.7}{2}{black};
	%\LogLogSlopeTriangle{0.9}{0.2}{0.5}{1}{black};
	
	\LogLogSlopeTriangle{0.3}{0.2}{0.2}{2}{black};
	\LogLogSlopeTriangle{0.3}{0.2}{0.5}{3}{black};

		\pgfplotsset{
    	after end axis/.code={
        	\node[above right] at (axis cs:1e-3,10){\texttt{G2F3}};
    	}
	}
	
	\end{axis}
\end{tikzpicture}

%% file: MIQMC.bbl
\begin{thebibliography}{10}

\bibitem{bachmayr2017representations}
M.~Bachmayr, A.~Cohen, and G.~Migliorati.
\newblock {Representations of Gaussian Random Fields and Approximation of
  Elliptic PDEs with Lognormal Coefficients}.
\newblock {\em Journal of Fourier Analysis and Applications}, 1:1--29, 2017.

\bibitem{barth2011multi}
A.~Barth, C.~Schwab, and N.~Zollinger.
\newblock {Multi-level Monte Carlo Finite Element Method for Elliptic PDEs with
  Stochastic Coefficients}.
\newblock {\em Numerische Mathematik}, 119(1):123--161, 2011.

\bibitem{boyle2009hsl}
J.~Boyle, M.~Mihajlovi\'c, and J.~Scott.
\newblock {HSL{\_}MI20: An Efficient AMG Preconditioner for Finite Element
  Problems in 3D}.
\newblock {\em International Journal for Numerical Methods in Engineering},
  2009.

\bibitem{bungartz2004sparse}
H.-J. Bungartz and M.~Griebel.
\newblock {Sparse Grids}.
\newblock {\em Acta Numerica}, 13:147--269, 2004.

\bibitem{chkifa2014high}
A.~Chkifa, A.~Cohen, and C.~Schwab.
\newblock {High-dimensional Adaptive Sparse Polynomial Interpolation and
  Applications to Parametric PDEs}.
\newblock {\em Foundations of Computational Mathematics}, 4(14):601--633, 2014.

\bibitem{cliffe2011multilevel}
K.~A. Cliffe, M.~B. Giles, R.~Scheichl, and A.~L. Teckentrup.
\newblock {Multilevel Monte Carlo Methods and Applications to Elliptic PDEs
  with Random Coefficients}.
\newblock {\em Computing and Visualization in Science}, 14(1):3--15, 2011.

\bibitem{collier2014continuation}
N.~Collier, A.-L. Haji-Ali, F.~Nobile, E.~Schwerin, and R.~Tempone.
\newblock {A continuation multilevel Monte Carlo algorithm}.
\newblock {\em BIT Numerical Mathematics}, 55(2):399--432, 2014.

\bibitem{dick2013higher}
J.~Dick, F.~Y. Kuo, Q.~T.~L. Gia, D.~Nuyens, and C.~Schwab.
\newblock {Higher order QMC Petrov--Galerkin discretization for affine
  parametric operator equations with random field inputs}.
\newblock {\em SIAM Journal on Numerical Analysis}, 52(6):2676--2702, 2014.

\bibitem{dick2013high}
J.~Dick, F.~Y. Kuo, and I.~H. Sloan.
\newblock {High-Dimensional Integration: The Quasi-Monte Carlo way}.
\newblock {\em Acta Numerica}, 22:133--288, 2013.

\bibitem{dick2010digital}
J.~Dick and F.~Pillichshammer.
\newblock {\em {Digital Nets and Sequences: Discrepancy Theory and Quasi--Monte
  Carlo Integration}}.
\newblock Cambridge University Press, 2010.

\bibitem{feischl2017fast}
M.~Feischl, F.~Y. Kuo, and I.~H. Sloan.
\newblock {Fast Random Field Generation with $H$-Matrices}.
\newblock in preparation, 2017.

\bibitem{gerstner2003dimension}
T.~Gerstner and M.~Griebel.
\newblock {Dimension--Adaptive Tensor--Product Quadrature}.
\newblock {\em Computing}, 71(1):65--87, 2003.

\bibitem{ghanem1991stochastic}
R.~G. Ghanem and P.~D. Spanos.
\newblock {\em {Stochastic Finite Elements: A Spectral Approach}}.
\newblock Springer New York, 1 edition, 1991.

\bibitem{giles2008multilevel}
M.~B. Giles.
\newblock {Multilevel Monte Carlo Path Simulation}.
\newblock {\em Operations Research}, 56(3):607--617, 2008.

\bibitem{giles2009multilevel}
M.~B. Giles.
\newblock {Multilevel Quasi-Monte Carlo Path Simulation}.
\newblock {\em Advanced Financial Modelling, Radon Series on Computational and
  Applied Mathematics}, pages 165--181, 2009.

\bibitem{giles2015multilevel}
M.~B. Giles.
\newblock {Multilevel Monte Carlo Methods}.
\newblock {\em Acta Numerica}, 24:259--328, 2015.

\bibitem{graham2014quasi}
I.~G. Graham, F.~Y. Kuo, J.~A. Nichols, R.~Scheichl, C.~Schwab, and I.~H.
  Sloan.
\newblock {Quasi-Monte Carlo Finite Element Methods for Elliptic PDEs with
  Lognormal Random Coefficients}.
\newblock {\em Numerische Mathematik}, 131(2):329--368, 2014.

\bibitem{graham2011quasi}
I.~G. Graham, F.~Y. Kuo, D.~Nuyens, R.~Scheichl, and I.~H. Sloan.
\newblock {Quasi-Monte Carlo Methods for Elliptic PDEs with Random Coefficients
  and Applications}.
\newblock {\em Journal of Computational Physics}, 230(10):3668--3694, may 2011.

\bibitem{griebel1992combination}
M.~Griebel, M.~Schneider, and C.~Zenger.
\newblock {A Combination Technique for the Solution of Sparse Grid Problems}.
\newblock In P.~de~Groen and R.~Beauwens, editors, {\em Iterative Methods in
  Linear Algebra}, pages 263--281. Elsevier, Amsterdam, 1992.

\bibitem{haji2016misc}
A.-L. Haji-Ali, F.~Nobile, L.~Tamellini, and R.~Tempone.
\newblock {Multi-Index Stochastic Collocation for Random PDEs}.
\newblock {\em Computer Methods in Applied Mechanics and Engineering},
  306:95--122, 2016.

\bibitem{haji2016multi}
A.-L. Haji-Ali, F.~Nobile, and R.~Tempone.
\newblock {Multi-Index Monte Carlo: When Sparsity Meets Sampling}.
\newblock {\em Numerische Mathematik}, 132(4):767--806, apr 2016.

\bibitem{heinrich2001multilevel}
S.~Heinrich.
\newblock {Multilevel Monte Carlo Methods}.
\newblock In I.~Lirkov, S.~D. Margenov, and J.~Wasniewski, editors, {\em
  Large-Scale Scientific Computing}, pages 58--67. Springer Verlag, Heidelberg,
  2001.

\bibitem{kuo2016application}
F.~Y. Kuo and D.~Nuyens.
\newblock Application of quasi-monte carlo methods to elliptic pdes with random
  diffusion coefficients: A survey of analysis and implementation.
\newblock {\em Foundations of Computational Mathematics}, 16(6):1631--1696,
  2016.

\bibitem{kuo2015multilevel}
F.~Y. Kuo, R.~Scheichl, C.~Schwab, I.~H. Sloan, and E.~Ullmann.
\newblock {Multilevel Quasi-Monte Carlo Methods for Lognormal Diffusion
  Problems}.
\newblock {\em Mathematics of Computation, in press}, 2016.

\bibitem{kuo2005lifting}
F.~Y. Kuo and I.~H. Sloan.
\newblock {Lifting the Curse of Dimensionality}.
\newblock {\em Notices of the AMS}, 52(11):1320--1328, 2005.

\bibitem{nichols2014fast}
J.~A. Nichols and F.~Y. Kuo.
\newblock {Fast CBC Construction of Randomly Shifted Lattice Rules Achieving
  Convergence for Unbounded Integrands over $\mathbb{R}^s$ in Weighted Spaces
  with POD Weights}.
\newblock {\em Journal of Complexity}, 30(4):444--468, aug 2014.

\bibitem{nobile2015convergence}
F.~Nobile, L.~Tamellini, and R.~Tempone.
\newblock {Convergence of Quasi-Optimal Sparse-grid Approximation of
  Hilbert-space-valued Functions: Application to Random Elliptic PDEs}.
\newblock {\em Numerische Mathematik}, pages 1--46, 2015.

\bibitem{nuyens2010magic}
D.~Nuyens.
\newblock {The magic point shop of QMC point generators and generating
  vectors}, 2010.

\bibitem{nuyens2006fast}
D.~Nuyens and R.~Cools.
\newblock {Fast Algorithms for Component-by-component Construction of Rank-1
  Lattice Rules in Shift-Invariant Reproducing Kernel Hilbert Spaces}.
\newblock {\em Mathematics of Computation}, 75(254):903--920, 2006.

\bibitem{QMC4PDE}
D.~Nuyens and F.~Y. Kuo.
\newblock {QMC4PDE: A practical guide to the software for constructing point
  sets and point generator code}, 2016.

\bibitem{paskov1995faster}
S.~H. Paskov and J.~F. Traub.
\newblock {Faster Valuation of Financial Derivatives}.
\newblock {\em The Journal of Portfolio Management}, 22(1):113--123, 1995.

\bibitem{pauli2015multilevel}
S.~Pauli, R.~N. Gantner, P.~Arbenz, and A.~Adelmann.
\newblock {Multilevel Monte Carlo for the Feynman--Kac Formula for the Laplace
  Equation}.
\newblock {\em BIT Numerical Mathematics}, 55(4):1125--1143, 2015.

\bibitem{robbe2016thesis}
P.~Robbe, D.~Nuyens, and S.~Vandewalle.
\newblock {\em {A Practical Multilevel Quasi-Monte Carlo Method for Elliptic
  PDEs with Random Coefficients}}.
\newblock masters thesis, ``Een Parallelle Multilevel Monte-Carlo-methode voor
  de Simulatie van Stochastische Parti\"ele {Differentiaalvergelijkingen}'' by
  P. Robbe, KU Leuven, 2016.

\bibitem{robbe2017dimension}
P.~Robbe, D.~Nuyens, and S.~Vandewalle.
\newblock {A Dimension-Adaptive Multi-Index Monte Carlo Method Applied to a
  Heat Exchanger}.
\newblock In {\em 12th International Conference on Monte Carlo and Quasi-Monte
  Carlo Methods in Scientific Computing}, submitted, 2017.

\bibitem{schwab2006karhunen}
C.~Schwab and R.~A. Todor.
\newblock {Karhunen--Lo\`{e}ve Approximation of Random Fields by Generalized
  Fast Multipole Methods}.
\newblock {\em Journal of Computational Physics}, 217(1):100--122, 2006.

\bibitem{smolyak1960interpolation}
S.~A. Smolyak.
\newblock {Interpolation and Quadrature Formulas for the Classes $W^a_s$ and
  $E^a_s$}.
\newblock In {\em Dokl. Akad. Nauk SSSR}, volume 131, pages 1028--1031. (In
  Russian, Engl. Transl.: Soviet Math. Dokl. 4, 240-243 (1963)), 1960.

\bibitem{wackerly2008mathematical}
D.~Wackerly, W.~Mendenhall, and R.~Scheaffer.
\newblock {\em Mathematical Statistics with Applications}.
\newblock Thomson Brooks/Cole, 2008.

\bibitem{xiu2002wiener}
D.~Xiu and G.~E. Karniadakis.
\newblock {The Wiener--Askey Polynomial Chaos for Stochastic Differential
  Equations}.
\newblock {\em SIAM Journal on Scientific Computing}, 24(2):619--644, 2002.

\end{thebibliography}
